\documentclass[trsc,nonblindrev]{informs3} % current default for manuscript submission
\OneAndAHalfSpacedXI % current default line spacing
\usepackage{amsfonts,epsfig,epstopdf}
\usepackage{xcolor}% http://ctan.org/pkg/xcolor
\usepackage{algpseudocode}% http://ctan.org/pkg/algorithmicx
\usepackage{algorithm}% http://ctan.org/pkg/algorithm
\PassOptionsToPackage{hyphens}{url}
\usepackage{url}
\usepackage{bm}
\usepackage{pifont}% http://ctan.org/pkg/pifont
\usepackage{changes}%[final]
\usepackage{float}
\usepackage{subfig}
\usepackage{lscape}
\usepackage{multirow}
\usepackage{booktabs}
\usepackage[para,online,flushleft]{threeparttable}

\usepackage{cleveref}
\usepackage{siunitx}
\usepackage{stfloats}
\usepackage{mathtools}
\usetikzlibrary{arrows}
\usepackage{tikz}
\usepackage{pgfplots}
\usepackage{natbib} % http://merkel.texture.rocks/Latex/natbib.php

\bibpunct[, ]{(}{)}{,}{a}{}{,}
\def\bibfont{\small}%

\TheoremsNumberedThrough     % Preferred (Theorem 1, Lemma 1, Theorem 2)
\EquationsNumberedThrough    % Default: (1), (2), ...
\begin{document}
%%%%%%%%%%%%%%%%

\TITLE{A branch-cut-and-price algorithm for the dial-a-ride problem with minimum disease-transmission risk}

\ARTICLEAUTHORS{%
\AUTHOR{Shuocheng Guo}
\AFF{Department of Civil, Construction and Environmental Engineering, The University of Alabama, Tuscaloosa, AL 35487,
\EMAIL{sguo18@ua.edu}}
\AUTHOR{Iman Dayarian}
\AFF{Culverhouse College of Business, The University of Alabama, Tuscaloosa, AL 35487, \EMAIL{idayarian@cba.ua.edu}}
\AUTHOR{Jian Li}
\AFF{School of Transportation Engineering, Tongji University, Shanghai, China, \EMAIL{jianli@tongji.edu.cn}}
\AUTHOR{Xinwu Qian}
\AFF{Department of Civil, Construction and Environmental Engineering, The University of Alabama, Tuscaloosa, AL 35487, \EMAIL{xinwu.qian@ua.edu}}
% Enter all authors
} % end of the block

\ABSTRACT{% 268 words< 300 words (limit).
This paper investigates a novel variant of the dial-a-ride problem (DARP), namely Risk-aware DARP (RDARP), which focuses on the bi-objective DARP by minimizing the total routing cost and the maximum exposed risk for onboard passengers. The RDARP requires that the exposed risk of each passenger is minimized along the trip while satisfying existing constraints of the DARP. We model the RDARP using a three-index arc-based formulation and reformulate it into an equivalent min-max trip-based formulation, which is solved optimally using a tailored Branch-Cut-and-Price (BCP) algorithm. The BCP algorithm adopts the Column Generation method by decomposing the problem into a master problem and a subproblem. The subproblem takes the form of an Elementary Shortest Path Problem with Resource Constraints and Minimized-Maximum Risk (ESPPRCMMR). To solve the ESPPRCMMR efficiently, we develop a new linear-time risk calibration algorithm and establish families of resource extension functions in compliance with the risk-related resources. Additionally, we extend the applicability of our setting and solution approach to address a generic bi-objective DARP, namely Equitable DARP (EDARP), which aims to distribute ride times among passengers equitably.
To demonstrate the effectiveness of our approach, we adopt a real-world paratransit trip dataset to generate the RDARP instances ranging from 17 to 55 heterogeneous passengers with 3 to 13 vehicles during the morning and afternoon periods.
Computational results show that our BCP algorithm can optimally solve 50 out of 54 real-world instances (up to 55 passengers) within a time limit of one hour.
}%
\maketitle

\fontsize{11}{17}\selectfont

\vspace{-1cm}
\section{Introduction}\label{sec:intro}
The \emph{dial-a-ride problem} (DARP) generalizes the \emph{pick-up and delivery problems with time windows} (PDPTWs) with the maximum ride time constraints~\citep{cordeau2007dial}, and has seen extensive real-world applications such as dial-a-ride services, on-demand mobility services, and paratransit services. For the DARP, each passenger must be picked up and dropped off within a specific time interval by the same vehicle. The focus of the DARP is primarily serving demand at a low cost~\citep{ho2018survey} by transporting more passengers at a time with limited resources.

Our study is motivated by the emerging real-world challenges posed by the COVID-19 pandemic on demand-responsive transit (DRT) operations. Specifically, the paratransit services, which provide essential mobility services to individuals eligible under the Americans with Disabilities Act, serve a vulnerable population such as people with disability, people with pre-existing health conditions, and senior citizens~\citep{nadtc2021}. 
The conventional DARP setting, which involves transporting more passengers at a time, facilitates the spread of infectious diseases and is unsuitable during pandemics, especially when considering the compound exposed risk towards the vulnerable population who require access to healthcare facilities~\citep{chen2020transportation,nie2022impact}. Consequently, there is an urgent need for DRT service providers such as paratransit, ridesharing, and flexible-route transit to develop risk-aware operation strategies that provide a safe travel environment for the vulnerable population in a cost-effective manner. 

In response to the emerging disease transmission concerns during the COVID-19 pandemic, we propose the risk-aware routing and scheduling framework for DRT services. Specifically, we introduce the Risk-aware DARP (RDARP), aiming to minimize both the total travel cost and the maximum individual exposed risk among onboard passengers (the reason for considering the \textit{min-max} exposed risk will be explained in Remark~\ref{remark:minimized_max_risk}). This framework has significant implications for mitigating the risk of ongoing waves of COVID-19 outbreak, seasonal communicable diseases (e.g., flu and respiratory infections), and future unknown diseases that may arise. Additionally, our model can be broadly extended to other practical applications that require minimal interactions among people or goods, such as food delivery service for fresh produce and frozen/hot food, the DRT service with co-rider preferences, patient transport service, and DARP with an equitable level of service (with the min-max objective of excessive ride time or detour rate).

This study proposes the mathematical model of the RDARP and a tailored algorithm to optimally solve it. The RDARP includes provisions for exposed risk among passengers with different levels of risk, subject to the propagation of disease risks along the route and the maximum cumulative exposed risk constraints, in addition to the conventional constraints in the DARP.

The main contributions of our study are summarized below: 
\begin{itemize}
    \item We propose the RDARP with a particular focus on the en-route exposed risk in addition to the total travel cost, which is formulated as a bi-objective optimization problem. Sets of risk-aware routing strategies on the Pareto front are provided to aid decision-making for DRT service.
    \item We solve the RDARP by means of the column generation method under a Branch-Cut-and-Price (BCP) framework. Our BCP model efficiently handles the non-linearity of the risk-related constraints using tailored valid inequalities and a labeling algorithm.
    \item We design a tailored labeling algorithm to provide an exact solution by leveraging the properties of the dynamic time window in the typical DARP and disease transmission among the onboard passengers. We consider the exposed risk among passengers as a function of onboard duration consisting of travel time, service time, and waiting time, and we introduce a family of risk-related resources with associated resource extension functions and dominance rules to accurately capture the minimum possible exposed risk. 
    \item We propose a model extension of our RDARP called Equitable DARP, which aims to distribute ride times equitably among the passengers with minor modifications to the RDARP. 
    \item We conduct extensive numerical experiments using the real-world paratransit trip log data in the central Alabama area. In particular, the trade-off between travel cost and maximum individual exposed risk is broadly discussed upon the results of sensitivity analyses and the Pareto fronts. 
\end{itemize}

The remainder of this paper is organized as follows. We review the related studies in Section~\ref{sec:lit}. Section~\ref{sec:prob_setting} formally introduces the arc-based formulations of the RDARP and the equivalent trip-based formulation, which is decomposed based on the Dantzig-Wolfe (DW) decomposition method. 
In Section~\ref{sec:Solution Approach}, we propose the labeling algorithm and a risk calibration algorithm to solve the RDARP.  Section~\ref{sec:branch_and_price} discusses the branch-and-price algorithm. The model extension is discussed in Section~\ref{sec:edarp}. Section~\ref{sec:results} introduces the data instances and presents the computational results and discussion. Finally, Section~\ref{sec:conclusion} concludes our study and proposes several future directions.

\section{Literature Review}\label{sec:lit}
Surveys on the DARP are provided by \cite{cordeau2007dial} and most recently, \cite{ho2018survey}.
In general, the DARP can be divided into static and dynamic modes based on whether the demand is pre-assigned~\citep{cordeau2007dial}. The static DARP is associated with pre-scheduled demand, e.g., dial-a-ride service and paratransit service~\citep{fu1999improving,gupta2010improving}. On the other hand, the dynamic DARP aims at providing high-quality service to immediate or short-term demand. Practical applications can be seen in the shared taxi service~\citep{hosni2014shared}, on-demand mobility system~\citep{sayarshad2018scalable}, and accessible taxi services~\citep{dikas2018scheduled}. The objectives of DARP mainly include minimizing the operating costs and/or user inconvenience~\citep[for a comprehensive survey of the objectives in DARP, see][]{paquette2009quality,molenbruch2017typology}, which can be handled by weighting or $\epsilon$-constraint methods~\citep{demir2014bi}. In addition, the min-max objectives have been investigated in the Vehicle Routing Problem with load balancing or workload equity~\citep{matl2018workload}, i.e., minimizing the maximum service completion time~\citep{perrier2008vehicle}, maximum routing cost~\citep{carlsson2009solving}.

The exact algorithms for the DARP mainly focus on the Branch-and-Cut (BC) and Branch-and-Price (BP) algorithms that are associated with the arc-based and trip-based formulations, respectively. In the arc-based formulation, each location is associated with the arrival time and load to ensure the feasibility regarding the time window, capacity, and ride time constraints. To tighten the linear programming (LP) relaxation of arc-based formulation, families of valid inequalities with separation heuristics were developed~\citep{cordeau2006branch,ropke2007models,ropke2009branch}, including subtour elimination constraints, rounded capacity constraints, and infeasible path elimination constraints (IPEC). Later, BC algorithms were applied to more practical cases in terms of different user types~\citep{parragh2011introducing}, multi-depot setting~\citep{braekers2014exact}, lunch-break for drivers~\citep{liu2015branch}, and driver consistency~\citep{paquette2013combining}. On the other hand, the trip-based formulation relies on the BP algorithms. The BP algorithm for PDPTW was first proposed by~\cite{dumas1991pickup} and extended by~\cite{savelsbergh1998drive}. \cite{ropke2009branch} incorporated a set of valid inequalities in the BC algorithm into the BP procedure, yielding a BCP algorithm. In addition, \cite{gschwind2015effective} proposed an exact BCP algorithm that handles the maximum ride time constraints. The feasibility check of the ride time constraint is conducted by inspecting the arrival time and the latest delivery time for each onboard passenger. \cite{luo2019two} applied the BCP algorithm to solve the DARP under sets of practical constraints using IPEC and Benders cuts.

Finally, it is worth noting that our concept of "risk while traveling" differs from existing studies. The term can be found in the cash transportation problem~\citep{allahyari2021novel}, where en-route risk relies only on simple propagation of the start-of-service time without contagion among onboard items and varying waiting times. On the contrary, our study considers potential waiting times and realistic risk-related constraints that are applicable to real-world operations. Additionally, accurate calculation of the minimal exposed risk depends on the exact start-of-service time, which has not been studied previously. 
Furthermore, although min-max DARP has been investigated in the equity issues~\citep{matl2018workload,lehuede2014multi}, no existing exact BCP algorithms are available to find the optimal solution nor the Pareto front.

\section{Problem Setting and Formulations}\label{sec:prob_setting}
We develop a risk-aware routing and scheduling framework, motivated by the emerging challenges posed by the COVID-19 pandemic on DRT operations. This section starts with the problem definition and the unique features of modeling the exposed risk among onboard passengers during travel. We then propose an arc-based formulation for the problem, which aids in understanding its structure and serves as a stepping stone for the trip-based formulation detailed next. Additionally, we underscore in Remark~\ref{remark:minimized_max_risk} the significance of the min-max exposed risk objective in the epidemic process over a physical contact network.

\subsection{Problem Definition}
Our RDARP is extended based on the static multi-vehicle DARP with additional consideration of exposed risk among onboard passengers. In our RDARP, passenger trips are made via reservations with known pick-up and drop-off locations, as well as the desired pick-up and/or drop-off time windows depending on outbound or inbound trips. Let $n$ be the number of {request locations}. The RDARP is defined on a complete directed graph denoted by $\mathcal{G} = (\mathcal{N},\mathcal{A)}$. The node set $\mathcal{N}$ includes four subsets: $\{0\}$ and $\{2n+1\}$ to represent the origin and final depots which may coincide; $\mathcal{P}= \{1,2,\ldots,n\}$ and $\mathcal{D}= \{n+1,n+2,\ldots,2n\}$ to denote the requests' pick-up and drop-off nodes, respectively. Hence, each request is associated with the pick-up node $i$ and drop-off node $n+i$. Each request may include one or multiple passengers, where $w_{i}$ indicates the number of passengers associated with location $i$, and therefore, $w_{i}=-w_{n+i},~\forall i\in\mathcal{P}$. For each node $i\in\mathcal{P}\cup\mathcal{D}$, the time window is represented by $\left[a_{i}, b_{i}\right]$, denoting the earliest and latest start-of-service time at location $i$. A non-negative service duration $s_{i}$ is considered for request locations $i\in\mathcal{P}\bigcup\mathcal{D}$, and $s_{0}=s_{2n+1}=0$ at the depots. 
Let $\mathcal{K}:=\{1,2,\ldots,K\}$ denote the set of identical vehicles with the capacity of $W_{\max}$. Each arc $(i,j)\in\mathcal{A}$ is associated with travel time $t_{ij}$. The maximum allowed ride time for passenger $i$ is $L^{i}_{\max}$.
In the context of disease exposed risk, we define a static risk score $r_{i}$ for every pick-up node $i\in\mathcal{P}$, which is a parameter estimated by the passenger's demographics and trip purpose. For passenger $i\in\mathcal{P}$, analogous to the passenger count ($w_{i}=-w_{n+i}$), we denote $r_{i}=-r_{n+i}$ associated with the pick-up and drop-off nodes, {which helps to track the onboard passengers' risk scores.}

For a specific vehicle $k\in\mathcal{K}$, let the binary variable $x_{ij}^{k}=1$ if vehicle $k$ traverses arc $(i,j)$. For each node $i\in\mathcal{N}$, let $A_{i}^{k}$ denote the start-of-service time and $L_{i}^{k}$ represents the ride time for request $i$. Note that, at node $i\in\mathcal{N}$, we allow the early arrival, where one will wait until the start of the time window $a_{i}$. Besides, one can wait and delay the start of the service at the latest to $b_{i}$ to align with the feasibility of maximum ride times. Let $W_{i}^{k}$ and $R_{i}^{k}$ denote the number of onboard passengers and total risk scores of onboard passengers at node $i\in\mathcal{N}$, respectively.  At each node $i\in\mathcal{N}$, $W_{i}^{k}$ and $R_{i}^{k}$ are propagated in a similar fashion by addition and subtraction of the onboard passengers' demand $w_{i}$ and risk score $r_{i}$ based on the pick-up and drop-off sequences. 

To formulate the risk components in the RDARP, we define two additional variables $Q_{i}^{k}$ and $H_{i}^{k}$. Specifically, $Q_{i}^{k}$ is the cumulative exposed risk along the vehicle $k$'s route after visiting node $i$.
Note that pathogen can remain in the vehicle environment for hours, the cumulative risk at destination depot $Q_{2n+1}^{k}$ can then be used to restrain excessive long trips (e.g., not exceeding $Q_{\max}$) and help guide the regular disinfection of vehicles. 
For each traversed arc $(i,j)$, the exposed risk is associated with onboard passengers' risk score $R_{i}^{k}$ and onboard duration on arc $(i,j)$, including travel time, service time, and potential waiting time or delay. 
Moreover, $H_{i}^{k}$ is the individual exposed risk for request $i$'s trip from $i$ to $n+i$ on vehicle $k$, which is equivalent to the increment between $Q_{i}^{k}$ and $Q_{n+i}^{k}$ after subtracting the self-transmission regarding the risk score $r_{i}$ and onboard duration from $i$ to $n+i$. Finally, let $\overline{H}$ be a nonnegative variable representing an upper bound of the exposed risk among all passengers. 

\subsection{Arc-based Formulation}\label{sec:arc-based_formulation}
We first present the RDARP as a bi-objective mixed-integer non-linear programming problem, established based on the three-index arc-based formulation~\citep{cordeau2006branch}.

\begin{subequations}\label{eq:arc_based_formulation}
\begin{align}
\min_{x}~&{\Biggl\{\sum_{k\in\mathcal{K}}\sum_{i\in\mathcal{N}}\sum_{j\in\mathcal{N}}t_{ij}x_{ij}^{k}, \  \overline{H} \Biggl\}} &&\label{eq:bi_objective_function}\\
\text{s.t.}~
    &\sum_{k\in\mathcal{K}} H_{i}^{k} \le \overline{H} && \forall i\in\mathcal{P}~{\label{eq:arc_min_max_risk}}\\
    &\sum_{k\in\mathcal{K}}\sum_{j\in\mathcal{N}}x_{ij}^{k}=1  && \forall i \in \mathcal{P}~\label{eq:served_once}\\
    &\sum_{j\in\mathcal{N}\setminus \{2n+1\}}x_{ji}^{k} - \sum_{j\in\mathcal{N}\setminus \{2n+1\}}x_{j,n+i}^{k} = 0  && \forall i \in \mathcal{P}, k \in \mathcal{K} ~\label{eq:served_by_one}\\
    &\sum_{j\in\mathcal{N}\setminus \{2n+1\}} x_{ji}^{k} - \sum_{j\in\mathcal{N}\setminus \{0\}} x_{ij}^{k} = 0 && \forall i \in \mathcal{P}\cup \mathcal{D},  k \in \mathcal{K} ~\label{eq:flow_conservation}\\
    &\sum_{j\in\mathcal{N}} x_{0j}^{k} = 1 && \forall k \in \mathcal{K}~\label{eq:depart_depot_once}\\
    &\sum_{j\in\mathcal{N}} x_{j,2n+1}^{k} = 1 && \forall k \in \mathcal{K}~\label{eq:arrive_depot_once}\\
    &A_{j}^{k}\ge A_{i}^{k} + s_{i} + t_{ij} + M_{1}\left(x_{ij}^{k} - 1\right) &&\forall i \in \mathcal{N}, j \in \mathcal{N},  k \in \mathcal{K} ~\label{eq:time_consistency}\\
    &a_{i} \le A_{i}^{k}\le b_{i} && \forall i \in \mathcal{N},  k \in \mathcal{K}~\label{eq:within_time_window}\\
    & L_{i}^{k} = A_{n+i}^{k} - \left(A_{i}^{k} + s_{i}\right) && \forall i \in \mathcal{P},  k \in \mathcal{K}~\label{eq:max_ride_time_calculation}\\
    & t_{i,n+i} \le L_{i}^{k} \le L^{i}_{\max} && \forall i \in \mathcal{P},  k \in \mathcal{K}~\label{eq:within_ride_time_limit}\\
    &W_{j}^{k}\ge W_{i}^{k} + w_{j} + M_{2}\left(x_{ij}^{k} - 1\right)&&\forall i \in \mathcal{N}, j \in \mathcal{N},  k \in \mathcal{K} ~\label{eq:capacity_consistency}\\
    &0 \le W_{i}^{k} \le W_{\max}  && \forall i \in \mathcal{N},  k \in \mathcal{K} ~\label{eq:within_capacity_limit}\\
    &R_{j}^{k}\ge R_{i}^{k} + r_{j} + M_{3}\left(x_{ij}^{k} - 1\right)&&\forall i \in \mathcal{N}, j \in \mathcal{N},  k \in \mathcal{K} ~\label{eq:on_vehicle_risk_consistency}\\
    &Q_{j}^{k}\ge Q_{i}^{k} + \left(A_{j}^{k}-A_{i}^{k}\right)R_{i}^{k} + M_{4}\left(x_{ij}^{k} - 1\right) &&\forall i \in \mathcal{N}, j \in \mathcal{N},  k \in \mathcal{K} ~\label{eq:risk_consistency}\\
    &0\le Q_{2n+1}^{k} \le Q_{\max}  &&  \forall k \in \mathcal{K} ~\label{eq:max_cumulative_risk_constraint}\\
    & H_{i}^{k} = Q_{n+i}^{k} - Q_{i}^{k} - \left(A_{n+i}^{k}-A_{i}^{k}\right)r_{i} && \forall i \in \mathcal{P},  k \in \mathcal{K}~\label{eq:direct_infection}\\
    &x_{ij}^{k} \in\{0,1\} &&\forall i \in \mathcal{N}, j \in \mathcal{N}, k \in \mathcal{K}~\label{eq:arc_based_decision_variable}
\end{align}
\end{subequations}

The objective function~\eqref{eq:bi_objective_function} aims at minimizing the total travel time and the maximum exposed risk among all passengers.
Constraints~\eqref{eq:arc_min_max_risk} serves as a min-max linearization, where the individual exposed risk for passenger $i\in\mathcal{P}$ is bounded by $\overline{H}$.
Constraints~\eqref{eq:served_once} and~\eqref{eq:served_by_one} ensure that each passenger is visited exactly once and one passenger's pick-up and drop-off nodes are served by the same vehicle. Constraints~\eqref{eq:flow_conservation} are the flow conservation constraints. Constraints~\eqref{eq:depart_depot_once} and~\eqref{eq:arrive_depot_once} enforce that each vehicle departs and arrives at the depot exactly once. Constraints~\eqref{eq:time_consistency} ensure the arrival time consistency by considering the travel time and service time, which is bounded by Constraints~\eqref{eq:within_time_window}. Constraints~\eqref{eq:max_ride_time_calculation} impose the ride time for each passenger, bounded by Constraints~\eqref{eq:within_ride_time_limit}. Constraints~\eqref{eq:capacity_consistency} define the conservation of the onboard number of passengers at each intermediate node, and Constraints~\eqref{eq:within_capacity_limit} limit the total number of onboard passengers by the vehicle capacity. The remaining constraints are associated with contagion risk calculation and propagation, and readers may refer to Appendix~\ref{sec:risk_measure_calculation} for an illustrative example. Specifically, Constraints~\eqref{eq:on_vehicle_risk_consistency} define the cumulative per-unit-time risk score at the intermediate node, which is further used in
Constraints~\eqref{eq:risk_consistency} to calculate the propagation of cumulative risk by multiplying the lapse of in-vehicle time between two consecutive nodes. To avoid excessive viruses suspended in the vehicle, we enforce a maximum cumulative risk limit $Q_{\max}$ by Constraints~\eqref{eq:max_cumulative_risk_constraint}. With cumulative risk in Constraints~\eqref{eq:risk_consistency}, the individual risk exposure can therefore be calculated by subtracting the amount of self-exposure from the total risk as in Constraints~\eqref{eq:direct_infection}.

Next, we remark on the rationale behind the inclusion of minimizing the \emph{maximum} individual exposed risk as one of the objectives and the challenges in solving the arc-based formulation above.
\begin{remark}\label{remark:minimized_max_risk}
(The role of minimized maximum individual risk)
To understand the use of $\overline{H}:=\max\limits_{i\in\mathcal{P}}(\sum\limits_{k\in\mathcal{K}}H_{i}^{k})$ in the objective function, one may consider the disease-spreading process in the dial-a-ride system as an epidemic process over physical contacts among all passengers. In such a contact network, all passengers are the nodes, and the contact duration between a node pair serves as the weight on edges. 
Note that while some of the passengers $i\in\mathcal{P}$ may not be physically adjacent at the same time, they may still be connected due to the existence of the vehicles $k\in\mathcal{K}${, e.g., the interaction of drivers {and cleaning technicians} at the depot.} The spread of an infectious disease can therefore be modeled as a percolation process~\citep{meyers2007contact} and the level of risk is determined by a network invasion threshold similar to the idea of {the basic reproduction number $\overline{r}_0$, which is used to estimate the speed of disease spreading}~\citep[for a comprehensive survey, see][]{guerra2017basic}.
In particular, the disease transmission between two nodes (passengers $i$ and $j$) takes place with a transmission rate $\zeta_{ij}$ that is positively correlated to the {contact} duration between two passengers and an infected individual can fully recover at rate $\overline{r}_c$. Specifically, \cite{qian2021scaling} showed that the disease can be eliminated if the degree distribution of the network satisfies:
\begin{equation}
\max_{i\in\mathcal{N}}\left(\sum_{j\in\mathcal{N}}\zeta_{ij}\right)<\overline{r}_{c}
\end{equation}
where the left-hand side denotes the maximum individual cumulative exposed risk among all passengers ~\citep[see Proposition 1 in][]{qian2021scaling}, which coincides $\overline{H}$ in our case. Minimizing $\overline{H}$ can be considered as minimizing the upper bound of the disease threshold, where the optimization of the latter is challenging~\citep{overton1988minimizing}. Instead, the use of $\overline{H}$ provides a linear function that can be embedded nicely into the RDARP framework and is shown to reach an exponential reduction in the number of infections~\citep{qian2021connecting}. 
\end{remark}

Finally, the arc-based formulation is known to be inefficient due to the weak LP relaxation~\citep{barnhart1998branch} in the BC process. Our RDARP problem further exacerbates the weakness due to the introduction of bilinear constraints (Constraints~\eqref{eq:risk_consistency}) as the product of two variables, arrival time and cumulative risk score.
This will introduce the spatial branch and bound algorithm to handle the bilinear constraints inside the BC framework, which requires exponential computation time for single-objective problems, letting alone the further challenge arising from the bi-objective function and the min-max term in our setting. We hence turn to the BCP framework and show that the problem can be efficiently solved with a new labeling algorithm by resolving the bilinear constraints within the pricing subproblems. 

\subsection{Trip-based Formulation}\label{sec:trip_based_formulation}
In this section, we reformulate the RDARP into a set-partitioning problem based on the DW Decomposition~\citep{dantzig1960decomposition}, which takes the form of a trip-based formulation. 

Let $\Omega$ be the set of elementary trips satisfying the RDARP constraints (see Eqs.~\eqref{eq:served_by_one}-\eqref{eq:direct_infection}). A feasible route $r:=\left(0,\cdots,2n+1\right)\in \Omega$ must satisfy the following conditions for node $i\in r$: (1) if $i\in \mathcal{P}$ then $n+i\in r$; if $i\in\mathcal{D}$ then $i-n\in r$, (2) $A_{i}\in[a_i,b_i]$, (3) $W_{i}\leq W_{\max}$, (4) $0<L_{i}\leq L_{\max}^{i}, \forall i\in\mathcal{P}$, and (5) $Q_{i}\leq Q_{\max}$. The travel cost of route $r\in\Omega$ is denoted by $c_{r}$; parameter $H_{ir}$ is the individual exposed risk of request $i\in \mathcal{P}$ on route $r$. If request $i$ is not part of route $r$, we set $H_{ir}=0$. The maximum exposed risk on route $r$ is bounded by variable $\overline{H}$.
Parameter $\alpha_{ir}$ is a binary parameter indicating whether request $i$ is visited in route $r$. 
Binary variable $\lambda_{r}$ equals 1 if route $r$ is selected to be part of the solution and $0$ otherwise. The trip-based formulation takes the following form.

\begin{subequations}\label{eq:trip_based_formulation}
\begin{align}
    \min& \Biggl\{\sum_{r\in\Omega} c_{r} \lambda_{r}, \ {\overline{H}}\Biggl\}  ~\label{obj:original_master_problem}\\
    \text{s.t.}& {\sum_{r\in\Omega} H_{ir} \lambda_{r} \le \overline{H}} &&\forall i \in \mathcal{P}\label{eq:trip_based_upperbound_risk} \\ 
    &\sum_{r\in\Omega} \alpha_{ir} \lambda_{r} = 1 &&\forall i \in \mathcal{P}\label{eq:set_partitioning}\\
    &\sum_{r\in\Omega}\lambda_{r} \leq K\label{eq:set_partitioning_fleet_size}\\
    &\lambda_{r} \in \{0,1\} && \forall r\in\Omega~\label{cons:integer_binary_mp}
\end{align}
\end{subequations}
where the exposed risk of each request is bounded by $\overline{H}$ in Constraints~\eqref{eq:trip_based_upperbound_risk}. Constraints~\eqref{eq:set_partitioning} guarantee that each request is served exactly once while Constraint~\eqref{eq:set_partitioning_fleet_size} poses a cap on the number of vehicles used, which cannot exceed the fleet size, {$K$}. The trip-based formulation decomposes the problem into trips $r\in\Omega$, each operated by one vehicle, with route-dependent risk information $H_{ir}$ for passenger $i$ when served on trip $r$. Note that in the arc-based formulation, the arrival time to each passenger is regulated using a decision variable, and therefore, the model can arrange the arrival times such that the maximum individual exposed risk is minimized. The trip-based formulation does not possess a mechanism to arrange the arrival times at passenger locations along the route to guarantee that the maximum individual exposed risk is minimized. Therefore, the time schedule along a trip $\Gamma_r$ shall be optimized to meet the Minimized-Maximum Risk (MMR) condition, defined as follows. 

\begin{definition}[\textbf{MMR condition}] \label{def:MMR}
    A route $r$, together with a given arrival time schedule $\Gamma_r$, are said to satisfy the MMR condition if the maximum total exposed risk of any individual on the route, $\max\limits_{i\in \mathcal{P}\cap r}H_{ir}$, can no further be reduced by changing $\Gamma_r$.
\end{definition}

\begin{proposition}
\label{prop:equivalence}
The arc-based and trip-based formulations are equivalent if every route $r\in\Omega$ satisfies the MMR condition.
\end{proposition}
\proof{Proof.} See Appendix~\ref{proof:equivalence}.

Based on the above statements, the key to the trip-based formulation is the accurate calibration of MMR, which can be further decomposed into the identification of a feasible $\Gamma_r$ that minimizes the maximum individual risk along the route. Note that a passenger's exposed risk is accumulated from the contact times with other passengers, while the contact time can be broken down into arc-level travel time, node-level service time, and additional waiting time due to the early arrival at a pick-up node. Given a route $r$, we observe that both travel and service times are constant. Therefore, the only component of a route that can potentially be altered to minimize the maximum risk is the waiting times before the start of service at a location. This can be made possible by delaying the service time at the predecessor nodes but at the risk of violating the latest arrival time and/or maximum trip length constraints for some of the passengers. 
Therefore, the calibration of MMR requires strategically adjusting the start-of-service time at the pick-up nodes while ensuring the trip remains feasible for all passengers. 

Nevertheless, ensuring the MMR condition for all serviced passengers given a route $r$ is a nontrivial task for the following two reasons. First, there is an inherent trade-off between earlier and later start-of-service times. Delaying a pick-up can reduce onboard passengers' waiting time and decrease exposed risk. However, it may also lead to delayed arrival at subsequent nodes hence the violation of the latest arrival time constraints. Second, due to the interlacing of pick-up and drop-off sequences, the optimal time delay becomes a dynamic programming problem that requires new algorithm development for the pricing subproblem. Existing DARP studies based on trip-based formulation~\citep{ropke2009branch,gschwind2015effective} only concern the feasibility of $\Gamma_r$ for generated trips, without considering optimization around $\Gamma_r$ which is essential to meet MMR condition at the trip level.
In the following section, we will first introduce the basic BCP components for RDARP, and we will then return to our specific algorithm for optimizing schedule $\Gamma_r$ in Section~\ref{sec:REF_rdarp}. 
\section{Solution Approach}\label{sec:Solution Approach}
Since the cardinality of $\Omega$ is extremely large, the approach used is based on a Column Generation (CG) algorithm to solve a restricted linear relaxation of the trip-based formulation, namely the \emph{restricted linear master problem} (RLMP). Specifically, the RLMP is limited to a subset of feasible trips, which is of moderate size and can be directly solved by the simplex algorithm. The CG method is an iterative procedure embedded in the branch-and-bound framework. At each branching node, the RLMP and the pricing subproblem (SP) are solved iteratively. The RLMP is first solved with a set of initial columns, and the resulting dual variables are used to generate feasible columns to be included in the basis in the SP. The RLMP is then resolved with the updated set of columns. This process is repeated until no valid columns with negative reduced costs are identified. 
If an integer solution is obtained at the root node, it is considered optimal. Otherwise, the branching procedure is initiated, where we compare and update the best integer solution with the incumbent solutions identified at branching nodes. During the branching procedure, invalid columns are eliminated and new columns may be generated in compliance with the branching rules. Finally, the best IP solution is reported to be the optimal solution of the original problem~\citep{dayarian2015column}.

\paragraph{\textbf{Restricted Linear Master Problem (RLMP)}.}\label{sec:RLMP}
The RLMP is constructed by relaxing the integrality constraints~\eqref{cons:integer_binary_mp} in the trip-based formulation. Instead of handling an exponential number of columns, the RLMP is restricted to a subset of feasible trips $\Omega'\subset\Omega$. 
In addition, we adopt the $\epsilon$-constraint method to handle the bi-objective nature of the objective function~\eqref{obj:original_master_problem}. The $\epsilon$-constraint method proceeds by optimizing one objective function while considering the other objective function as a constraint with a predefined bound. Let $\epsilon^{\text{risk}}$ and $\epsilon^{\text{cost}}$ be the predetermined upper bounds of the maximum exposed risk and total travel cost. We will consider two RLMPs, $\BFP^{\text{cost}}(\epsilon^{\text{risk}})$ and $\BFP^{\text{risk}}(\epsilon^{\text{cost}})$, aiming to minimize the total travel cost and maximum exposed risk while imposing the other objectives with the upper bounds of $\epsilon^{\text{risk}}$ and $\epsilon^{\text{cost}}$, respectively. The RLMPs take the form of 

\begin{subequations}\label{eq:RLMP}
\begin{align}
    \BFP^{\text{cost}}(\epsilon^{\text{risk}}): & \min~ \sum_{r\in\Omega'} c_{r} \lambda_{r}&&  &&&\BFP^{\text{risk}}(\epsilon^{\text{cost}}): \min~\overline{H} \label{cons:obj_RLMP}&&&& &&&&&\\
    \text{s.t.}&~\sum_{r\in\Omega'} \alpha_{ir} \lambda_{r} = 1 &&\forall i \in \mathcal{P} &&&\text{s.t.}~\sum_{r\in\Omega'} \alpha_{ir} \lambda_{r} = 1 &&&& \forall i \in \mathcal{P} &&&&& {(\pi_{i} \in \mathbb{R})}~\label{cons:visit_once_dual_pi}\\
    &\sum\limits_{r\in\Omega'}\lambda_{r} \leq K &&  &&& \sum\limits_{r\in\Omega'}\lambda_{r} \leq K &&&& &&&&& {(\mu \le 0)}~\label{cons:fleet_size_constraint}\\
    & \sum_{r\in\Omega'} H_{ir} \lambda_{r} \leq \epsilon^{\text{risk}} &&\forall i \in \mathcal{P} &&& \sum_{r\in\Omega'} H_{ir} \lambda_{r} \leq \overline{H} &&&& \forall i \in \mathcal{P}&&&&& {(\rho_{i} \le 0)} ~\label{cons:max_risk_dual_rho_node_level} \\
    & && &&& \sum_{r\in \Omega'}c_{r}\lambda_{r} \le \epsilon^{\text{cost}}&&&& &&&&&  {(\xi\leq 0)}~\label{cons:RLMP_total_travel_cost}\\ 
    &\lambda_{r} \geq 0 &&  \forall r\in\Omega' &&& \lambda_{r} \geq 0 &&&& \forall r\in\Omega'&&&&& 
\end{align}
\end{subequations}
where $\pi:= \{\pi_{i}\in\mathbb{R}:\forall i \in \mathcal{P}\}$, $\mu\le0$, $ \rho:= \{\rho_{i}\le0:\forall i \in \mathcal{P}\}$, and $\xi\leq 0$ are the dual variables associated with Constraints~\eqref{cons:visit_once_dual_pi}-\eqref{cons:RLMP_total_travel_cost}, respectively. 
In $\BFP^{\text{cost}}(\epsilon^{\text{risk}})$, one can convert the set partitioning constraints~\eqref{cons:visit_once_dual_pi} to set covering constraints. However, in $\BFP^{\text{risk}}(\epsilon^{\text{cost}})$, using the set covering constraint may result in multiple visits in some extreme cases. For instance, given the riskiest passenger with strictly positive risk $\overline{H}>0$, the passengers with no co-riders ($H_{ir}=0$) can be visited multiple times while satisfying Constraints~\eqref{cons:fleet_size_constraint}-\eqref{cons:RLMP_total_travel_cost}. 
Therefore, we keep the set partitioning form.

We use the $\epsilon$-constraint method to obtain the exact Pareto optimal solutions. We adopt an exact method~\citep{berube2009exact} by iteratively solving $\BFP^{\text{cost}}(\epsilon^{\text{risk}})$ and $\BFP^{\text{risk}}(\epsilon^{\text{cost}})$. At each iteration, the objective function value of the former is considered as the upper bound of $\epsilon^{\text{cost}}$ for the latter, and the latter is solved to avoid the dominated solutions~\citep{berube2009exact}. For brevity, we present the detailed algorithm with an illustration in Appendix~\ref{sec:Pareto_front_approximation}.

\paragraph{\textbf{Pricing Subproblem (SP).}}\label{sec:SP}
The columns in the subset $\Omega'$ are generated iteratively by solving the SP. The SP aims to identify the columns with negative reduced cost or prove the optimality of RLMP if none is found.
{The reduced cost functions $\tilde{c}_{r}^{\text{cost}}$ and $\tilde{c}_{r}^{\text{risk}}$, associated with route $r\in\Omega$, are shown in Eqs.~\eqref{eq:reduced_cost_sum_P1} and \eqref{eq:reduced_cost_sum_P2}.} Further, the reduced cost functions can be rewritten back to arc-level to incorporate the RDARP Constraints~\eqref{eq:served_by_one}-\eqref{eq:direct_infection}. Given an independent trip $r$, we omit the dimension $k$ in $x_{ij}^{k}$ and Constraints~\eqref{eq:served_by_one}-\eqref{eq:direct_infection}. The reduced costs take the form of

\begin{subequations}
\begin{align}
    & {\mathbf{P}_{\text{cost}}(\epsilon^{\text{risk}}): \tilde{c}_{r}^{\text{cost}} = c_{r} - \sum_{i\in\mathcal{P}} \alpha_{ir}\pi_{i} - \sum_{i\in\mathcal{P}} H_{ir}  \rho_{i} - \mu=\sum_{(i,j)\in r}\tilde{c}_{ij}^{\text{cost}}x_{ij}- \sum_{i\in\mathcal{P}} H_{ir} \rho_{i} - \mu}~\label{eq:reduced_cost_sum_P1}\\
    & {\mathbf{P}_{\text{risk}}(\epsilon^{\text{cost}}): \tilde{c}_{r}^{\text{risk}} = -\sum_{i\in\mathcal{P}} \alpha_{ir}\pi_{i} - \sum_{i\in\mathcal{P}} H_{ir}  \rho_{i} - \mu - \xi c_{r}=\sum_{(i,j)\in r}\tilde{c}_{ij}^{\text{risk}}x_{ij}- \sum_{i\in\mathcal{P}} H_{ir}  \rho_{i} - \mu}~\label{eq:reduced_cost_sum_P2}
\end{align}\label{eq:reduced_cost}
\end{subequations}
where $\tilde{c}_{ij}^{\text{cost}} = t_{ij}-\pi_{i}, \forall i \in \mathcal{P}$ and $\tilde{c}_{ij}^{\text{cost}} = t_{ij}$ otherwise; $\tilde{c}_{ij}^{\text{risk}} = -\xi t_{ij}-\pi, \forall i \in \mathcal{P}$ and $\tilde{c}_{ij}^{\text{risk}} = -\xi t_{ij}$ otherwise. Based on the expressions above, the pricing subproblem can be formulated as follows.

\begin{subequations}
\begin{align}
    \min & {\sum\limits_{(i,j)\in r}\tilde{c}_{ij}x_{ij}-\sum_{i\in\mathcal{P}}H_{ir}\rho_{i}-\mu}~\label{eq:reduced_cost_objective_funtion} \\
    \text{s.t.}& \eqref{eq:served_by_one}-\eqref{eq:direct_infection}
\end{align}
\label{eq:SP}
\end{subequations}
where the arc-level reduced cost $\Tilde{c}_{ij}=\Tilde{c}_{ij}^{\text{cost}}$ in $\mathbf{P}_{\text{cost}}(\epsilon^{\text{risk}})$ and $\Tilde{c}_{ij}=\Tilde{c}_{ij}^{\text{risk}}$ in $\mathbf{P}_{\text{risk}}(\epsilon^{\text{cost}})$.

The SP is solved in the same graph $\mathcal{G} = (\mathcal{N},\mathcal{A})$ by replacing the arc cost $t_{ij}$ with $\Tilde{c}_{ij}$, for arc $(i,j) \in \mathcal{A}$. The SP with  Constraints~\eqref{eq:served_by_one}-\eqref{eq:within_capacity_limit} is a typical Elementary Shortest Path Problem with Resource Constraints (ESPPRC), which is strongly NP-hard~\citep{dror1994note}. To guarantee the MMR condition (see Definition~\ref{def:MMR} and Proposition~\ref{prop:equivalence}), we impose additional constraints on our SP, resulting in the ESPPRC with the MMR condition (ESPPRCMMR). In contrast to previous studies that consider a \textit{linear cost function} of travel time $t_{ij}$~\citep{fukasawa2016branch,luo2017branch}, our study considers the exposed risk, which depends on the sequence of pick up and drop off of passengers and the varying onboard times including the potential waiting time due to early arrivals. To address this problem, we propose a dynamic programming-based algorithm, specifically the labeling algorithm, which is customized to satisfy the MMR condition in our RDARP. 

\subsection{Labeling Algorithm}\label{sec:labeling_algo}
We solve the SP (formulated as ESPPRCMMR) by means of the labeling algorithm~\citep{dumas1991pickup,feillet2004exact}.
In the labeling algorithm, each label represents a partial trip from origin $0$ to node $i\in\mathcal{N}\setminus\{0\}$ that is extended along the arcs in graph $\mathcal{G}=(\mathcal{N},\mathcal{A})$. The labeling algorithm aims to generate a complete trip that satisfies the RDARP constraints and contributes to the LP relaxation of the RLMP. 
Specifically, for each potential feasible label extension along arc $(i,j)\in\mathcal{A}$, the stored information (namely \textit{resources}) in the label are updated following \textit{resource extension functions} (REFs)~\citep{desaulniers1998unified} (see Sections \ref{sec:REFPDPTW}-\ref{sec:REF_rdarp}). Based on the updated resources, the feasibility of the label is checked regarding the RDARP constraints. However, the enumeration of label extensions can still result in an exponential size of explorations. In this regard, label dominance rules are required to terminate unpromising labels during the extension procedure (see Section \ref{subsec:dominance_darp}).

In our labeling algorithm, we define a label as follows.

\begin{definition}\label{def:label_definition}
    Let $l$ denote a label, representing a partial trip that starts from origin depot $0$ and resides at the current node $\eta_{l}$. The label can be encoded as 
    \begin{equation*}
    l:=\left(\underbrace{\eta_{l},\Tilde{c}_{l},A_{l},B_{l},W_{l},\mathcal{V}_{l},\mathcal{O}_{l}}_\text{PDPTW-related},\underbrace{\mathcal{O}_{l}^{a}, R_{l}, Q_{l}, \{h_{l}^{o}, d_{l}^{o}:o\in\mathcal{O}_{l}\cup\mathcal{V}_{l}\}}_\text{Risk-related},\underbrace{\{B_{l}^{o}, D_{l}^{o}(A_{l}),D_{l}^{o}(B_{l}^{o}): o \in \mathcal{O}_{l}\}}_\text{Ride time-related}\right)    
    \end{equation*}
\begin{itemize}
    \item $\eta_{l}$: the node where label $l$ resides
    \item $\Tilde{c}_{l}$: reduced cost after serving node $\eta_{l}$ 
    \item $A_{l}$: the earliest possible start-of-service time at node $\eta_{l}$ given a partial path $l$
    \item $\Gamma_{l}$: a sequence of time schedules that record the feasible start-of-service times, denoted as $\Gamma_{l}^{0}<\Gamma_{l}^{i}<\Gamma_{l}^{j}<\cdots$ for route $(0,i,j,\cdots)$
    \item $B_{l}$: the latest possible start-of-service time at node $\eta_{l}$ given a partial path $l$
    \item $W_{l}$: number of requests on board at node $\eta_{l}$
    \item $\mathcal{V}_{l}$: the set of serviced requests along the partial trip until node $\eta_{l}$
    \item $\mathcal{O}_{l}$: the set of open requests, where request $o\in\mathcal{O}_{l}$ has been picked up but not dropped off
    \item $\mathcal{O}_{l}^{a}$: the set of associated requests at $\eta_{l}$, which includes the requests that have been dropped off but \emph{were exposed to the open requests $o\in\mathcal{O}_{l}$}
    \item {$R_{l}$: sum of risk scores for all onboard requests $o\in\mathcal{O}_{l}$ at node $\eta_{l}$} 
    \item $Q_{l}$: cumulative risk at node $\eta_{l}$
    \item $h_{l}^{o}$: individual exposed risk of request $o\in\mathcal{O}_{l}\cup \mathcal{V}_{l}$ at node $\eta_{l}$
    \item $d_{l}^{o}$: available delay buffer for request $o\in\mathcal{O}_{l}\cup\mathcal{V}_{l}$ when starting from node $\eta_{l}$ (see Section \ref{sec:risk_min_without_associated})
    \item $D_{l}^{o}(\tau)$: latest possible drop-off time for request $o\in\mathcal{O}_{l}$ as a function of start-of-service time $\tau$
    \item $B_{l}^{o}$: a breaking point representing a feasible start-of-service time at node $\eta_{l}$ for an open request $o\in\mathcal{O}_{l}$; If $A_{l}\geq B_{l}^{o}$, it implies that $D_{l}^{o}(A_{l})$ is bounded by the latest time window for pick-up ($b_{o}$) or drop-off ($b_{o+n}$).
\end{itemize}
\end{definition}

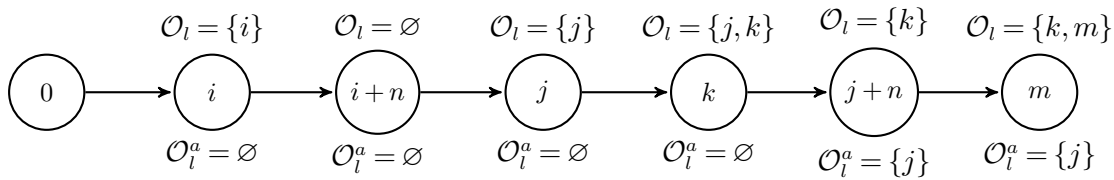
\begin{figure}[H]
\caption{Illustration for the set of open and associated requests on an interlaced trip}
\label{fig:illustration_open_associated}
    \begin{tikzpicture}[->,>=stealth',shorten >=1pt,auto,node distance=2.2cm,  thick,main node/.style={circle,draw,minimum size=1cm,font=\sffamily\small\bfseries}]
      \node[main node] (1) {$0$};
      \node[main node, label=above:{$\mathcal{O}_{l}=\{i\}$}, label=below:{$\mathcal{O}_{l}^{a}=\varnothing$}] (2) [right of=1] {$i$};
      \node[main node, label=above:{$\mathcal{O}_{l}=\varnothing$}, label=below:{$\mathcal{O}_{l}^{a}=\varnothing$}] (3) [right of=2] {$i+n$};
      \node[main node, label=above:{$\mathcal{O}_{l}=\{j\}$}, label=below:{$\mathcal{O}_{l}^{a}=\varnothing$}] (4) [right of=3] {$j$};
      \node[main node, label=above:{$\mathcal{O}_{l}=\{j,k\}$}, label=below:{$\mathcal{O}_{l}^{a}=\varnothing$}] (5) [right of=4] {$k$};
      \node[main node, label=above:{$\mathcal{O}_{l}=\{k\}$}, label=below:{$\mathcal{O}_{l}^{a}=\{j\}$}] (6) [right of=5] {$j+n$};
      \node[main node, label=above:{$\mathcal{O}_{l}=\{k,m\}$}, label=below:{$\mathcal{O}_{l}^{a}=\{j\}$}] (7) [right of=6] {$m$};
      \path[every node/.style={font=\sffamily\small}]
        (1) edge node [above]  {} (2)
        (2) edge node [above]  {} (3)
        (3) edge node [above]  {} (4)
        (4) edge node [above]  {} (5)
        (5) edge node [above]  {} (6)
        (6) edge node [above]  {} (7);
    \end{tikzpicture}
\end{figure}

In particular, $\eta_{l}$ is recorded for label $l$ to identify the labels with the same ending node, and $\Tilde{c}_{l}$ is used to prioritize the labels with lower reduced cost. The resources $A_{l}$ and $W_{l}$  are considered to keep track of capacity and time window constraints. $\mathcal{V}_{l}$ and $\mathcal{O}_{l}$ are designed for pairing and precedence constraints. Set $\mathcal{O}_{l}^{a}$ is introduced to ensure that the adjustment for the time schedule is in compliance with the MMR condition for both open requests as well as the associated requests, where an illustrative example is shown in Figure~\ref{fig:illustration_open_associated}. For example, at node $j+n$, passenger $k$ is the only onboard passenger hence $\mathcal{O}_{l}=\{k\}$. Due to exposure between $j$ and $k$ prior to the arrival, $j$ is an associated request to $k$ and $\mathcal{O}_{l}^a=\{j\}$. 

In addition, resources {$B_{l}$}, $D_{l}^{o}(\tau)$, and $B_{l}^{o}$ are introduced to handle the maximum ride time constraint. The main idea is to delay the start-of-service time such that potential waiting time can be reduced and total ride time is more likely to be feasible for the unknown subsequent nodes.

To ensure feasibility with respect to the risk-related components, we propose four risk-related resources, $R_{l}$, $\{h_{l}^{o}, d_{l}^{o}: o \in \mathcal{O}_{l}\cup\mathcal{V}_{l}\}$, and $Q_{l}$. The resource $R_{l}$ keeps track of total risk scores of visited requests $i\in\mathcal{O}_{l}\cup\mathcal{V}_{l}$. The individual exposed risk $h_{l}^{o}$ is propagated at the arc level along the partial path of label $l$ by considering the risk scores of co-riders and the duration of onboard travel, which can be calibrated by utilizing the delay buffer $d_{l}^{o}$. Lastly, $Q_{l}$ is used to verify the feasibility of the maximum cumulative risk constraints ($Q_{l}\leq Q_{\max}$).

\begin{algorithm}
\footnotesize %\small
\caption{Labeling algorithm}
\label{alg:label_generation}
$\mathcal{G}=(\mathcal{N},\mathcal{A})$, \{$a_{i},b_{i},n_i, s_i, r_i: \forall i\in\mathcal{N}\}, \{t_{ij}: \forall i,j \in \mathcal{A}\}$ \Comment{Input}\\
Label $l^{\star}$ with most negative reduced cost \Comment{Output}
\begin{algorithmic}[1] %[1] enables line numbers
    \State Initiate the label, {$l_{0}\leftarrow \left(0,0,a_{0},b_{0}, 0, \{0,b_{0}-a_{0}\}, \varnothing, \varnothing, 0, 0, \varnothing ,\varnothing\right)$}
    \State Initialize a priority queue: $\mathcal{Q}\leftarrow (\Tilde{c}_{l_{0}},l_{0})$. Initialize a collection of sets of associated labels with the same ending node, $\mathcal{L} = \{\mathcal{L}\left(i\right):i\in\mathcal{N}\}$. Initialize $\Tilde{c}^{\star}=0$.
    \While{$\mathcal{Q}$ is not empty}
    \State $l\leftarrow \mathcal{Q}.\text{pop}()$
        \ForAll{$j \in \mathcal{N}$} 
            \State \textbf{Update} $\eta_{l'}, A_{l'},W_{l'},\mathcal{V}_{l'},\mathcal{O}_{l'}$ following the REF rules in Eqs.~\eqref{eq:ref_ending_node}-\eqref{ref:open_request_set}.\Comment{PDPTW-related Updates}
            \If{$(\eta_{l},\eta_{l'})$ is a PDPTW-feasible extension}\Comment{DARP-related Updates}
                \ForAll{$o\in\mathcal{O}_{l'}$}
                    \State \textbf{Update} $B_{l'}^{o},D_{l'}^{o}(A_{l'}),D_{l'}^{o}(B_{l'}^{o})$ following REFs in Eqs.~\eqref{REF:B}-\eqref{ref:latest_possible_arrival_time}
                \EndFor
                \If{$(\eta_{l},\eta_{l'})$ is a DARP-feasible extension}\Comment{RDARP-related Updates}
                    \State \textbf{Update} $\mathcal{O}_{l'}^{a}, R_{l'}, Q_{l'},\{h_{l'}^{o}, d_{l'}^{o}:o\in\mathcal{O}_{l}\cup\mathcal{O}_{l}^{a}\},\Tilde{c}_{l'}$ following REFs in Eqs.~\eqref{ref:associated_request_set}-\eqref{REF:reduced_cost} and Algorithm~\ref{alg:contact_duration_calibriation}
                    \If{$(\eta_{l},\eta_{l'})$ is a RDARP-feasible extension}
                        \State  $l'\leftarrow\left(\eta_{l'},\Tilde{c}_{l'},A_{l'},W_{l'},\mathcal{V}_{l'},\mathcal{O}_{l'},\mathcal{O}_{l'}^{a}, \{h_{l'}^{o},d_{l'}^{o}:o\in\mathcal{O}_{l}\cup\mathcal{O}_{l}^{a}\}, R_{l'}, Q_{l'},\{B_{l'}^{o}, D_{l'}^{o}(A_{l'}),D_{l'}^{o}(B_{l'}^{o}): o \in \mathcal{O}_{l'}\}\right)$
                        \State Check if $l'$ is dominated by the existing labels in $\mathcal{L}\left(\eta_{l'}\right)$ following Eqs.~\eqref{dom:dominance_rule_darp0}-\eqref{dom:dominance_rule_darp3}
                        \If{$l'$ is not dominated}
                            \State $\mathcal{Q}\leftarrow \mathcal{Q}\cup\{(\Tilde{c}_{l'},l')\}$, $\mathcal{L}\left(\eta_{l'}\right)\leftarrow \mathcal{L}\left(\eta_{l'}\right)\cup\{l'\}$
                        \EndIf
                        \ForAll{$l_{k}\in \mathcal{L}\left(\eta_{l'}\right)$}
                            \If{$l'$ dominates $l_{k}$}
                                \State $\mathcal{Q}\leftarrow \mathcal{Q}\setminus\{(\Tilde{c}_{k},l_{k})\}$, $\mathcal{L}\left({\eta_{l'}}\right)\leftarrow \mathcal{L}\left(\eta_{l'}\right)\setminus\{l_{k}\}$
                            \EndIf
                        \EndFor
                        \If{$\eta_{l'}=2n+1$ and $\Tilde{c}_{\eta_{l'}}<\Tilde{c}^{\star}$}
                            \State Update the current best reduced cost and label: $\Tilde{c}^{\star}=\Tilde{c}_{\eta_{l'}}$, $l^{\star}=l'$
                        \EndIf
                        % \EndIf
                    \EndIf
                \EndIf
            \EndIf
        \EndFor
    \EndWhile
\end{algorithmic}
\end{algorithm}

Our labeling algorithm proceeds as shown in Algorithm~\ref{alg:label_generation}. The labeling algorithm describes a label extension procedure from the origin depot $0$ to final depot $2n+1$. For the label extension from $l$ to $l'$, the feasibility checks and the updates of related resources consist of three levels, which are detailed as follows in compliance with line numbers in Algorithm~\ref{alg:label_generation}. 
\textbf{PDPTW} (line 7): The PDPTW-related resources are first updated to the current node $\eta_{l'}$. An extension is said to be PDPTW-feasible if it satisfies the time window, capacity, as well as pairing and precedence constraints, corresponding to $A_{l}+s_{\eta_{l}}+t_{\eta_{l},\eta_{l'}}\leq b_{\eta_{l'}}$, $W_{l}+w_{\eta_{l'}}\leq W_{\max}$, $\eta_{l'}-n\in\mathcal{O}~\text{if}~ \eta_{l'}\in\mathcal{D}$ $\eta_{l'}\not\in\mathcal{V}_{l}\text{~if~}\eta_{l'}\in\mathcal{P}$, and $\mathcal{O}=\varnothing$\text{~if~}$\eta_{l'}=2n+1$, respectively. 
\textbf{DARP} (line 11): For the ride time constraint, we verify the existence of a feasible time schedule by checking if there is a valid drop-off time $A_{l'}$ within the feasible range of the time window and ride time constraints, such that $A_{l'}\leq D_{l'}^{o}(A_{l'})$. 
\textbf{RDARP} (line 13): The risk-related resources are updated to obtain the individual risks and cumulative exposed risk along the route. The extension is feasible with respect to RDARP if $Q_{l'}\leq Q_{\max}$.
Finally, the label dominance rules (lines 16-26) are conducted to eliminate unpromising labels, thus accelerating the labeling procedure. The labeling algorithm terminates when all non-dominated labels are extended and reached node $2n+1$. The final output is a set of RDARP-feasible trips. If no trip with negative reduced cost $\Tilde{c}_{l}<0$ is found, then the RLMP is proven to be optimal. 

We note that the proposed labeling algorithm is tailored to handle the risk calibration process in the SP with three major differences as compared to existing ones for DARP~\citep{ropke2009branch, qu2015branch, gschwind2015effective}. First, we add a family of risk-related resources, $R_{l}$, $Q_{l}$, $h_{l}^{i}$, and $d_{l}^{i}$, to capture the risk propagation dynamics and update individual exposed risk $h_{l}^{i}$ with label extensions. Second, we introduce new resources, dominance rule, and risk calibration algorithm collectively reduce the search space and guarantee the simultaneous satisfaction of risk-related constraints and the standard DARP constraints (ensuring feasibility). Finally,  we propose a risk calibration algorithm (proceed in line 14 in Algorithm 1, see details in Algorithm~\ref{alg:contact_duration_calibriation}) that tracks and updates the delaying potential for each request, which is used to determine the optimal $\Gamma_r$ during label extension and guarantee the satisfaction of the MMR condition (ensuring optimality).

\paragraph{\textbf{Resource Extension Functions.}}\label{subsec:REF_darp}
The REFs either initiate the resources or update the resources following a label extension. By updating the resources in the RDARP, the feasibility of the partial route (the PDPTW constraints, maximum ride time constraints, and maximum cumulative risk constraints) is sequentially checked. The REFs for the RDARP are described in three parts: PDPTW-related resources, DARP-related resources, and RDARP-related resources in Sections \ref{sec:REFPDPTW}-\ref{sec:REF_rdarp}.

\subsubsection{\textbf{REF for PDPTW.}}\label{sec:REFPDPTW}
While extending labels $l$ to $l'$, we update the first set of resources following the REFs~\eqref{eq:ref_ending_node}-\eqref{ref:open_request_set}, which have been applied to the PDPTW.

\begin{subequations}
\begin{align}
    \eta_{l'} &\leftarrow \eta_{l}\label{eq:ref_ending_node}\\
    A_{l'} &\leftarrow \max\{A_{l}+s_{\eta_{l}}+t_{\eta_{l},\eta_{l'}}, a_{\eta_{l'}}\}\label{REF:early_time}\\
    W_{l'} &\leftarrow W_{l} + w_{\eta_{l'}}\label{REF:capacity}\\
    \mathcal{V}_{l'} &\leftarrow \mathcal{V}_{l}\cup \{\eta_{l'}-n\}\text{~if~}\eta_{l'}\in\mathcal{D}\label{REF:visited_nodes}\\
    \mathcal{O}_{l'} &\leftarrow
        \begin{cases}
            \mathcal{O}_{l}\cup\{\eta_{l'}\} & \text{if}~\eta_{l'}\in \mathcal{P}\\
            \mathcal{O}_{l}\setminus\{\eta_{l'}-n\} & \text{if}~\eta_{l'}\in \mathcal{D}\\
        \end{cases}\label{ref:open_request_set}
\end{align}
\label{eq:ref_pdptw}
\end{subequations}
Eq.~\eqref{eq:ref_ending_node} updates the current node from $\eta_{l}$ to $\eta_{l'}$. Eqs.~\eqref{REF:early_time}-\eqref{REF:capacity} update the earliest arrival time and number of requests onboard. 
Eq.~\eqref{REF:visited_nodes} updates the set of serviced nodes by adding the corresponding request $\eta_{l'}-n$ after dropping off at node $\eta_{l'}\in\mathcal{D}$. Finally, Eqs.~\eqref{ref:open_request_set} record the set of open requests after serving node $\eta_{l'}$.

\subsubsection{\textbf{REF for DARP.}}
If the arc extension $(\eta_{l},\eta_{l'})$ is feasible according to the PDPTW constraints, we update the set of resources for all open requests $o\in\mathcal{O}_{l'}$ to account for constraints related to maximum ride time and risk. For each $o\in\mathcal{O}_{l'}$, we update the resources related to ride time as the REFs below.

\begin{subequations}\label{eq:ref_darp}
\begin{align}
    B_{l'}^{o}&\leftarrow
        \begin{cases}
            \min\{B_{l'}, b_{\eta_{l'}+n}-s_{\eta_{l'}}-L_{\max}^{\eta_{l'}}\} & \text{if}~ \eta_{l'}=o\\
            \max\{A_{l'},\min\{B_{l}^{o}+s_{\eta_{l}}+t_{\eta_{l},\eta_{l'}},B_{l'}\}\} & \text{otherwise}\\
        \end{cases}\label{REF:B}\\
    D_{l'}^{o}(A_{l'})&\leftarrow
        \begin{cases}
            \min\{A_{l'} + s_{\eta_{l'}} + L_{\max}^{\eta_{l'}}, b_{\eta_{l'}+n}\}& \text{if}~ \eta_{l'}=o\\
            D_{l}^{o}(A_{l}) + \min\{A_{l'}-s_{\eta_{l}}-t_{\eta_{l},\eta_{l'}} - A_{l}, B_{l}^{o}-A_{l}\}& \text{otherwise}\\
        \end{cases}\label{REF:ldT}\\
    D_{l'}^{o}(B_{l'}^{o})&\leftarrow
        \begin{cases}
            B_{l'}^{o} + s_{\eta_{l'}} + L_{\max}^{\eta_{l'}}& \text{if}~ \eta_{l'}=o\\
            D_{l}^{o}(B_{l}^{o}) - \max\{0, B_{l}^{o}+s_{\eta_{l}}+t_{\eta_{l},\eta_{l'}} - B_{l'}\}& \text{otherwise}\\
        \end{cases}\label{REF:ldB}\\
    B_{l'}&\leftarrow
    \begin{cases}
        b_{\eta_{l'}} & \text{if} ~\eta_{l'}\in \mathcal{P}\\
        \min\{b_{\eta_{l'}}, D^{\eta_{l'}-n}_{l}(B_{l}^{\eta_{l'}-n})\} &\text{if}~ \eta_{l'}\in \mathcal{D}\\
    \end{cases}\label{ref:latest_possible_arrival_time}
\end{align}
\end{subequations}

The REFs in Eqs.\eqref{REF:B}--\eqref{ref:latest_possible_arrival_time} were originally proposed by \cite{gschwind2015effective} to handle the maximum ride time constraints. The main idea is to update two distinct time stamps for each open request on the partial route, for every label extension, to specify the \textit{range} of the latest possible drop-off time of the particular open request. These time stamps are referred to as \textit{dynamic time windows} (DTW) proposed by~\cite{gschwind2015effective}. Specifically, in Appendix~\ref{appendix:description_ride_time_resources}, we present two particular cases to better understand the properties of DTWs: (1) initiating the DTW for a new open request, and (2) updating the DTW for an existing open request.

\subsubsection{REF for RDARP}\label{sec:REF_rdarp} 
In the labeling algorithm, the risk calibration process is invoked when a feasible DARP extension is identified along a partial path. The purposes of the risk calibration, upon the extension, are to identify possible waiting due to early arrival, compute the maximum available resources for adjusting start-of-service time, and update the up-to-date value of MMR. 
To accomplish the tasks, we will make use of three auxiliary variables as follows:

\begin{itemize}
    \item $\Delta d_{\eta_{l},\eta_{l'}}^{i}$: the start-of-service delay buffer for request $i$ on label extension $(\eta_{l},\eta_{l'})$, which corresponds to the maximum amount of time that an open or associated request $i\in\mathcal{O}_{l}\cup \mathcal{O}_{l}^{a}$ can delay its start-of-service at node $\eta_{l}$, while extending to $\eta_{l'}$ without violating the arrival time and maximum ride time constraints. We use $\Delta d_{\eta_{l},\eta_{l'}}$ without $i$ in the superscript to represent the set $\{\Delta d_{\eta_{l},\eta_{l'}}^{i}:i\in\mathcal{O}_{l}\cup \mathcal{O}_{l}^{a}\}$
    \item $\Delta \omega_{\eta_{l},\eta_{l'}}$: waiting time due to the early arrival at node $\eta_{l'}$ that can be offset by delaying the start-of-service at node $\eta_{l}$
    \item $\delta_{\eta_{l},\eta_{l'}}^{i}$: the actual amount of start-of-service delay for each request $i\in\mathcal{O}_{l}\cup\mathcal{O}_{l}^{a}$ to achieve the MMR condition along the partial route associated with label $l'$. We use $\delta_{\eta_{l},\eta_{l'}}$ without $i$ in the superscript to represent the set $\{\delta_{\eta_{l},\eta_{l'}}^{i}:i\in\mathcal{O}_{l}\cup \mathcal{O}_{l}^{a}\}$
\end{itemize}

With the above notation, for a label $l$ extending from $\eta_l$ to $\eta_{l'}$, the anticipated maximum risk $\BFH_{\eta_{l},\eta_{l'}}(\cdot)$ can be minimized following the Bellman Equation as below:
\begin{equation}
    \BFH_{\eta_{l},\eta_{l'}}(l, \Delta d_{\eta_{l},\eta_{l'}})= \min_{\delta_{\eta_{l},\eta_{l'}}} \big\{\Delta \BFH_{\eta_{l},\eta_{l'}}(\delta_{\eta_{l},\eta_{l'}}) + \BFH_{\eta_{l'},*} ( l', \Delta d_{\eta_{l',*}})\big\}
    \label{eq:bellman_mmr}
\end{equation}
where we use $\Delta \BFH_{\eta_{l},\eta_{l'}}$ to denote the realized risk increment to the maximum risk by extending from $\eta_l$ to $\eta'_l$, and $\BFH_{\eta_{l'},*}$ as the future realization of maximum risk by further extension from $\eta'_l$ with remaining delaying resources $\Delta d_{\eta_{l',\star}}$. 
In general, solving Eq.~\eqref{eq:bellman_mmr} needs to consider the trade-off between spending resources (e.g., the allowable start-of-service delay $\Delta d_{\eta_{l},\eta_{l'}}$) for immediate risk reduction and the potential gain from future extensions. Nevertheless, we can leverage problem-specific structures to decompose the dependency between current and future risk realizations based on the following series of propositions.

\begin{proposition}\label{prop:R_H_non_increasing}
    Both $\Delta \BFH_{\eta_{l},\eta_{l'}}$ and $\BFH_{\eta_{l'},\star}$ are  non-increasing function with respect to the increase in $\delta_{\eta_{l},\eta_{l'}}$. 
\end{proposition}
\proof{Proof.} See Appendix~\ref{proof:R_H_non_increasing}.

\begin{proposition}\label{prop:optimality}
    Let $\bar{\delta}_{\eta_{l},\eta_{l'}}$ be the minimum possible delay that minimizes $\Delta \BFH_{\eta_{l},\eta_{l'}}$, and $\delta_{\eta_{l},\eta_{l'}}^{\star}$ be one feasible delay that minimizes Eq.~\eqref{eq:bellman_mmr}. Then it holds true that $\delta_{\eta_{l},\eta_{l'}}^{\star}\geq \bar{\delta}_{\eta_{l},\eta_{l'}}$. Moreover, delaying only $\bar{\delta}_{\eta_{l},\eta_{l'}}$ at node $\eta_{l'}$ and delaying  $\delta_{\eta_{l},\eta_{l'}}^{\star}-\bar{\delta}_{\eta_{l},\eta_{l'}}$ at the subsequent node extended from $\eta_{l'}$ also minimize Eq.~\eqref{eq:bellman_mmr}.
\end{proposition}
\proof{Proof.} See Appendix~\ref{proof:optimality}.

Following Proposition~\ref{prop:optimality}, it is evident that one can always bank the extra amount of delay beyond $\bar{\delta}_{l'}$ for future extensions while still satisfying the MMR condition. This gives rise to an appealing approach where we only offset the minimum amount of delay ($\bar{\delta}_{l'}$) at each visited node and wait until further extensions for additional delay adjustment. It allows us to perform risk calibration and delay adjustment without seeing future extensions, hence allowing seamless integration of maximum risk minimization procedures into existing labeling algorithms. More importantly, this avoids the risk of violating DARP feasibility constraints and the risk of eliminating feasible label extensions in future extensions due to excessive local delay. The reason is that the incurred minimum amount of delay will not drive the start-of-service time to exceed the earliest possible start-of-service time following DARP label extensions. As a consequence, delaying by exactly the minimum amount $\bar{\delta}_{l'}$ at each label extension $(\eta_{l},\eta_{l'})$ results in one such solution that guarantees the minimization of maximum risk while ensuring the feasibility of the generated paths.

With the above understanding, the risk calibration and optimization will proceed by comparing potential waiting time ($\Delta \omega_{\eta_{l},\eta_{l'}}$) with the maximum allowable delay buffer ($\Delta d_{\eta_{l},\eta_{l'}}$). The goal is to identify the optimal amount of delay $\delta_{\eta_{l},\eta_{l'}}^{i}$ for each request $i\in\mathcal{O}_{l}\cup\mathcal{O}_{l}^{a}$ to satisfy the MMR condition. The rules for delay optimization will differ depending on if the set of associated requests is empty, which can be divided into the following two scenarios: 
\begin{figure}[H]
    \centering
    \caption{Two cases for risk calibration}
    \includegraphics[width=1.0\textwidth]{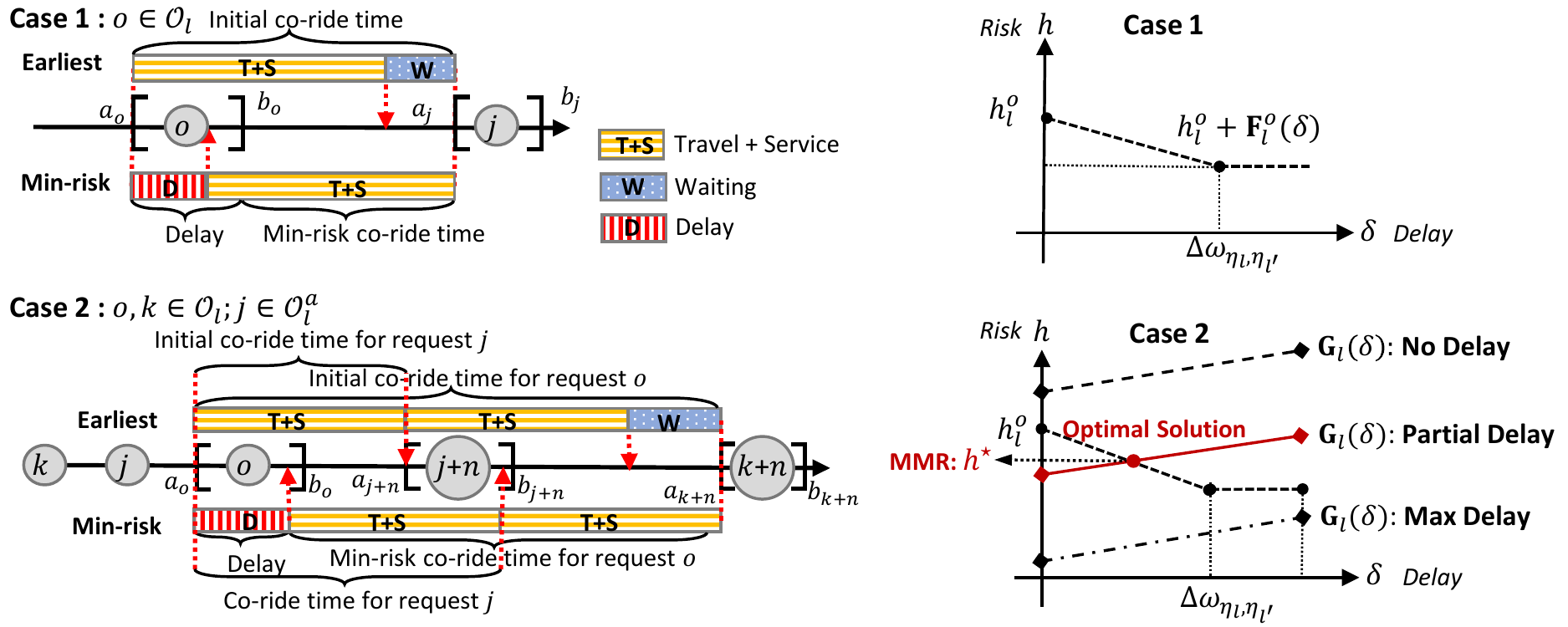}
    \label{fig:risk_calibration_3cases}
\end{figure}

\paragraph{\textbf{Risk minimization without associated requests.}}\label{sec:risk_min_without_associated}
\begin{subequations}
The most straightforward scenario corresponds to the case where $\mathcal{O}_{l}\neq\varnothing$ and $\mathcal{O}_{l}^{a}=\varnothing$. Under this setting, when extending a label $(\eta_{l}, \eta_{l'})$ associated with a particular node $\eta_{l}$, we only focus on the resources available from the present open requests. We introduce a characteristic function, denoted as $\BFF_{l}^{o}(\delta)$, to capture the change of individual exposed risk for an open request $o\in\mathcal{O}_{l}$ as follows:

\begin{equation}
    \BFF_{l}^{o}(\delta):= \sum\limits_{i\in\mathcal{O}_{l}\setminus\{o\}}\left(\underbrace{A_{l'}-A_{l}-\min\{\delta,\omega_{\eta_{l},\eta_{l'}}\}}_\text{Co-ride travel time with $o$}\right)\cdot \underbrace{r_{i}}_\text{Risk of co-rider}, \forall o\in\mathcal{O}_{l}
\end{equation}
where the function $\BFF_{l}^{o}(\delta)$ is a piecewise non-increasing function with the increased amount of delay $\delta$. $\BFF_{l}^{o}(\delta)$ will first decrease for $\delta\leq \omega_{\eta_{l},{\eta}_{l'}}$ and remain fixed for $\delta>\omega_{\eta_{l},\eta_{l'}}$, as shown in Case 1 in Figure~\ref{fig:risk_calibration_3cases}. 
Note that the amount of delay $\delta$ at node $\eta_{l}$ is upper bounded by the available delay buffer $\Delta d_{\eta_{l},\eta_{l'}}^{o}$. 
To minimize the increment of exposed risk on arc $(\eta_{l},\eta_{l'})$, we need to maximally utilize the available delay buffer to offset the waiting time. In this case, the actual delay $\delta_{\eta_{l},\eta_{l'}}^{o}$ can be computed as:
     \begin{equation}
    \delta_{\eta_{l},\eta_{l'}}^{o}=\min\{\Delta \omega_{\eta_{l}, \eta_{l'}}, \Delta d_{\eta_{l},\eta_{l'}}^{o}\}, \forall o\in \mathcal{O}_l
        \label{eq:delta_update_algo2}
     \end{equation}
Specifically, the adjustable waiting time on arc $(\eta_{l},\eta_{l'})$ is denoted by:
    \begin{equation}
        \Delta \omega_{\eta_{l},\eta_{l'}}= \min\{A_{l'}-\left(A_{l}+s_{\eta_l}+t_{\eta_{l},\eta_{l'}}\right),b_{\eta_{l}}-A_{l}\}\label{eq:waiting_time_algo2}
    \end{equation}    
where $A_{l} + s_{\eta_{l}} + t_{\eta_{l},\eta_{l'}}$ corresponds to the actual arrival time at node $\eta_{l'}$. The waiting time is measured as the difference between the start-of-service time $A_{l'}$ and the actual arrival time at node $\eta_{l'}$, which is further bounded by the latest time window at node $\eta_{l}$, denoted by $b_{\eta_{l}}$.

Meanwhile, the amount of delay buffer $\Delta d_{\eta_{l},\eta_{l'}}^{o}$ is bounded by the DARP feasibility constraints. For each node-to-extend $\eta_{l'} \in \mathcal{P}$, the most relaxed $\Delta d_{\eta_{l},\eta_{l'}}^{\eta_{l'}}$ is capped by the difference between the earliest and latest possible start-of-service time, denoted as $\Delta d_{\eta_{l},\eta_{l'}}^{\eta_{l'}} = B_{l'} - A_{l'}$.
On top of that, for each open request $o \in \mathcal{O}_{l}$, the available resource $\Delta d_{\eta_{l},\eta_{l'}}^{o}$ is inherited from the previous extension, which is further restricted by the DTW related to the maximum ride time, denoted by $B_{l'}^{o} - A_{l'}$:
\begin{equation}\label{eq:Delta_d_range}
    \Delta d_{\eta_{l},\eta_{l'}}^{o}= d_{l}^{o}-\max\{B_{l'}^{o}-A_{l'},B_{l'}-A_{l'}\}
\end{equation}

The delay buffer $d_{l}^{o}$ is initially set as the time window length of request $o\in\mathcal{O}_{l}$, and then continuously updated to account for the consumption of actual delay $\delta^{o}_{\eta_{l},\eta_{l'}}$, which is updated as follows:
\begin{equation}
    d^{o}_{l'} \leftarrow d^{o}_{l}-\delta_{\eta_{l},\eta_{l'}}^{o}\label{REF:delay_updates}
\end{equation}

Finally, we define the function $\BFF_{l}(\delta)$ to measure the increment of \emph{maximum} exposed risk considering all open requests $o\in\mathcal{O}_{l}$, which is written as follows:
\begin{equation}
    \BFF_{l}(\delta):= \max\limits_{o\in\mathcal{O}_{l}}\big\{h_{l}^{o}+\BFF_{l}^{o}(\delta)\big\}
\end{equation}
where $\BFF_{l}(\delta)$ takes the upper bounds of piecewise linear functions $h_{l}^{o}+\BFF_{l}^{o}(\delta)$ for each open request $o\in\mathcal{O}_{l}$, and it is also a piecewise linear function. This function will be used to analyze the maximum exposed risk in the case of $\mathcal{O}_{l}^{a}\neq \varnothing$ in the next subsection.
\end{subequations}

\paragraph{\textbf{Risk minimization with associated requests.}}\label{sec:risk_min_associated}
Unlike the previous case, the scenario with non-empty associated requests is more involved since the changes made to open requests can potentially lead to an impact on the already-dropped-off associated requests.
Therefore, we need to examine the MMR condition by balancing the trade-off between the increased exposed risk among associated requests and the reduced exposed risk among open requests.
We define the function $\BFG_{l}(\delta)$ as the maximum risk increase value among associated requests $i\in\mathcal{O}_{l}^{a}$:
\begin{equation}
\BFG_{l}(\delta) := \max_{i\in\mathcal{O}_{l}^{a}} \{h_{l}^{i}+\BFG_{l}^{i}(\delta)\}
\end{equation}
where the function $\BFG_{l}^{i}(\delta)$ is a piecewise linear non-decreasing function constructed based on the sum of co-riders' risk scores and the co-ride time on each traversed arc between the pick-up and drop-off nodes $i$ and $i+n$. Further details and a formal definition of the characteristic function can be found in Appendix~\ref{sec:char_risk_associated}.

Our goal is to find the optimal amount of delay, denoted by $\delta^{\star}$. Based on that, we can obtain the actual delay, calculated as $\delta_{\eta_{l},\eta_{l'}}^{i}=\min\{\Delta d_{\eta_{l},\eta_{l'}}^{i},\delta^{\star}\}$, at the corresponding pick-up node $i$ that satisfied the MMR condition for all requests $i\in\mathcal{O}_{l}\cup\mathcal{O}_{l}^{a}$. The amount of delay $\delta^{\star}$ can be explicitly written as:
\begin{equation}
\delta^{\star} := \arg\min_{\delta}\{\BFF_{l}(\delta),\BFG_{l}(\delta)\}\label{eq:argmin_F_G}
\end{equation}
Since both $\BFF_{l}(\delta)$ and $\BFG_{l}(\delta)$ are piecewise linear functions, the optimal amount of delay $\delta^{\star}$ can be efficiently determined by searching the intersection of two line segments. In certain scenarios, we can make definitive decisions without an exact line-searching procedure, which reduces the complexity. We present three subcases to illustrate these features (see Case 2 in Figure~\ref{fig:risk_calibration_3cases}). 

\begin{enumerate}
    \item \textbf{No Delay:} If $\BFG_{l}(0)\geq \BFF_{l}(0)$, no delay is required. Any further delay would only increase the maximum exposed risk among the associated requests.
    \item \textbf{Max Delay:} If $\BFG_{l}(\max\limits_{o\in\mathcal{O}_{l}}{d_{l}^{o}})\leq \BFF_{l}(\max\limits_{o\in\mathcal{O}_{l}}{d_{l}^{o}})$, one can offset as much waiting time as possible by setting $\delta= \max\limits_{o\in\mathcal{O}_{l}}{d_{l}^{o}}$. In this case, the maximum exposed risk among the associated requests will not exceed the maximum exposed risk among the open requests.
    \item \textbf{Partial Delay:} In this scenario, we strategically determine an optimal solution by balancing the maximum exposed risk between the open and associated requests. The goal is to minimize the overall risk while satisfying the MMR condition. To achieve this, we decompose the problem in Eq.~\eqref{eq:argmin_F_G} into finding the intersection of two line segments, which can be effectively solved using a line-searching algorithm. Detailed examples are provided in Example~\ref{example:piecewise}, and further methodological details can be found in Appendix~\ref{sec:char_risk_associated}.
\end{enumerate}

\begin{example}\label{example:piecewise}
In Figure~\ref{fig:lp_example}, we consider an interlaced segment $(i, j, k, o, i+n)$ extended from node $i+n$ to a new request $m \in \mathcal{P}$, where a waiting time $\Delta \omega_{i+n, m}>0$ needs to be offset. 
\begin{figure}[H]
    \caption{Illustration of interlaced segment $(i,j,k,o,i+n,\cdots)$ and the piecewise functions $\BFF_{l}(\delta)$ and $\BFG_{l}(\delta)$.}
    \label{fig:lp_example}
    \begin{tikzpicture}
        % Axis
        \draw[->] (-1,0) -- (5,0) node[right] {$\delta$};
        \draw[->] (0,-1) -- (0,4) node[above] {$h$};
        
        \draw[] (0,1) -- (1,1) -- (2,1) -- (3,1.75) -- (4.5,3.75);
        \draw[dashed] (0,3.75) -. (4.5,2.5) node[midway,above=3mm]{};%${$h_{l}^{i}-\left(R_{l}-r_{i}\right)\delta$};
    
        % Vertical lines and labels
        \draw[dotted] (1,0) -- (1,1) node[right]{};% {$G_{l}^{o}(\Delta d_{l}^{i})$};
        \draw[dotted] (2,0) -- (2,1) node[right]{};% {$G_{l}^{o}(\Delta d_{l}^{j})$};
        \draw[dotted] (3,0) -- (3,1.75) node[right]{};% {$G_{l}^{o}(\Delta d_{l}^{k})$};
        \draw[dotted] (4.5,0) -- (4.5,3.75) node[right]{};% {$G_{l}^{o}(\Delta d_{l}^{k})$};

        % # compare M^a and M
        \node[draw,inner sep=2.5pt,label={right:$\BFG_{l}(\Delta d_{\eta_{l},\eta_{l'}}^{o})$}] at (4.5,3.75) {};
        \node[draw,inner sep=2.5pt,label={right:$\BFF_{l}(\Delta d_{\eta_{l},\eta_{l'}}^{o})$}] at (4.5,2.5) {};
        \node[draw,circle,inner sep=2pt] (OptimalPt) at (3.675,2.72) {};
        \node[above=2mm] at (OptimalPt.west) {$(\delta^{*},h^{*})$};
        
        % draw its location
        \draw[dash dot dot] (3.675,0) -- (3.675,2.72) node[right]{};%{$\delta^{\star}$};% {$G_{l}^{o}(\Delta d_{l}^{k})$};

        \draw[dash dot dot] (0,2.72) -- (3.675,2.72) node[left]{};%{$h^{\star}$};% {$G_{l}^{o}(\Delta d_{l}^{k})$};

        \node[left] at (0,2.72) {$h^{\star}$};
        \node[below] at (3.675,0) {$\delta^{\star}$};
        
        % \node[left] at (0,0) {$0$};
        \node[left] at (0,1) {$\BFG_{l}(0)$};
        % \node[left] at (0,2) {$2$};
        \node[left] at (0,3.75) {$\BFF_{l}(0)$};
        
        \node[below] at (1,0) {$\Delta d_{\eta_{l},\eta_{l'}}^{i}$};
        \node[below] at (2,0) {$\Delta d_{\eta_{l},\eta_{l'}}^{j}$};
        \node[below] at (3,0) {$\Delta d_{\eta_{l},\eta_{l'}}^{k}$};
        \node[below] at (4.5,0) {$\Delta d_{\eta_{l},\eta_{l'}}^{o}$};    
        % add scope:
        \draw[] (0.5,1) node[above] {$0$};
        \draw[] (1.5,1) node[above] {$0$};
        
        \draw[] (2.15,1.6) node[right] {$r_j$};
        \draw[] (2.15,2.35) node[right] {$r_j+r_k$};

        \draw[] (0.5,0.5) node[right] {$\BFG_{l}(\delta)$};
        \draw[dashed] (0.5,4) node[right] {$\BFF_{l}(\delta)$};
    \end{tikzpicture}
    \includegraphics[width=0.55\textwidth]{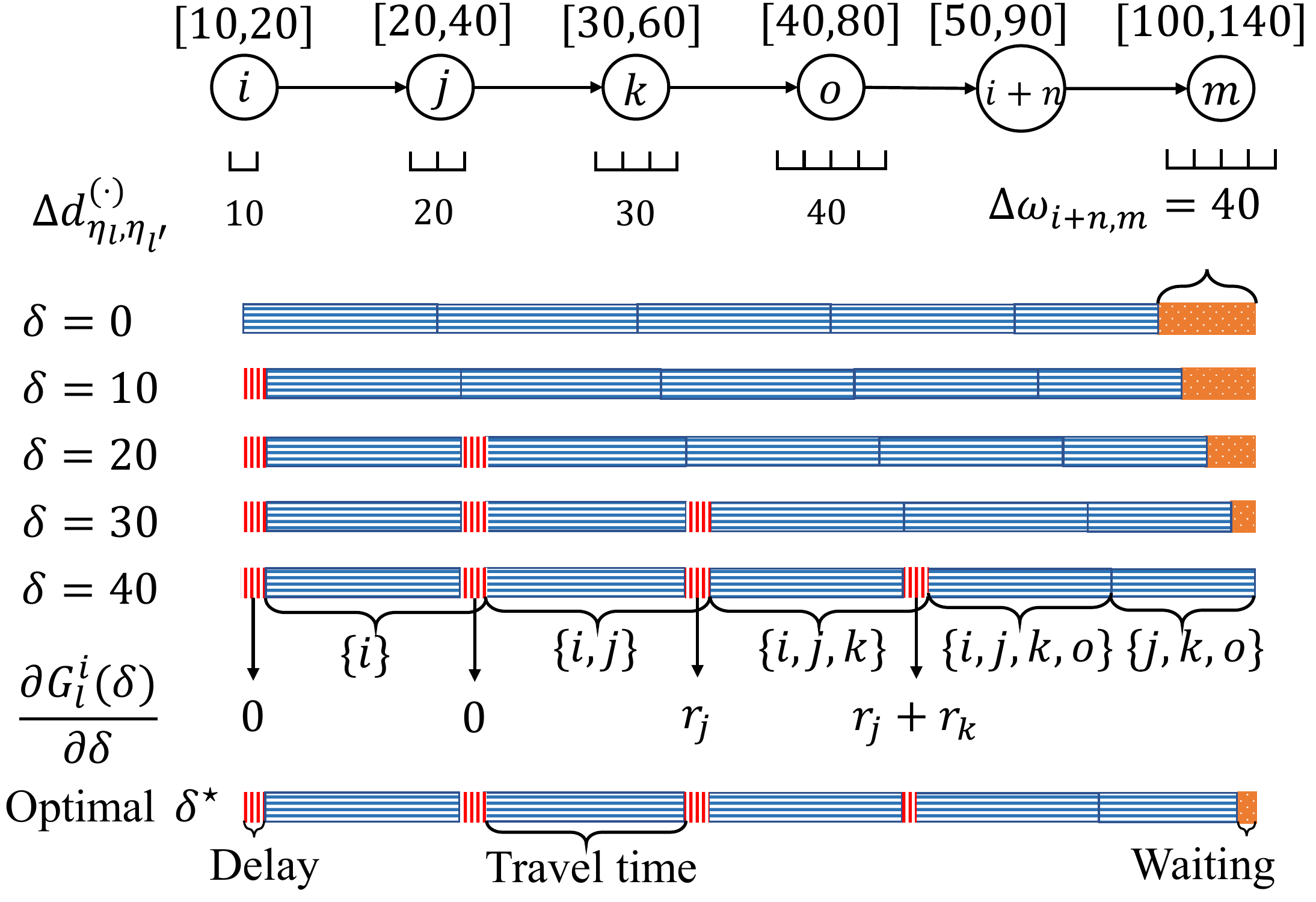}
\end{figure}
Assuming $\Delta d_{\eta_{l},\eta_{l'}}^{i} \leq \Delta d_{\eta_{l},\eta_{l'}}^{j} \leq \Delta d_{\eta_{l},\eta_{l'}}^{k} \leq \Delta d_{\eta_{l},\eta_{l'}}^{o}$, if we only focus on minimizing the exposed risk for open requests, we can delay the request $o$ by $\delta = \Delta d_{\eta_{l},\eta_{l'}}^{o} = 40$ to fully offset the waiting time of $\Delta \omega_{i+n, m} = 40$. This adjustment results in four distinct delay segments of length 10 at nodes $i$, $j$, $k$, and $o$, as shown in the scenario $\delta=40$ in Figure~\ref{fig:lp_example}. Consequently, the exposed risk for request $o$ is reduced from $\BFF_{l}(0)$ to $\BFF_{l}(\Delta d_{\eta_{l},\eta_{l'}}^{o})$.

On the other hand, the characteristic function $\BFG_{l}(\delta)$ consists of four line segments representing the request $i$'s traversed arcs $(i,j)$, $(j,k)$, $(k,o)$, and $(o,i+n)$. Based on this, we update the exposed risk for associated request $i$, such that $\BFG_{l}(\Delta d_{\eta_{l},\eta_{l'}}^{o})>\BFF_{l}(\Delta d_{\eta_{l},\eta_{l'}}^{o})$. However, we observe that $\BFG_{l}(0)<\BFF_{l}(0)$ and $\BFG_{l}(\Delta d_{\eta_{l},\eta_{l'}}^{o})>\BFF_{l}(\Delta d_{\eta_{l},\eta_{l'}}^{o})$, indicating that delaying open request $o$ to the maximum extent leads to a higher maximum risk of associated request $i$. Therefore, we aim to find the amount of delay $\delta^{\star}:=\argmin\limits_{\delta}\{\BFF_{l}(\delta),\BFG_{l}(\delta)\}$ that partially offsets the waiting time while satisfying the MMR condition. This solution is represented by $(\delta^{\star}, h^{\star})$, which achieves a compromise between minimizing the risk increments and meeting the MMR condition.
\end{example}

\paragraph{\textbf{Risk Calibration Algorithm.}}\label{sec:risk_minimization}
Having enumerated all the possible scenarios in risk calibration, we will now present the comprehensive procedure of the risk calibration in Algorithm~\ref{alg:contact_duration_calibriation}. The goal of risk calibration algorithm is to calculate the amount of delay $\delta_{\eta_{l},\eta_{l'}}^{o}$ that satisfies the MMR condition for each open and associated request $o\in\mathcal{O}_{l}\cup\mathcal{O}_{l}^{a}$. 
\begin{algorithm}[H]
\footnotesize%\small, \scriptsize, or \tiny
\caption{Risk Calibration and Optimization Algorithm}
\label{alg:contact_duration_calibriation}
Label $l$, node-to-extend $\eta_{l'}$, $\mathcal{O}^{a}_{l}$, $\mathcal{O}_{l}$, $A_{l},A_{l'},B_{l'},\{B_{l'}^{o}:o\in\mathcal{O}_{l}\}, \{r_i,d_{l}^{i}, h_{l}^{i}:i\in\mathcal{O}_{l}\cup\mathcal{O}_{l}^{a}\}$ \Comment{Input}\\
Onboard travel time $\Delta A_{\eta_{l},\eta_{l'}}^{o}, \forall o\in\mathcal{O}_{l}$; actual amount of delay at pick-up node $i$, $\delta_{\eta_{l},\eta_{l}'}^{i}, \forall i\in\mathcal{O}_{l}\cup\mathcal{O}_{l}^{a}$ \Comment{Output}
\begin{algorithmic}[1] %[1] enables line numbers
    % \State \Comment{Step 1: Min-Risk Calculation (Maximally delay the open requests)}
    \If{$\mathcal{O}_{l}^{a}\not=\varnothing$}\Comment{Step 1: Min-Risk Calculation for open requests}
    \State \textbf{Calculate} the potential time to wait $\Delta \omega_{\eta_{l},\eta_{l'}}$ following Eq.~\eqref{eq:waiting_time_algo2}
    \State \textbf{Initiate} the maximum amount of delay at node $\eta_{l'}$: $\Delta d_{\eta_{l},\eta_{l'}}^{\eta_{l'}}=B_{l'}-A_{\eta_{l'}}, \text{if~}\eta_{l'}\in\mathcal{P}$
    \State \textbf{Calculate} the remaining amount of delay $\Delta d_{\eta_{l},\eta_{l'}}^{o}, \forall o\in\mathcal{O}_{l}$ following Eq.~\eqref{eq:Delta_d_range}
    % \State \textbf{Calculate} minimal onboard duration: $\Delta A_{\eta_{l},\eta_{l'}}^{o}, \forall o\in\mathcal{O}_{l}$, where
    \State \textbf{Calculate} the optimal amount of delay:  $\delta_{\eta_{l},\eta_{l'}}^{o}$ following Eq.~\eqref{eq:delta_update_algo2}
    % \State */\textbf{Step 2: MMR Calibration (Examine the MMR condition if the associated request exists)}*/
    \ElsIf{$\mathcal{O}_{l}^{a}\not=\varnothing$}\Comment{Step 2: MMR Calibration if associated requests exist}
        \State \textbf{Calculate} the interval of maximum exposed risk of open requests: $\BFF_{l}(0)$ and $\BFF_{l}(\max\limits_{o\in\mathcal{O}_{l}}{\Delta d_{\eta_{l},\eta_{l'}}^{o}})$
        \State \textbf{Calculate} the interval of maximum exposed risk of associated requests: $\BFG_{l}(0)$ and $\BFG_{l}(\max\limits_{o\in\mathcal{O}_{l}}{\Delta d_{\eta_{l},\eta_{l'}}^{o}})$
        \If{$\BFF_{l}(0)\leq\BFG_{l}(0)$}\Comment{No Delay}
            \State \textbf{Obtain} the amount of delay: $\delta_{\eta_{l},\eta_{l'}}^{o}=0,\forall o\in\mathcal{O}_{l}\cup\mathcal{O}_{l}^{a}$
        \ElsIf{$\BFF_{l}(\max\limits_{o\in\mathcal{O}_{l}}{\Delta d_{\eta_{l},\eta_{l'}}^{o}}))\geq\BFG_{l}(\max\limits_{o\in\mathcal{O}_{l}}{\Delta d_{\eta_{l},\eta_{l'}}^{o}}))$}\Comment{Max Delay}
            \State \textbf{Calculate} the amount of delay: $\delta_{\eta_{l},\eta_{l'}}^{o}, \forall o\in\mathcal{O}_{l}\cup\mathcal{O}_{l}^{a}$ following Eq.~\eqref{eq:delta_update_algo2}
        \Else \Comment{Partially Delay to meet MMR condition}
            \State \textbf{Obtain} the optimal amount of delay $\delta^{\star}$ and \textbf{output} the amount of delay $\delta_{\eta_{l},\eta_{l'}}^{o}=\min\{\Delta d_{\eta_{l},\eta_{l'}}^{o},\delta^{\star}\}$ 
        \EndIf
    \EndIf    
    \State \textbf{Update} the remaining amount of delay: $ d^{o}_{l'}$ following the REF in Eq.~\eqref{REF:delay_updates}
    \State \textbf{Output} the onboard duration that meets the MMR condition:  $\Delta A_{\eta_{l},\eta_{l'}}^{o} \leftarrow t_{\eta_{l},\eta_{l'}}+s_{\eta_l} +\Delta \omega_{\eta_{l},\eta_{l'}}-\delta_{\eta_{l},\eta_{l'}}^{o}$
\end{algorithmic}
\end{algorithm}

We summarize the updating rules in Algorithm~\ref{alg:contact_duration_calibriation} as follows:
\begin{itemize}
    \item \textbf{Step 1:} Calculate the minimal onboard travel time for all open requests $o \in \mathcal{O}_{l}$.
    \item \textbf{Step 2:} If there is a non-empty set of associated requests, calculate the increase in exposed risk for the associated requests and examine if the MMR condition still holds.
    \item \textbf{Updating Resources:} Based on the utilization of the resources, we update the remaining adjustable amount of delay $ d_{l'}^{o}$ and further calibrate the onboard duration $\Delta A_{\eta_{l},\eta_{l'}}^{o}$ based on the optimal amount of delay $\delta_{\eta_{l},\eta_{l'}}^{o}$ in lines 18--19.
\end{itemize}
In Algorithm~\ref{alg:contact_duration_calibriation}, the time complexity is bounded by $O(W_{\max}|\mathcal{N}|)$, where $W_{\max}$ is the maximum number of open and associated requests that need to be examined (line 4), and $|\mathcal{N}|$ is the total number of nodes. This indicates that the contact duration calibration is performed on the Hamiltonian path with precedence constraints. However, in practice, it is unlikely to consider the longest route due to the high travel cost and increased exposed risk. Therefore, we can conclude that our algorithm has a polynomial time complexity, which is more efficient compared to solving non-linear constraints (e.g., Constraints~\eqref{eq:risk_consistency}) in the arc-based formulation.

We now present the REF for risk-related resources. The set of associated requests $\mathcal{O}_{l'}^{a}$ is initially empty. It is reset to empty when all open requests have been serviced. When dropping off a request $\eta_{l'}\in\mathcal{D}$, check if it has been exposed to any open requests $o\in\mathcal{O}_{l}$. If it has been exposed, add the serviced request $\eta_{l'}-n$ to the set $\mathcal{O}_{l}^{a}$.
\begin{subequations}\label{eq:ref_rdarp}
\begin{equation}
    \mathcal{O}_{l'}^{a} \leftarrow
    \begin{cases}
        \mathcal{O}_{l}^{a}\cup\{\eta_{l'}-n\} & \text{if}~\eta_{l'}\in \mathcal{D} \text{~and~} \min\limits_{o\in\mathcal{O}_{l}}\{\Gamma_{l}^{o}\}<\Gamma_{l}^{\eta_{l'}}\\
        \varnothing & \text{if}~\eta_{l'}\in \mathcal{D} \text{~and~} \mathcal{O}_{l}=\{\eta_{l'}-n\}\\
    \end{cases}\label{ref:associated_request_set}
\end{equation}

Next, we update the total risk scores of onboard requests and the cumulative exposed risk at node $\eta_{l'}$ as $R_{l'}$ and $Q_{l'}$ respectively. To ensure the MMR condition, we calibrate the onboard travel time $\Delta A_{\eta_{l},\eta_{l'}}^{o}$ and amount of delay $\delta_{\eta_{l},\eta_{l'}}^{o}$ for each open request $o\in\mathcal{O}_{l}$ on arc $(\eta_{l},\eta_{l'})$. The increment of the cumulative exposed risk of open requests is calculated as the sum of the risk score of request $o$, $r_{o}$, and the onboard travel time  $\Delta A_{\eta_{l},\eta_{l'}}^{o}$. Notably, during the calibration procedure in Algorithm~\ref{alg:contact_duration_calibriation}, there may exist a potential increased risk among associated requests due to the adjustment on the open requests. With slight abuse of notation, $\delta_{\eta_{l},\eta_{l'}}^{i}$ for associated requests $i\in\mathcal{O}_{l}^{a}$ represents the \emph{equivalent amount of delay}. Therefore, the cumulative exposed risk is calculated based on the equivalent delay and the risk scores of the associated requests.
    \begin{align}
        R_{l'} &\leftarrow R_{l} + r_{\eta_{l'}} \label{REF:total_risk_score}\\
        Q_{l'} &\leftarrow Q_{l} + \sum_{o\in\mathcal{O}_{l}} r_{o}\cdot \Delta A_{\eta_{l},\eta_{l'}}^{o} + \sum_{i\in\mathcal{O}_{l}^{a}} r_{i}\cdot \delta_{\eta_{l},\eta_{l'}}^{i}\label{REF:cumulative_risk}
    \end{align}

We update the exposed risk for each open request $o\in\mathcal{O}_{l}$ based on the service sequence. Co-riders $i$ with $\Gamma_{l}^{i}<\Gamma_{l}^{o}$ (who were visited before picking up request $o$) will have the same onboard time as request $o$, while co-riders $j$ visited after request $o$, i.e., $\Gamma_{l}^{j}>\Gamma_{l}^{o}$, may have a shorter co-ride time with request $o$ due to a longer delay at the subsequent pick-up nodes $j$.
\begin{equation}
    h_{l'}^{o} \leftarrow 
    h_{l}^{o} + \sum_{i\in\mathcal{O}_{l},\Gamma_{l}^{i}<\Gamma_{l}^{o}} r_{i} \Delta A_{\eta_{l},\eta_{l'}}^{o} + \sum_{j\in\mathcal{O}_{l},\Gamma_{l}^{j}>\Gamma_{l}^{o}} r_{j} \Delta A_{\eta_{l},\eta_{l'}}^{j}, \forall o\in\mathcal{O}_{l} \label{ref:individual_risk}
\end{equation}

Next, we update the potential increased exposed risk for each associated request $i\in\mathcal{O}_{l}^{a}$. 
\begin{equation}
    h_{l'}^{i} \leftarrow h_{l}^{i} + \BFG^{i}_{l}\left({\delta}_{\eta_{l},\eta_{l'}}^{i}\right), \forall i\in\mathcal{O}_{l}^{a}\label{ref:individual_associated_risk}
\end{equation}

Finally, we introduce the REF for the reduced cost, which involves two terms on arc $(\eta_{l},\eta_{l'})$: the reduced cost $\Tilde{c}_{\eta_{l},\eta_{l'}}$ and the sum of arc-level increment of individual exposed risk, weighted by the associated dual variables $\rho_{o}$. The REF can be written as follows:
\begin{equation}
    \Tilde{c}_{l'} \leftarrow \Tilde{c}_{l} + \Tilde{c}_{\eta_{l},\eta_{l'}} - \sum_{o\in\mathcal{O}_{l}\cup\mathcal{O}_{l}^{a}}\rho_{o}\left(h_{l'}^{o}-h_{l}^{o}\right)\label{REF:reduced_cost}
\end{equation}

\end{subequations}

\begin{remark}\label{remark:equivalency}
    The proposed risk calibration algorithm goes beyond the condition in Proposition~\ref{prop:equivalence} by not only minimizing the maximum exposed risk for the entire trip but also ensuring that the risk increments at each traversed arc are minimized. This also results in one MMR solution that consumes the minimum amount of delay buffer. MMR solutions are in general not unique and other possible solutions can be identified by further delaying at certain nodes if buffer permits. 
\end{remark}

\subsubsection{Dominance rule for RDARP}\label{subsec:dominance_darp}
In the label extension procedure, the number of labels will increase exponentially. To avoid exhaustive enumeration, we adopt the dominance rule, which identifies the unpromising labels during the label extension and discards them. For simplicity, we differentiate the associated variables of labels $l_{1}$ and $l_{2}$ by using the subscripts $1$ and $2$.

\begin{proposition}\label{prop:dominance_rule_darp}
Given two RDARP-feasible labels $l_{1}$ and $l_{2}$ residing at the same node, i.e., $\eta_{1}=\eta_{2}\in\mathcal{N}$, $l_{1}$ is said to dominate $l_{2}$ if the following conditions are satisfied:
\begin{subequations}\label{eq:dominance_rule}
    \begin{align}
        &A_{1}\leq A_{2}, W_{1}\leq W_{2}, {Q_{1}\leq Q_{2}},\mathcal{V}_{1} \subseteq \mathcal{V}_{2}~,\mathcal{O}_{1}\subseteq \mathcal{O}_{2}, \mathcal{O}_{1}^{a} \subseteq \mathcal{O}_{2}^{a}~\label{dom:dominance_rule_darp0}\\
        &\Tilde{c}_{1} \leq \Tilde{c}_{2}~\label{dom:dominance_rule_darp1}\\
        &D^{o}_{1}(A_{1}) - A_{1}\geq D^{o}_{2}(A_{2}) - A_{2}, &&\forall o \in \mathcal{O}_{1}~\label{dom:dominance_rule_darp2}\\ 
        &D^{o}_{1}(B^{o}_{1}) \geq D^{o}_{2}(B^{o}_{2}), &&\forall o \in \mathcal{O}_{1}~\label{dom:dominance_rule_darp3}\\
        & d_{1}^{o} \geq d_{2}^{o}, &&\forall o \in \mathcal{O}_{1}\cup\mathcal{O}_{1}^{a}~\label{dom:dominance_rule_darp4}
    \end{align}
\end{subequations}
\end{proposition}

\proof{Proof.} 
See Appendix~\ref{proof:dominance_rule_darp}.
\endproof
We note that the delivery triangle inequality (DTI) holds in the reduced cost matrix for the components associated with the reduced cost $\sum\limits_{(i,j)\in r}\Tilde{c}_{ij}x_{ij}$ and individual risk $\sum\limits_{i\in\mathcal{P}}H_{ir}\rho_{r}$ in Eq.~\eqref{eq:reduced_cost_objective_funtion}. While the former can be transformed into a reduced cost matrix that satisfies the DTI following \cite{ropke2009branch}, the latter cannot be determined without a given (partial) route. However, inserting a new delivery node between two nodes will not result in a lower reduced cost. In this case, there exists no cycle of negative reduced cost, which implies that the DTI always holds.
Therefore, the dominance rule in Eq.~\eqref{dom:dominance_rule_darp0} inspects the subset relationship of the open requests. 

\section{Branch-cut-and-price algorithm}\label{sec:branch_and_price}
This section describes the BCP algorithm to solve the RDARP. The idea of adopting cutting planes in the BP algorithm has been successfully applied to solve large-sized instances of the DARP~\citep{gschwind2015effective,luo2019two} and PDPTW~\citep{ropke2009branch,qu2015branch}.
Specifically, the cutting planes are generated by conducting the separation heuristics before exploring the root node of the branching tree. Therefore, the procedure of the column and cut generations in the BB framework results in a BCP algorithm.

\subsection{Valid inequalities}
Three types of valid inequalities are considered in our study, including the IPEC, two-path constraints, and rounded capacity (RC) constraints. The separation heuristics are adopted to identify the violated inequalities based on a node set, which is implemented using constructive heuristic methods {(for the IPEC, see Section 5.3.4 in~\cite{cordeau2006branch} and for the two-path and RC constraints, see Section 6.2 in~\cite{ropke2009branch})}.
The formula of valid inequalities and the approach to incorporate the associated dual variables into the RMP are introduced in Appendix~\ref{appendix:dual_valid_inequalities}. An experiment is conducted to investigate the impact of valid inequalities on the LP relaxation of the trip-based formulation. Details are presented in Appendix~\ref{appendix:impact_of_valid_inequality}

\subsection{Branching strategy}

The BCP algorithm uses two types of branching strategies without changes in the problem structure: first on the number of vehicles, then on the outflow of a node pair. Linear constraints are added to the RMP, and their dual variables are included in the reduced cost function. A best-first search is used to explore the branching tree. See Appendix~\ref{appendix:branching_cuts} for more details.

\subsection{Acceleration technique}
To speed up the CG procedure, we consider a heuristic relaxation of the labeling algorithm by weakening the dominance rule. Specifically, for the weakened dominance rule, we relax the comparison between the visited nodes of two labels, which results in more discarded labels. Thus, several negative reduced cost columns can be found rapidly. However, this procedure cannot guarantee that all columns with negative reduced costs are generated, yielding a heuristic solution. In light of this, when no more columns can be generated under the weakened dominance rule, we reiterate the labeling procedure, this time with the complete dominance rule to ensure optimality of the solution \citep{dayarian2015branch}. 
\section{Model extension - Equitable DARP}\label{sec:edarp}
This section presents an extension of our RDARP that can ensure an equitable distribution of ride times among the passengers, a variant called Equitable DARP (EDARP). In the EDARP, the min-max objective is replaced by minimizing the maximum detour rate, defined by the ratio of onboard duration $A_{n+i}-A_{i}$ (variable) and direct trip time $t_{i,n+i}$ (parameter). We note that the EDARP can be solved by the same BCP framework. However, several modifications are required to adapt the min-max objective of the detour rate in the EDARP as detailed next.

\subsection{Modification on the Problem Setting}

In the EDARP, the onboard duration for passenger $i$ can be understood as the contact duration between passenger $i$ and a dummy passenger traveling from origin $0$ to destination depot $2n+1$. Hence, we make the following changes to effectively solve the EDARP without structural changes.

\begin{itemize}
    \item We insert a zero-volume dummy passenger $p_{0}$, who travels from $p_{0}^{o}=0$ to $p_{0}^{d}=2n+1$ and can be visited multiple times. The time windows of $p_{0}^{o}$ and $p_{0}^{d}$ are the same as the depots.
    \item The only successor of depot $0$ is modified to be $p_{0}^{o}$, and the only predecessor of depot $2n+1$ is $p_{0}^{d}$. Hence, the vehicle must visit $p_{0}^{o}$ first, serve real passengers $i\in\mathcal{P}$, and end with $p_{0}^{d}$.
    \item We neglect the exposed risk among the real passengers by setting $r_{0}=1$ for the dummy passenger and $r_{i}=0$ for all real passengers $ i\in\mathcal{P}$. In this case, the exposed risk can be understood to be mono-directional from the dummy passenger to real passengers with the rate of $r_{0}=1$. And the real passengers are free of getting infected by all passengers except the dummy one.
\end{itemize}

To enhance the understanding on our EDARP, we present an illustrative example in Appendix~\ref{appendix:detour_rate_calculation}. 
Based on the EDARP setting, we show in the following proposition that, the measure of $H_i$ in RDARP and onboard duration in EDARP are equivalent. Hence, the onboard duration in EDARP is equivalently ensured to be minimized along a given route.
\begin{proposition}
    Under the EDARP setting, the measure of individual exposed risk $H_{i}$ for a real passenger $i\in\mathcal{P}$ is equivalent to the onboard duration $A_{n+i}-A_{i}$.\label{prop:rdarp_edarp}
\end{proposition}
\proof{Proof.}
See Appendix~\ref{proof:rdarp_edarp}.
\endproof

Following the notations above, the maximum detour rate among all real passengers $i\in\mathcal{P}$ can be expressed as $\max\limits_{i\in\mathcal{P}}\{\frac{A_{n+i}-A_{i}}{t_{i,n+i}}\}=\max\limits_{i\in\mathcal{P}}\{\frac{H_{i}}{t_{i,n+i}}\}$. Let $\epsilon^{\text{dt}}$ be the upper bound of the maximum detour rate. We impose the min-max objective of detour rate with the upper bounds of $\epsilon^{\text{dt}}$ (parameter) or $\overline{D}$ (variable), which takes similar forms as in Constraints~\eqref{cons:max_risk_dual_rho_node_level}:
\begin{equation}
    \sum_{r\in\Omega'} \frac{H_{ir}}{t_{i,n+i}}\alpha_{ir}\lambda_{r} \le \epsilon^{\text{dt}} ~\text{or}~ \sum_{r\in\Omega'} \frac{H_{ir}}{t_{i,n+i}}\alpha_{ir}\lambda_{r} \le \overline{D}, \forall i \in \mathcal{P}
\end{equation}
 % we describe more general problems as a simplified version of the RDARP.
\section{Computational Study}~\label{sec:results}
This section introduces the test instances for the RDARP, including (1) the transmission risk extension of existing DARP benchmark instances and (2) the instances generated based on real-world paratransit trip data. 
We next conduct extensive experiments using our BCP algorithm and present the results. Further, detailed Pareto fronts together with the sensitivity analyses are presented to draw insights into the trade-off between travel cost and exposed risk mitigation.

The BCP algorithm was coded via Python 3.8, where the RMP is solved by Gurobi 9.5. The labeling algorithms were implemented in Julia Programming Language 1.7.1 as a callable function in Python. All experiments were performed on a server with a 24-core AMD Ryzen Threadripper 3960X CPU, 128 GB RAM, and the Linux operating system.
\subsection{Test Instances}
\subsubsection{Benchmark Instances}
We first adopt the DARP benchmark instances provided by \cite{cordeau2006branch} \citep[later augmented by][]{ropke2007models}. {For the benchmark DARP instances, there are type ``a" and ``b" instances, indicating the ride-sharing and minibus scenarios. As the cases of the minibus are highly analogous to our real-world paratransit instances, we will only consider the type-a instances in our experiment.} The type-a instances are originally labeled as ``$aK\text-N$", indicating up to $K$ vehicles and $N$ passengers. There is one passenger per node, the vehicle capacity is $3$, the time window size is $15$, and the maximum ride time is $30$ \citep[for the characteristics of DARP instances, see][]{ropke2007models}.
Next, to incorporate the risk components for RDARP, we {assume that risk score $r_{i}$ per node $i\in\mathcal{N}$ is proportional to the number of passengers, such that $r_{i}=w_{i}$.}
We perform a similar pre-processing procedure following {\cite{cordeau2006branch}}, including arc elimination and time window tightening in order to reduce the problem size by removing the infeasible arcs from the initial complete graph $\mathcal{G}$.

\subsubsection{Real-world Paratransit Trip Data}
We use the paratransit trip log data in the central Alabama area covering Jefferson County, Shelby County, and Walker County. The paratransit service is operated by Central Alabama's Specialized Transit (ClasTran), a non-profit corporation that provides specialized mobility services to people with disabilities and limited mobility, seniors, and people who reside in and travel to or from rural areas~\citep{clastran2021}. Interested readers are referred to our previous study~\citep{nie2022impact} for more details on the historical trip data.
The trip records provided by ClasTran include the exact information of the passengers (i.e., age, travel purpose) and trips (i.e., trip date, coordinates of the origin and destination, and the scheduled and reported pick-up and drop-off time). We particularly select one-day trip data on March 11, 2020, before the national emergency regarding the COVID-19 pandemic was announced. At that time, the demand remained at a high level (336 versus roughly 100 daily trips after the national emergency). {In particular, those 336} passengers were originally assigned to 66 different trips without the consideration of risk mitigation. Hence, it makes considerable space for us to optimize the paratransit scheduling while ensuring a low level of exposed risk on the route.

Specifically, we generate 20 different instances by considering 10 different problem sizes from one-day trip data during the morning {(AM)} and afternoon {(PM)} periods. The problem size ranges from 17 passengers with 3 vehicles to 55 passengers with 13 vehicles.{For example, instance ``$AM\text-N\text-K$" indicates the instance during the AM period, where $N$ passengers are served by up to $K$ vehicles.} 
% Note that the selected passengers are observed to conduct two trips during the morning and afternoon periods. 
We further categorize those trips to be \textit{outbound} and \textit{inbound} trips, indicating the trips from home locations to destinations and the return trip, respectively. 
Following \cite{cordeau2006branch}, we assume that the passengers either reserve a desired arrival time at the destinations (e.g., the outbound trips) or a desired departure time at the origins (e.g., inbound trips).
In this case, we conduct the same preprocessing procedures as described in \cite{cordeau2006branch}. In addition, the width of the time window is set to be 30 minutes according to the policy~\citep{clastran2021}, the maximum ride time is set based on the historical trip duration $\Bar{L}_{i}$ for each passenger $i\in\mathcal{P}$ such that $L_{\max}^{i}=1.25 \Bar{L}_{i}$. 
To better represent the practical operation, we consider heterogeneous passengers with different risk scores ($r_{i}$). The risk score is assumed to range from 0 (lowest) to 0.9 (highest), which is estimated by three equally-weighted factors: the county of origin, travel purpose, and age group (for details, see Appendix~\ref{appendix:risk_assessment}). As observed in the historical routes, there may exist multiple passengers associated with one location $w_{i}>1, \forall i \in\mathcal{P}$. We will consider the risk score as the sum of all associated passengers and neglect the risk transmission within the group traveling together. This is because those passengers traveling together as a group are likely to spend time together longer than the onboard ride time. It is also worth noting that the risk score only serves for the comparison purpose, where a higher risk score of a passenger is an indication of greater vulnerability to disease compared to the other co-riders. 
Further, the identical maximum cumulative risk for each vehicle $Q_{\max}$ is set to be $(b_{2n+1}-a_{0})\times\frac{1}{n}\sum\limits_{i\in\mathcal{P}}r_{i}$ by considering the trip duration and average risk score.
Since the risk score $r_{i}$ is considered to be a constant, the units of $Q_{\max}$, $Q_{i}$, and $H_{i}$ are ``minutes". 
Finally, all real-world instances are available at \url{https://github.com/sguo28/rdarp_instance}.
\subsection{Experimental Setting}
This section presents the results of the numerical experiments. 
The BCP algorithm is first applied to the extended DARP benchmark instances, and then the real-world instances. The key performance metrics  are listed below.
\begin{itemize}
    \item Cost: the travel cost in the best integer solution (unit: minutes);
    \item $\overline{H}$: maximum individual risk in the best integer solution (unit: minutes);
    \item $T_{MP}$: Computational time for the RMP (unit: seconds);
    \item $T_{SP}$: Computational time for the SP (unit: seconds).
    \item{DARP(\%), Hist(\%): the percentage change to the baseline (e.g., DARP or historical routes), calculated as $(z^{*}-\Bar{z})/\Bar{z}$, where $z^{*}$ and $\Bar{z}$ denote the selected metrics (e.g., Cost, $\overline{H}$, $T_{MP}$, or $T_{SP}$) of current model and the baseline, respectively.}
\end{itemize}

\subsection{Experiments on Extended DARP Benchmark Instance}
This section presents the computational results of DARP and RDARP on the extended benchmark instances. 
First, the DARP case ($\epsilon^{\text{risk}}\rightarrow\infty$) serves as a baseline, which helps to understand the potential improvement space for risk mitigation. We can also confirm the validity of our approach with respect to the DARP by comparing our results with the best-known solutions in the literature. Further, we will present the results of RDARP ($\epsilon^{\text{risk}}=15$ and $\epsilon^{\text{risk}}=30$) to explore the impacts of risk factor on the routing strategy and total travel cost. The case of $\epsilon^{\text{risk}}=15$ is closely related to the definition of close contact~\citep{cdc2022appendices}, which indicates the one rider's contact within 6 feet for a total of 15 minutes or more. 
The case of $\epsilon^{\text{risk}}=30$ suggests a less restrictive version of RDARP, which helps to understand the trade-off between the exposed risk and total travel cost.

In Table~\ref{table:type_a_benchmark_V1}, all 21 DARP instances can be solved within one hour CPU time, and the optimal integer solutions regarding the travel cost are consistent with the best-known solutions as reported in \cite{gschwind2015effective}. As a by-product of our BCP algorithm, we can also obtain the maximum exposed risk $\overline{H}$ in the DARP solutions. With the minimum total travel cost, it is reported that only 6 of 21 instances are associated with a maximum exposed risk of less than $30$ min, and none of the instances are found with an exposed risk of less than $15$ min. Such a high level of risk exposure in the DARP solutions indicates a significant potential for implementing a risk-aware routing strategy. Specifically, for the RDARP problem, we observed that 17 out of 21 instances under $\epsilon^{\text{risk}}=30$ and $\epsilon^{\text{risk}}=15$ were solved optimally within the given CPU time limit. There was one instance ("a3-30") under $\epsilon^{\text{risk}}=15$ that could not be solved due to the feasibility issue under the relatively strict $\epsilon^{\text{risk}}$ value. However, we can still obtain meaningful IP solutions in all other 41 feasible instances under $\epsilon^{\text{risk}}=30$ and $\epsilon^{\text{risk}}=15$. Hence, we report the high efficiency of identifying an applicable risk-aware routing strategy with joint consideration of travel cost and exposed risk. 

When comparing the average performance of CPU times for RDARP, we observed a significant reduction in CPU time for the cases with $\epsilon^{\text{risk}}=15$, particularly in the instances of ``a3-36", ``a4-48", ``a5-60", and ``a8-80". This is due to the fact that, with a more restrictive risk threshold of $\epsilon^{\text{risk}}=15$, the feasible region for the labeling approach solving the ESPPRCMMR becomes smaller, leading to a potentially smaller number of feasible label extensions.
Observing the trade-off in the RDARP, we report that one does not necessarily need to compromise the travel cost heavily to improve the maximum exposed risk. Specifically, 1.08\% and 6.86\% more travel costs can lead to a reduction of exposed risk of 44.40\% and 55.81\% in the cases of $\epsilon^{\text{risk}}=30$ and $15$, respectively (for details, see the last row ``DARP(\%)" in Table~\ref{table:type_a_benchmark_V1}). 
\begin{table}[H]
\TABLE
{Computational results for the extended DARP benchmark instances
\label{table:type_a_benchmark_V1}}
{
\footnotesize
\begin{tabular}{c|rrrr|rrrr|rrrr}\toprule
    \multirow{2}{*}{Instance} & \multicolumn{2}{c}{DARP}& \multirow{2}{*}{$T_{MP}$} & \multirow{2}{*}{$T_{SP}$} & \multicolumn{2}{c}{RDARP($\epsilon^{\text{risk}}=30$)}& \multirow{2}{*}{$T_{MP}$} & \multirow{2}{*}{$T_{SP}$} & \multicolumn{2}{c}{RDARP($\epsilon^{\text{risk}}=15$)} & \multirow{2}{*}{$T_{MP}$} & \multirow{2}{*}{$T_{SP}$} \\\cmidrule(lr){2-3}\cmidrule(lr){6-7}\cmidrule(lr){10-11}
    & Cost & $\overline{H}$ &&& Cost & $\overline{H}$ & && Cost & $\overline{H}$ & & \\ \hline
a2-16 & 294.25  & 19.00   & 0.1     & 0.1   & 294.25  & 19.00 & 0.1    & 0.2    & 318.63  & 13.13  & 0.1   & 0.3    \\
a2-20 & 344.83  & 15.33   & 0.1   & 0.8   & 344.83  & 15.33 & 0.2    & 2.4    & 380.12  & 12.69  & 0.2   & 1.7    \\
a2-24 & 431.12  & 36.57   & 0.1   & 1.2   & 441.06  & 17.75 & 0.1    & 5.5    & 441.57  & 13.72  & 0.1   & 3.7    \\
a3-24 & 344.83  & 20.01   & 0.1   & 1.1   & 344.83  & 20.01 & 0.1    & 2.0    & 353.09  & 14.74  & 0.1   & 1.7    \\
a3-30 & 494.84  & 20.76   & 0.1   & 2.8   & 494.84  & 20.76 & 0.4    & 4.2     & $\emptyset$       & $\emptyset$     &  $\emptyset$   & $\emptyset$      \\
a3-36 & 583.18  & 47.08   & 1.4   & 14.6  & 595.09  & 21.40 & 2.2    & 18.1   & 636.70  & 13.89  & 0.7   & 16.9   \\
a4-32 & 485.49  & 30.01   & 0.1   & 1.2   & 487.43  & 24.90 & 0.2    & 2.2    & 515.04  & 12.96  & 0.9   & 9.1    \\
a4-40 & 557.68  & 22.23   & 0.5   & 8.2   & 557.68  & 22.23 & 0.9    & 10.2   & 644.19  & 14.16  & 0.5   & 8.2    \\
a4-48 & 668.81  & 45.11   & 3.4   & 40.4  & 678.58  & 19.01 & 38.2   & 641.3  & 686.74  & 13.32  & 5.8   & 159.9  \\
a5-40 & 498.40   & 30.48   & 0.5   & 4.7   & 499.16  & 26.34 & 0.3    & 9.4    & 529.46  & 14.45  & 1.2   & 26.9   \\
a5-50 & 686.62  & 29.66   & 7.9   & 86.2  & 686.62  & 29.66 & 1.5    & 55.4   & 707.98  & 14.44  & 1.8   & 58.5   \\
a5-60 & 808.42  & 44.90   & 6.7   & 77.5  & 817.44  & 25.97 & 79.5   & 1811.9 & 862.72  & 14.69  & 3.6   & 334.9  \\
a6-48 & 604.12  & 34.16   & 0.8   & 8.6   & 607.94  & 23.94 & 0.6    & 24.4   & 641.83  & 12.48  & 0.3   & 15.3   \\
a6-60 & 819.24  & 46.50   & 4.4   & 35.1  & 833.90  & 28.22 & 3.5    & 86.6   & 891.66  & 14.26  & 1.2   & 54.0   \\
a6-72 & 916.04  & 49.21   & 30.3  & 442.7 & 922.97  & 29.28 & 18.7   & 664.5  & 1019.97 & 14.83  & --  & -- \\
a7-56 & 724.03  & 37.40   & 2.5   & 31.1  & 725.23  & 28.96 & 0.6    & 36.6   & 770.26  & 15.00  & 17.8  & 254.7  \\
a7-70 & 889.11  & 46.25   & 9.0     & 70.1  & 921.80  & 26.11 & --     & -- & 967.76  & 14.90  & 41.2  & 2401.5 \\
a7-84 & 1033.36 & 37.11   & 94.4  & 710.9 & 1042.00 & 22.75 & --  & -- & 1100.99 & 14.68  & -- & -- \\
a8-64 & 747.46  & 35.13   & 3.3   & 26.7  & 751.95  & 28.25 & 1.7    & 79.7   & 796.89  & 15.00  & 3.0   & 89.3   \\
a8-80 & 945.72  & 35.99   & 27.6  & 330.6 & 993.04  & 28.29 & --   & -- & 978.23  & 14.99  & 2.6   & 198.7  \\
a8-96 & 1229.65 & 45.68   & 120.7 & 774.2 & 1250.05 & 27.96 & --  & -- & 1303.43 & 14.82  & --  & -- \\\hline
DARP(\%) & 0       & 0       & 0     & 0     & 1.08    & -44.40 & 163.11 & 560.60 & 6.86    & -55.81 & 78.95 & 385.02 \\\bottomrule
\end{tabular}}
{\vspace{-0.5cm}
\begin{itemize}
    \item [--] CPU time limit of one hour
    \item [$\emptyset$] No feasible solution at root node
\end{itemize}
}
\end{table}

\subsection{Experiments on Paratransit Trip Data}
In this section, we investigate the performance of our BCP algorithm on the instances extracted from the real-world dataset.{Specifically, we first compare the optimal routing under the DARP ($\epsilon^{\text{risk}}\rightarrow\infty$) and the RDARP ($\epsilon^{\text{risk}}=15$ and $30$) to investigate the trade-off between travel cost and exposed risk. The historical trip records are shown as a baseline to demonstrate the potential improvement in our RDARP cases.  Next, we present the Pareto fronts by solving four medium-size real-world instances for both AM and PM periods. Further, we conduct sensitivity analyses on the time window width ($b_i-a_i$). Finally, the EDARP results are presented under two different upper bounds $\epsilon^{\text{dt}}=2$ and $4$, which will be compared with the DARP results and the historical routes.}

\subsubsection{Comparison between Historical routes, DARP, and RDARP}
\begin{table}
\TABLE
{Computational results for the real-world instances
\label{table:type_a_para_V1}}
{\footnotesize
\begin{tabular}{r|rr|rrrr|rrrr|rrrr}\toprule
    \multirow{2}{*}{Instance} & \multicolumn{2}{c}{Historical}& \multicolumn{2}{c}{DARP}& \multirow{2}{*}{$T_{MP}$} & \multirow{2}{*}{$T_{SP}$} & \multicolumn{2}{c}{RDARP($\epsilon^{\text{risk}}=30$)}& \multirow{2}{*}{$T_{MP}$} & \multirow{2}{*}{$T_{SP}$} & \multicolumn{2}{c}{RDARP($\epsilon^{\text{risk}}=15$)} & \multirow{2}{*}{$T_{MP}$} & \multirow{2}{*}{$T_{SP}$} \\\cmidrule(lr){2-3}\cmidrule(lr){4-5}\cmidrule(lr){8-9}\cmidrule(lr){12-13}
    & Cost & $\overline{H}$ & Cost & $\overline{H}$&&& Cost & $\overline{H}$ & && Cost & $\overline{H}$ & & \\ \hline
AM-17-3  & 713.29  & 138.10 & 526.76  & 42.50 & 0.0 & 0.3   & 539.21  & 18.92 & 0.0 & 0.2   & 541.68  & 10.22 & 0.0 & 0.2    \\
AM-20-4  & 930.94  & 138.10 & 611.17  & 55.43 & 0.0 & 0.3   & 614.91  & 20.41 & 0.0 & 0.3   & 630.24  & 9.51  & 0.0 & 0.4    \\
AM-25-5  & 1093.53 & 138.10 & 747.60  & 42.50 & 0.1 & 1.1   & 788.75  & 20.67 & 0.4 & 2.0   & 808.84  & 12.18 & 0.6 & 3.0    \\
AM-25-6  & 1010.40 & 138.10 & 857.97  & 25.97 & 0.6 & 0.1   & 857.97  & 25.94 & 0.1 & 0.6   & 886.31  & 12.93 & 0.1 & 0.7    \\
AM-32-7  & 1314.34 & 138.10 & 885.23 & 60.46 & 0.2 & 17.2  & 930.62  & 29.54 & 0.3 & 9.7   & 966.44  & 14.27 & 0.9 & 27.4   \\
AM-39-8  & 1739.28 & 138.10 & 1187.73 & 65.69 & 4.1 & 48.4  & 1206.00 & 28.44 & 5.9 & 45.1  & 1240.03 & 13.94 & 0.8 & 17.7   \\
AM-43-10 & 2037.03 & 138.10 & 1267.27 & 22.25 & 0.4 & 16.0  & 1264.56 & 27.99 & 5.5 & 78.0  & 1352.26 & 14.19 & 4.3 & 61.7   \\
AM-46-11 & 2145.44 & 138.10 & 1317.09 & 22.25 & 0.7 & 49.5  & 1317.09 & 22.25 & 5.0 & 93.5  & 1402.92 & 13.71 & 0.5 & 10.3   \\
AM-50-12 & 2253.21 & 138.10 & 1367.72 & 22.25 & 0.6 & 111.2 & 1367.72 & 22.25 & 8.0 & 209.9 & 1453.91 & 14.19 & 0.8 & 38.5   \\
AM-55-13 & 2493.91 & 138.10 & 1494.76 & 37.85 & 1.5 & 300.6 & 1507.15 & 28.44 & 1.1 & 98.7  & 1569.30 & 13.94 & 1.3 & 53.2   \\
PM-17-3  & 728.16  & 26.18  & 610.85  & 11.06 & 0.0 & 0.2   & 644.59  & 22.36 & 0.0 & 0.3   & 646.18  & 12.32 & 0.0 & 0.4    \\
PM-20-4  & 928.20  & 26.18  & 661.51  & 30.17 & 0.0 & 1.5   & 664.34  & 29.90 & 0.0 & 1.7   & 686.71  & 14.74 & 0.1 & 1.9    \\
PM-25-5  & 1151.16 & 26.18  & 819.88  & 36.64 & 0.6 & 11.7  & 819.88  & 28.46 & 0.1 & 1.6   & 885.36  & 14.75 & 0.1 & 1.8    \\
PM-25-6  & 1364.82 & 119.83 & 817.77  & 39.75 & 1.7 & 0.1   & 825.44  & 26.39 & 0.0 & 2.1   & 843.76  & 14.36 & 0.0 & 1.4    \\
PM-32-7  & 1471.93 & 49.23  &985.84 & 38.18 & 0.1 & 117.3   & 1130.74 & 23.26 & 0.1 & 1.6   & 1145.22 & 9.72  & 0.0 & 0.9    \\
PM-39-8  & 1797.39 & 49.23  & 1192.92 & 95.75 & 0.2 & 69.4  & 1219.34 & 29.90 & 0.2 & 20.9  & 1281.63 & 12.81 & 0.3 & 22.9   \\
PM-43-10 & 2261.12 & 119.83 & 1327.88 & 65.24 & 0.2 & 534.9 & 1341.68 & 29.90 & 0.2 & 64.5  & 1441.66 & 14.74 & 0.7 & 69.7   \\
PM-46-11 & 2448.27 & 119.83 & 1390.73 & 30.13 & 1.4 & 952.6 & 1393.56 & 29.90 & 1.3 & 525.8 & 1442.72 & 12.81 & 3.5 & 418.6  \\
PM-50-12 & 2069.71 & 49.23  & -       & -     & -   & -     & -       & -     & -   & -     & 1530.21 & 13.50 & 8.7 & 1504.8 \\
PM-55-13 & 2668.55 & 119.83 & -       & -     & -   & -     & -       & -     & -   & -     & 1663.26 & 14.74 & -   & -     \\\hline
Hist(\%) & 0 & 0 & -32.37 &-48.54 & \textbackslash & \textbackslash & -31.32 & 	-63.03 & \textbackslash & \textbackslash& 	-29.06 & -81.32 & \textbackslash & \textbackslash \\
DARP(\%) & \textbackslash & \textbackslash& 0 & 0 & 0 & 0 & 1.47 & 	-18.23& -28.16 &	-32.41& 	5.43&-57.73& 	21.65&	-34.78\\
\bottomrule
\end{tabular}}
{--:  CPU time limit of one hour;\textbackslash:  Not applicable}
\end{table}

Table~\ref{table:type_a_para_V1} shows the computational results of DARP and RDARP under $\epsilon^{\text{risk}}=15$ and $30$. We also present the ground-truth information (total travel cost and exposed risk) based on the historical routes, which helps to understand the improvement on our risk-aware routing strategy. Considering the historical routes or DARP solution as the baseline, we summarize the percentage changes of the RDARP solutions in Hist(\%) and DARP (\%), respectively.

In general, 19 of 20 instances can be solved to the optimal within the time limit of one hour with the most restrictive upper bound $\epsilon^{\text{risk}}=15$, and meaningful IP solutions can be obtained in all 20 instances under $\epsilon^{\text{risk}}=15$. In both DARP and RDARP ($\epsilon^{\text{risk}}=30$), 18 of 20 instances are solved to the optimality within the CPU time limit. 

Compared with the historical routes, we report that one can achieve both improvements in both total travel cost and exposed risk. On average, the DARP solutions suggest a travel cost savings of 32.37\% and maximum risk reduction of 48.54\%. The reduction on the maximum exposed risk can be further enhanced to 81.32\% (from over \SI{100}{min} to less than \SI{15}{min}) in the case of $\epsilon^{\text{risk}}=15$, where a 29.06\% savings on travel cost is still obtained compared with the historical routes. The significant reduction in exposed risk roots in the fact that the riskiest passenger (denoting $\overline{H}$=\SI{138.10}{min}) is alternatively assigned to other routes with few co-riders or served by a dedicated door-to-door paratransit service in our DARP or RDARP solution. 
Furthermore, compared with the DARP, the average performance of RDARP under $\epsilon^{\text{risk}}=15$ can be translated to 5.43\% more travel cost, but yielding 57.73\% fewer maximum exposed risk (from more than \SI{40}{min} to about \SI{13}{min} on average). The average improvements are found to be more marginal in the case of $\epsilon^{\text{risk}}=30$, where 1.47\% more travel cost can lead to an 18.23\% reduction on the maximum exposed risk.

By comparing the CPU time between the AM and PM instances, we report more computational efforts of the PM instances than the AM instances. One possible reason is the different distribution of requests and the property of trips (e.g., inbound for PM vs. outbound for AM), which makes it harder to identify the feasible trips with a low-risk level in the PM period.
For instance, the trips during the PM period serve as the trips back home, which involve a strict time window at their origins and a loose time window at the destinations (e.g., home locations). It is likely that a group of passengers reserved a similar pick-up time window to align with their activity schedule (e.g., the end time for social activity and routine appointments at the clinics). In this case, the paratransit shuttle must collect the passengers in one shot to meet the strict time window and sequentially drop them off, leading to a heavily overlapped onboard time and exposed risk.
\subsubsection{Exact Pareto Fronts and Sensitivity Analyses}

\begin{figure}[H]
    \centering
    \caption{Exact Pareto fronts for total travel cost and maximum risk}\label{fig:pf_risk}
    \subfloat[Pareto front and historical record (AM)]{\includegraphics[width=0.5\linewidth]{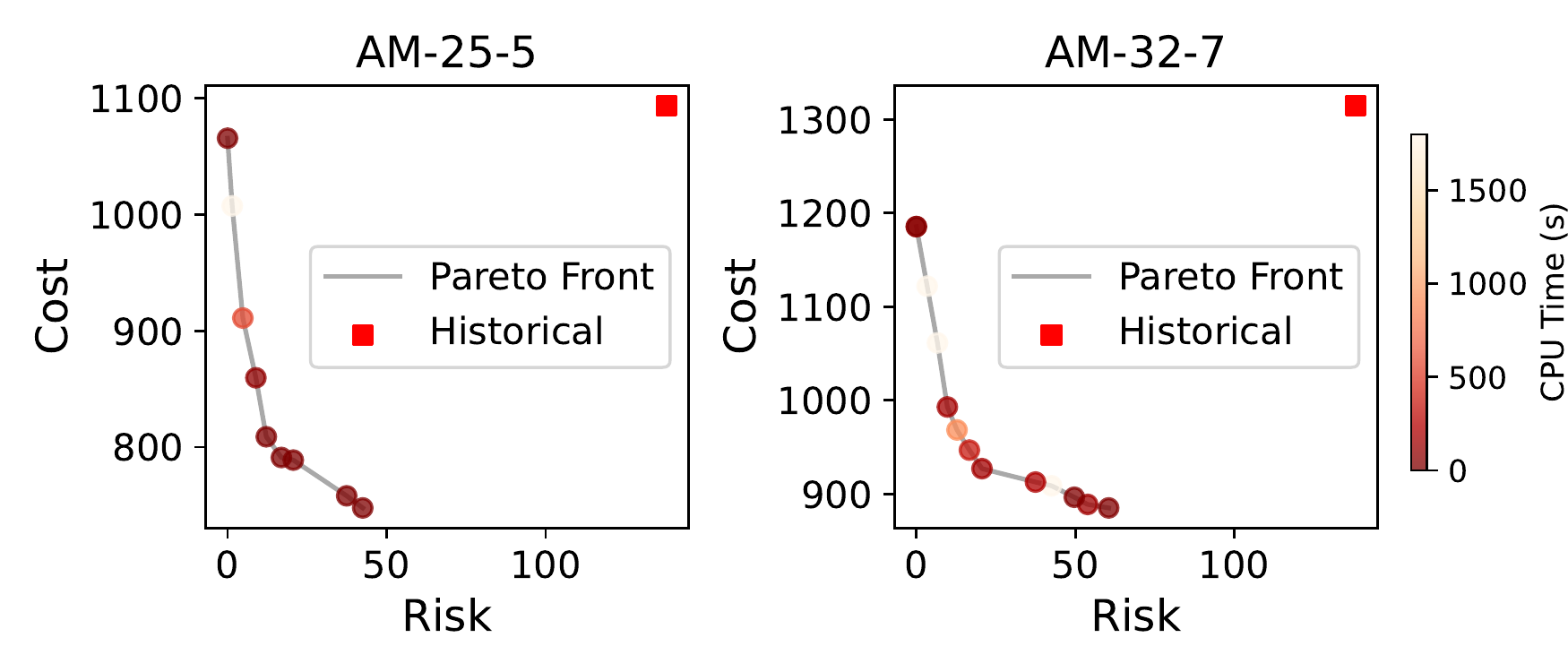}\label{fig:pf_am}}
    \subfloat[Pareto front and historical record (PM)]{\includegraphics[width=0.5\linewidth]{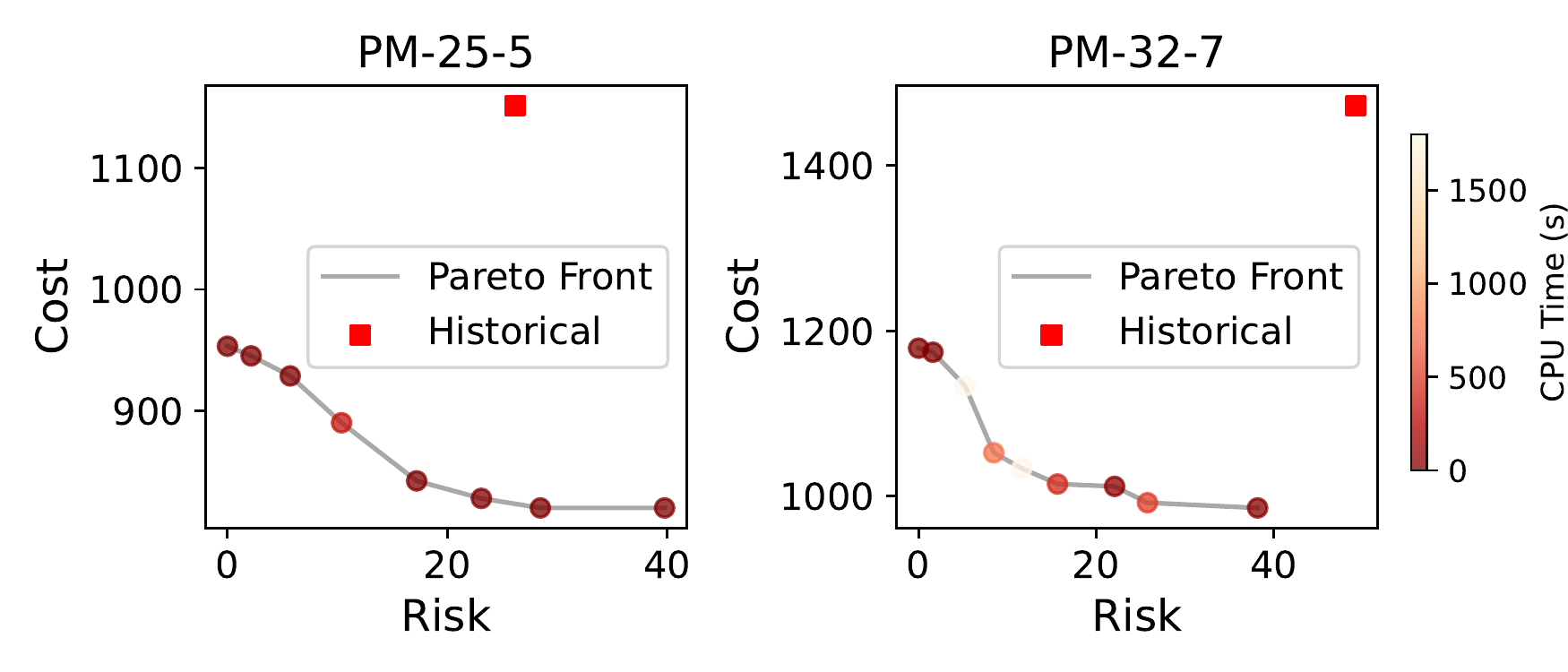}\label{fig:pf_pm}}\vfill
    \subfloat[Sensitivity analyses on time window (AM)]{\includegraphics[width=0.45\linewidth]{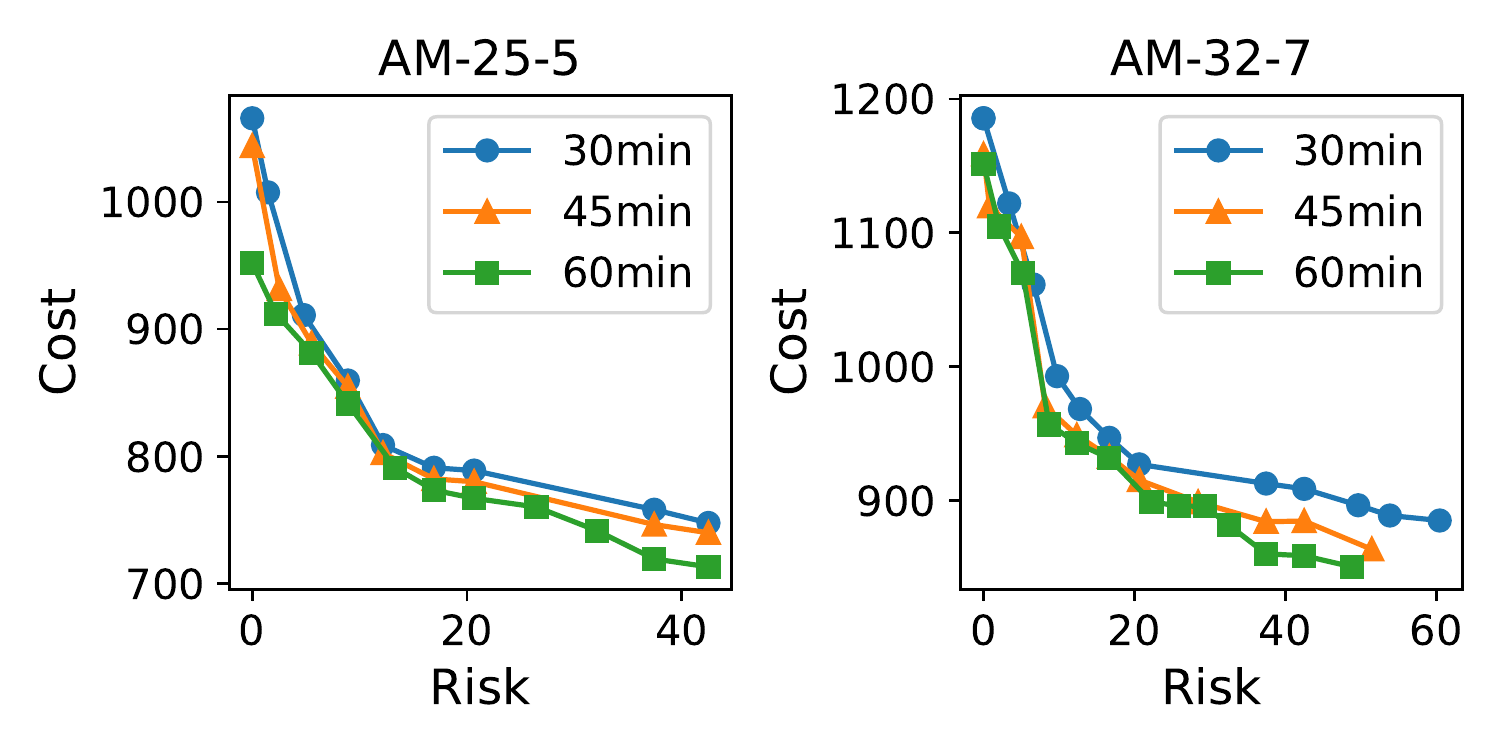}\label{fig:pf_am_tw}}
    \subfloat[Sensitivity analyses on time window (PM)]{\includegraphics[width=0.45\linewidth]{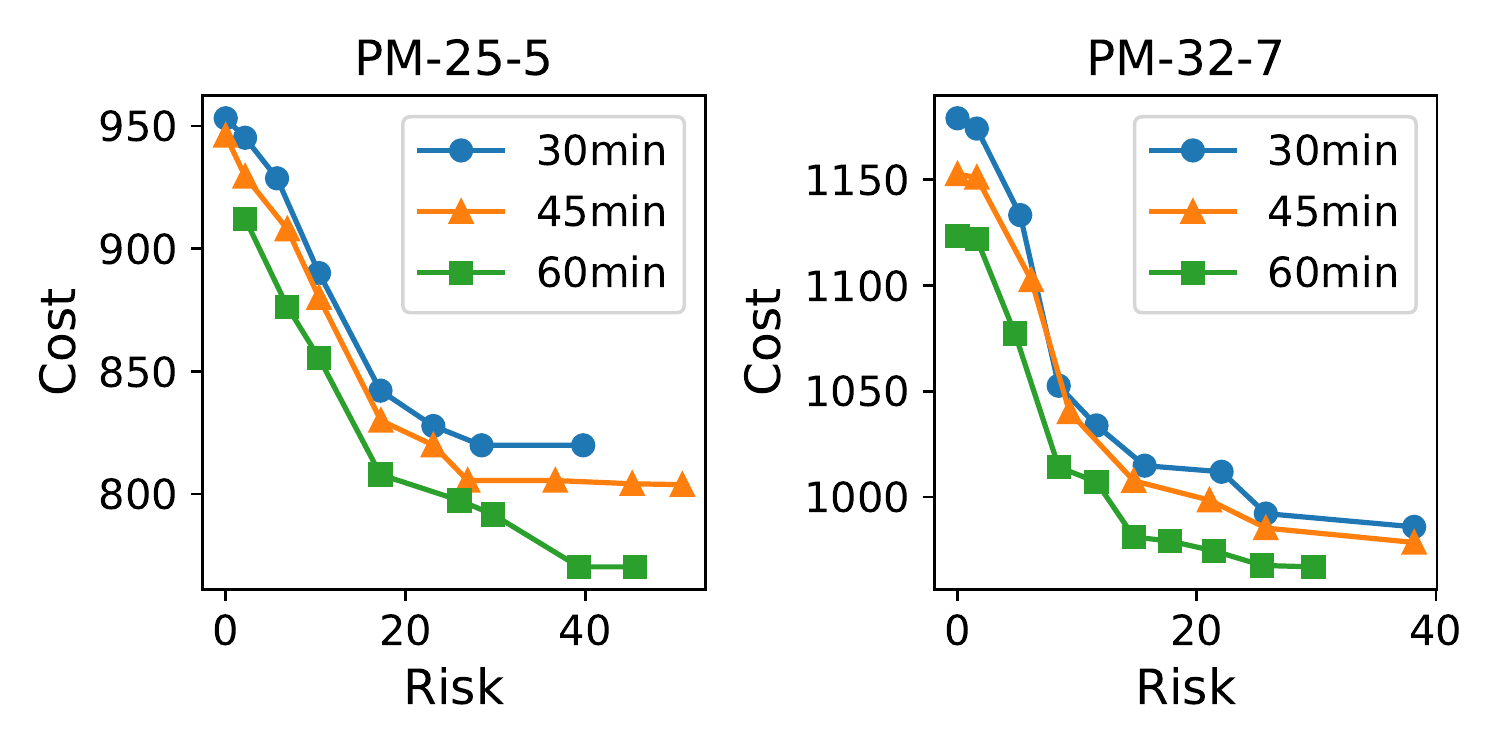}\label{fig:pf_pm_tw}}
\end{figure}

Figure~\ref{fig:pf_risk} displays the exact Pareto frontiers of four instances (AM-25-5, AM-32-7, PM-25-5, and PM-32-7). Specifically, Figures~\ref{fig:pf_am} and \ref{fig:pf_pm} compare the Pareto fronts with the historical records, where detailed information (cost- and risk-minimizing solutions and total CPU time to obtain all solutions on the Pareto fronts) are shown in Table~\ref{table:pareto_front_real-world_V1}. In Figures~\ref{fig:pf_am} and \ref{fig:pf_pm}, the historical records locate at the upper-right corner in all instances except PM-25-5, which indicates at least 31\% and 22\% improvements on the total travel cost and maximum risk in our risk-aware routing strategies.
On the other hand, the Pareto fronts illustrate the major trade-off between total travel cost and maximum risk. Taking AM-32-7 as an example, the total travel cost ranges between 885.23 and 1185.55 (33.93\% higher), whereas the associated maximum risk varies between \SI{60.46}{min} and \SI{0}{min} (a risk-free routing). Similar trends are observed in the other three instances. In particular, for the case PM-32-7, the risk-free routing strategy (a maximum risk reduction of \SI{38.18}{min}) can be obtained with a slight increase of travel cost of \SI{193.17}{min} (19.59\%).
In addition to the two extreme points, the Pareto fronts in all instances exhibit a strong non-linear relationship, which can be adopted to gain insights on the operational impacts of implementing a certain level of service, e.g., maximum exposed risk in our case. For instance, in AM-32-7, tangible improvements in risk mitigation from \SI{60.46}{min} to \SI{20.63}{min} (-65.88\%) are observed by increasing the total travel cost from 885.27 to 927.15 (+4.7\%). In this case, the mid-point ($20.63,927.15$) suggests a cost-effective routing strategy, and the risk-minimizing point ($0,1185.55$) may be not favorable by the transit service operator due to the high budget for travel cost ($33.93$\% higher cost). Hence, with granular information of the operational impacts, the Pareto front is deemed to be useful for the transit service operator's decision-making (e.g., DRT agencies, paratransit operators, ride-sharing companies). 

In Figures~\ref{fig:pf_am_tw} and~\ref{fig:pf_pm_tw}, we conduct the sensitivity analyses on the relaxed time window width, and the experiments are performed on the same instances in Figures~\ref{fig:pf_am} and \ref{fig:pf_pm} (for detailed computational results, see Table~\ref{table:pareto_front_real-world_TW_V1}). The realistic consideration is that passengers may choose to afford more flexibility for potentially lower exposed risk in the pandemic world.  
In this regard, we consider enlarging the original time window width of \SI{30}{min} to be \SI{45}{min} and \SI{60}{min}.

From Figures~\ref{fig:pf_am_tw} and~\ref{fig:pf_pm_tw}, we observe that the larger time window size contributes to the savings of total travel cost and reduction of the maximum risk in both AM and PM instances. In particular, the performance of PM instances is observed to be more sensitive to the time window width. Regarding the risk-free routing solution ($\overline{H}=0$) in instance PM-32-7, the savings on travel cost can achieve up to 2.2\% and 4.7\% by relaxing the time window to \SI{45}{min} and \SI{60}{min}, respectively. 
However, in AM-32-7, the savings of travel cost is only up to 2.0\% and 2.2\% under the time windows of \SI{45}{min} and \SI{60}{min} in the risk-minimizing routing strategy. We note that a more significant improvement in the PM instances is due to the property of inbound trips (from target locations back to home locations), where a relaxed time window at target points makes it possible for the vehicle to serve the passengers in sequence. 

\vspace{-1cm}

\begin{table}
\TABLE
{Computational results of Pareto fronts for the real-world instances
\label{table:pareto_front_real-world_V1}}
{\footnotesize
\begin{tabular}{rr|rr|rr|rr|rr|rr}\toprule
    \multirow{2}{*}{Instance} & \multirow{2}{*}{\# Sol} & \multicolumn{2}{c}{Historical} &\multicolumn{2}{c}{DARP sol} & \multicolumn{2}{c}{Min-Risk RDARP sol} & \multicolumn{2}{c}{DARP vs. RDARP} & \multirow{2}{*}{$T_{MP}$} & \multirow{2}{*}{$T_{SP}$} \\\cmidrule(lr){3-4}\cmidrule(lr){5-6}\cmidrule(lr){7-8}\cmidrule(lr){9-10}
    && Cost & $\overline{H}$ & Cost & $\overline{H}$ & Cost & $\overline{H}$ & Cost (\%) & $\overline{H}$(\%) & & \\ \hline
  AM-25-5 &    9 &            1093.53 &        138.10 & 747.60 & 42.50 &   1065.78 &       0.0 &     42.56 &     100.0 & 1256.9 &  1930.1 \\
  AM-32-7 &      13 &            1314.34 &        138.10 & 885.23 & 60.46 &   1185.55 &       0.0 &     33.93 &     100.0 & 7415.8 & 11348.3 \\
  PM-25-5 &     8 &            1151.16 &         26.18 & 819.88 & 39.75 &    953.15 &       0.0 &     16.25 &     100.0 &   43.1 &   350.1 \\
  PM-32-7 &      9 &            1471.93 &         49.23 & 985.84 & 38.18 &   1179.01 &       0.0 &     19.59 &     100.0 & 6452.2 & 30330.3 \\\bottomrule
\end{tabular}}{}
\end{table}

\subsubsection{EDARP Results of Real-world Instances}
We finally present the EDARP results, which serve as one of the general DARPs with the min-max objectives. The EDARP aims to minimize the travel cost while imposing the upper bounds on the maximum detour rate ($\epsilon^{\text{dt}}$). In EDARP, the upper bounds $\epsilon^{\text{dt}}$ are set to $4$ and $2$, based on our observation that $\epsilon^{\text{dt}}$ varies between 2.80 and 6.01 in the DARP solutions (see Columns 2-3 in Table~\ref{table:edarp_realworld_V1}). Also, for the weights $t_{i,n+i}$ in the detour rate term $\frac{H_i}{t_{i,n+i}}$, we round up the weights $t_{i,n+i}<15$ to \SI{15}{min} such that extreme high detour rate due to short trips can be avoided. In this case, the detour rate is $\frac{H_i}{\max\{15,t_{i,n+i}\}}$. 
For instance, a ride time of \SI{15}{min} and a direct trip time of \SI{3}{min}, result in a detour rate of $\frac{15}{\max(15,3)}=1$ instead of $\frac{15}{3}=5$. The preprocessing step may produce detour rates below 1, such as a short trip of \SI{10}{min} yielding a rate of $\frac{2}{3}<1$. While these rates are not directly relevant to our analysis, our aim is to develop a metric that quantifies detours and focuses on relatively long trips. Additionally, based on historical routes, detour rates greater than 2 and 4 are always present. Hence, our objective is still to obtain a min-max solution with a maximum detour rate of no more than $\epsilon^{\text{dt}}=2$ or $4$.

\begin{table}
\TABLE
{Computational results of the EDARP for the real-world instances
\label{table:edarp_realworld_V1}}
{\footnotesize
\begin{tabular}{r|rrrr|rrrr|rrrr}\toprule
    \multirow{2}{*}{Instance} & \multicolumn{2}{c}{DARP sol}& \multirow{2}{*}{$T_{MP}$} & \multirow{2}{*}{$T_{SP}$} & \multicolumn{2}{c}{EDARP ($\epsilon^{\text{dt}}=4$)}& \multirow{2}{*}{$T_{MP}$} & \multirow{2}{*}{$T_{SP}$} & \multicolumn{2}{c}{EDARP($\epsilon^{\text{dt}}=2$)} & \multirow{2}{*}{$T_{MP}$} & \multirow{2}{*}{$T_{SP}$} \\\cmidrule(lr){2-3}\cmidrule(lr){6-7}\cmidrule(lr){10-11}
    & Cost & $\overline{D}$ &&& Cost & $\overline{D}$ & && Cost & $\overline{D}$ & & \\ \hline
AM-17-3  & 526.76  & 2.95 & 0.0 & 0.2   & 526.76  & 2.95 & 0.0 & 0.2   & 585.06  & 1.82 & 0.0 & 0.1   \\
AM-20-4  & 611.17  & 3.14 & 0.0 & 0.3   & 611.17  & 3.14 & 0.0 & 0.3   & 647.39  & 1.82 & 0.0 & 0.2   \\
AM-25-5  & 747.60  & 4.48 & 0.0 & 2.4   & 750.83  & 2.95 & 0.0 & 2.0   & 805.58  & 1.99 & 0.0 & 1.2   \\
AM-25-6  & 857.97  & 5.40 & 0.0 & 0.8   & 858.82  & 3.16 & 0.0 & 0.8   & 876.91  & 1.99 & 0.0 & 0.5   \\
AM-32-7  & 880.40  & 3.14 & 0.0 & 12.7  & 880.40  & 3.14 & 0.0 & 9.9   & 972.40  & 1.99 & 0.0 & 4.2   \\
AM-39-8  & 1187.73 & 3.63 & 4.1 & 48.4  & 1187.73 & 3.63 & 0.1 & 13.6  & 1271.71 & 1.99 & 0.0 & 6.6   \\
AM-43-10 & 1267.27 & 2.80 & 0.8 & 39.6  & 1267.27 & 2.80 & 0.1 & 14.4  & 1355.24 & 1.99 & 0.1 & 9.3   \\
AM-46-11 & 1317.09 & 2.80 & 0.7 & 49.5  & 1317.09 & 2.80 & 0.1 & 23.6  & 1408.36 & 1.99 & 0.1 & 13.3  \\
AM-50-12 & 1367.72 & 2.80 & 0.6 & 86.0  & 1367.72 & 2.80 & 0.0 & 52.2  & 1464.17 & 1.99 & 0.2 & 33.9  \\
AM-55-13 & 1494.76 & 3.48 & 1.5 & 300.6 & 1494.76 & 3.48 & 0.1 & 101.4 & 1607.57 & 1.99 & 5.0 & 866.6 \\
PM-17-3  & 644.59  & 3.06 & 0.0 & 0.5   & 644.59  & 3.06 & 0.0 & 0.6   & 646.18  & 1.43 & 0.0 & 0.1   \\
PM-20-4  & 661.51  & 5.30 & 0.0 & 2.5   & 687.93  & 3.62 & 0.0 & 1.7   & 705.79  & 1.84 & 0.0 & 1.5   \\
PM-25-5  & 819.87  & 6.01 & 0.0 & 4.5   & 822.52  & 3.61 & 0.0 & 6.1   & 841.94  & 1.86 & 0.0 & 0.8   \\
PM-25-6  & 817.77  & 3.59 & 0.0 & 2.6   & 817.77  & 3.59 & 0.0 & 3.3   & 843.76  & 1.95 & 0.0 & 1.0   \\
PM-32-7  & 1130.96 & 4.11 & 0.1 & 1.5   & 1131.48 & 3.40 & 0.0 & 1.3   & 1173.97 & 1.87 & 0.0 & 1.4   \\
PM-39-8  & 1192.92 & 4.90 & 0.2 & 69.4  & 1207.57 & 3.62 & 0.1 & 126.0 & 1246.90 & 1.84 & 0.1 & 29.7  \\
PM-43-10 & 1327.88 & 4.11 & 0.2 & 534.9 & 1341.40 & 3.89 & 0.0 & 958.5 & 1388.48 & 1.94 & 0.1 & 80.9  \\
PM-46-11 & 1390.73 & 4.90 & 1.4 & 952.6 & 1404.71 & 3.72 & 0.1 & 630.8 & 1452.12 & 1.94 & 0.2 & 252.9 \\
PM-50-12 & -       & -    & -   & -     & -       & -    & -   & -     & 1493.38 & 1.95 & 0.1 & 973.1 \\
PM-55-13 & -       & -    & -   & -     & -       & -    & -   & -     & -       & -    & -   & -    \\\hline
DARP(\%) & 0 & 0 & 0 & 0 & 0.45 & -12.22 & -90.89& -9.62 & 5.73 & -48.49 & -36.31 & -44.88\\
\bottomrule
\end{tabular}}
{--: CPU time limit of one hour\\
CPU times less than \SI{0.05}{sec} for $T_{MP}$ or $T_{SP}$ are rounded up to 0.1. }
\end{table}

Table~\ref{table:edarp_realworld_V1} presents three sets of results: Columns 2-5 represents the DARP solution as a baseline, where the maximum detour rate $\overline{D}$ serves as a by-product while solving the DARP. In addition, we present the EDARP results under $\epsilon^{\text{dt}}=2$ (Columns 6-9) and $\epsilon^{\text{dt}}=4$ (Columns 10-13), including the total travel cost, maximum detour rate, and CPU times for MP and SP.
Compared with the DARP results, the EDARP solutions under $\epsilon^{\text{dt}}=4$ are obtained with only a 0.45\% increase of travel cost on average. With a more restrictive requirement that the ride time cannot exceed twice of the trip time ($\epsilon^{\text{dt}}=2$), we report an average of 5.73\% higher travel cost, resulting in an improvement of maximum ride time by 48.49\%. 
Regarding the CPU times, less computational efforts are reported in all AM instances for the cases of $\epsilon^{\text{dt}}=4$, and the CPU times are further reduced under the case of $\epsilon^{\text{dt}}=2$ except the largest instance ``AM-55-13". This indicates that more restrictive constraints will help to truncate infeasible label extensions in advance, resulting in reduced searching space in the CG procedure and less computational effort. As for the PM instances, all solvable instances, except ``PM-17-3" and ``PM-25-6", are associated with a maximum detour rate of $\overline{D}>4$. In this case, more computational efforts are spent to identify the EDARP solutions that are aligned with $\overline{D}\le \epsilon^{\text{dt}}=4$. We also report that a more restrictive cap $\epsilon^{\text{dt}}=2$ will significantly reduce the CPU times, especially for larger PM instances (e.g., PM-39-8, PM-43-10, and PM-46-11).
\section{Conclusion}~\label{sec:conclusion}
In this study, we proposed the RDARP, a variant of the DARP that minimizes total travel cost and maximum exposed risk. 
To solve RDARP efficiently, we developed a BCP algorithm and designed a CG procedure. The SP was formulated as an ESPPRCMMR, to adapt to the MMR condition in our RDARP, and was solved using a labeling algorithm with a tailored design for exposed risk minimization. We introduced additional resources to capture the individual exposed risk and proposed dominance rules, REFs, and a risk calibration algorithm to effectively obtain a feasible time schedule that satisfies the MMR condition given a partial trip. We also extended the model to EDARP, a general DARP with min-max objectives regarding the detour rate.

Computational results demonstrated the efficiency of our BCP algorithm, which optimally solved large-scale benchmark instances (``a8-96") and real-world instances (up to 55 passengers) within the time limit of one hour.
Further, we compared our RDARP results with the historical ground-truth paratransit trip data, where our RDARP suggests at least 29.06\% savings of travel cost and an 81.32\% reduction of maximum exposed risk.
Compared with the DARP solutions, our risk-aware routing strategies can significantly reduce the maximum exposed risk without heavily sacrificing the travel cost. On average, a 5.43\% higher travel cost can lead to a 57.73\% reduction of maximum exposed risk. 
Using the $\epsilon$-constraint method, sets of Pareto frontiers are presented to shed light on the operational impacts of implementing a certain level of maximum risk cap. Based on that, the transit service operators can make a better decision on the most effective routing strategy that aligns with the requirement of risk mitigation and the budget limit.

Future research will be devoted to heuristic algorithms that can further help to solve the instances with more passengers. Meanwhile, the DRT operation can be considered from a multi-objective optimization perspective, e.g., limited co-ride duration~\citep{cdc_coride_time2021}, reduced capacity due to social distancing policy~\citep{cdc_coride_time2021}, and multiple depots~\citep{braekers2014exact} that are located near clinics or hospitals. All those factors will be particularly practical for the DRT services operating in the urban and rural areas during the demand rebound in the post-pandemic period.

\ACKNOWLEDGMENT{%
We thank the Federal Transit Agency for funding this study through the Public Transportation COVID-19 Research Demonstration Grant Program, ClasTran for providing the necessary trip log data, and the Alabama Transportation Institute for their support.
}

\begingroup\footnotesize
\let\section\subsubsection

\bibliographystyle{informs2014trsc} 
\bibliography{ref.bib}

\begin{thebibliography}{49}
\providecommand{\natexlab}[1]{#1}
\providecommand{\url}[1]{\texttt{#1}}
\providecommand{\urlprefix}{URL }

\bibitem[{Allahyari, Yaghoubi, \protect\BIBand{}
  Van~Woensel(2021)}]{allahyari2021novel}
Allahyari S, Yaghoubi S, Van~Woensel T, 2021 \emph{A novel risk perspective on
  location-routing planning: An application in cash transportation}.
  \emph{Transportation Research Part E: Logistics and Transportation Review}
  150:102356.

\bibitem[{Barnhart et~al.(1998)Barnhart, Johnson, Nemhauser, Savelsbergh,
  \protect\BIBand{} Vance}]{barnhart1998branch}
Barnhart C, Johnson EL, Nemhauser GL, Savelsbergh MW, Vance PH, 1998
  \emph{Branch-and-price: Column generation for solving huge integer programs}.
  \emph{Operations research} 46(3):316--329.

\bibitem[{B{\'e}rub{\'e}, Gendreau, \protect\BIBand{}
  Potvin(2009)}]{berube2009exact}
B{\'e}rub{\'e} JF, Gendreau M, Potvin JY, 2009 \emph{An exact
  $\epsilon$-constraint method for bi-objective combinatorial optimization
  problems: Application to the traveling salesman problem with profits}.
  \emph{European journal of operational research} 194(1):39--50.

\bibitem[{Braekers, Caris, \protect\BIBand{}
  Janssens(2014)}]{braekers2014exact}
Braekers K, Caris A, Janssens GK, 2014 \emph{Exact and meta-heuristic approach
  for a general heterogeneous dial-a-ride problem with multiple depots}.
  \emph{Transportation Research Part B: Methodological} 67:166--186.

\bibitem[{Carlsson et~al.(2009)Carlsson, Ge, Subramaniam, Wu, \protect\BIBand{}
  Ye}]{carlsson2009solving}
Carlsson J, Ge D, Subramaniam A, Wu A, Ye Y, 2009 \emph{Solving min-max
  multi-depot vehicle routing problem}. \emph{Lectures on global optimization}
  55:31--46.

\bibitem[{CDC(2020)}]{cdc_coride_time2021}
CDC, 2020 \emph{Public health guidance for community-related exposure}.
  \url{https://www.cdc.gov/coronavirus/2019-ncov/php/public-health-recommendations.html},
  accessed 31 July 2021.

\bibitem[{{CDC}(2022)}]{cdc2022appendices}
{CDC}, 2022 \emph{Appendices}.
  \url{https://www.cdc.gov/coronavirus/2019-ncov/php/contact-tracing/contact-tracing-plan/appendix.html},
  accessed on October 2022.

\bibitem[{Chen et~al.(2020)Chen, Brozen, Rollman, Ward, Norris, Gregory,
  \protect\BIBand{} Zimmerman}]{chen2020transportation}
Chen KL, Brozen M, Rollman JE, Ward T, Norris K, Gregory KD, Zimmerman FJ, 2020
  \emph{Transportation access to health care during the covid-19 pandemic:
  Trends and implications for significant patient populations and health care
  needs}. Technical report, Institute of Transportation Studies, UCLA.

\bibitem[{Clastran(2021)}]{clastran2021}
Clastran, 2021 \emph{Clastran: Central alabama’s specialized transit}.
  \url{https://clastran.com/}, accessed on July 2021.

\bibitem[{Cordeau(2006)}]{cordeau2006branch}
Cordeau JF, 2006 \emph{A branch-and-cut algorithm for the dial-a-ride problem}.
  \emph{Operations Research} 54(3):573--586.

\bibitem[{Cordeau \protect\BIBand{} Laporte(2007)}]{cordeau2007dial}
Cordeau JF, Laporte G, 2007 \emph{The dial-a-ride problem: models and
  algorithms}. \emph{Annals of operations research} 153(1):29--46.

\bibitem[{Dantzig \protect\BIBand{} Wolfe(1960)}]{dantzig1960decomposition}
Dantzig GB, Wolfe P, 1960 \emph{Decomposition principle for linear programs}.
  \emph{Operations research} 8(1):101--111.

\bibitem[{Dayarian et~al.(2015{\natexlab{a}})Dayarian, Crainic, Gendreau,
  \protect\BIBand{} Rei}]{dayarian2015branch}
Dayarian I, Crainic TG, Gendreau M, Rei W, 2015{\natexlab{a}} \emph{A
  branch-and-price approach for a multi-period vehicle routing problem}.
  \emph{Computers \& Operations Research} 55:167--184.

\bibitem[{Dayarian et~al.(2015{\natexlab{b}})Dayarian, Crainic, Gendreau,
  \protect\BIBand{} Rei}]{dayarian2015column}
Dayarian I, Crainic TG, Gendreau M, Rei W, 2015{\natexlab{b}} \emph{A column
  generation approach for a multi-attribute vehicle routing problem}.
  \emph{European Journal of Operational Research} 241(3):888--906.

\bibitem[{Demir, Bekta{\c{s}}, \protect\BIBand{} Laporte(2014)}]{demir2014bi}
Demir E, Bekta{\c{s}} T, Laporte G, 2014 \emph{The bi-objective
  pollution-routing problem}. \emph{European Journal of Operational Research}
  232(3):464--478.

\bibitem[{Desaulniers et~al.(1998)Desaulniers, Desrosiers, Solomon, Soumis,
  Villeneuve et~al.}]{desaulniers1998unified}
Desaulniers G, Desrosiers J, Solomon MM, Soumis F, Villeneuve D, et~al., 1998
  \emph{A unified framework for deterministic time constrained vehicle routing
  and crew scheduling problems}. \emph{Fleet management and logistics}, 57--93
  (Springer).

\bibitem[{Dikas \protect\BIBand{} Minis(2018)}]{dikas2018scheduled}
Dikas G, Minis I, 2018 \emph{Scheduled paratransit transport enhanced by
  accessible taxis}. \emph{Transportation Science} 52(5):1122--1140.

\bibitem[{Dror(1994)}]{dror1994note}
Dror M, 1994 \emph{Note on the complexity of the shortest path models for
  column generation in vrptw}. \emph{Operations Research} 42(5):977--978.

\bibitem[{Dumas, Desrosiers, \protect\BIBand{} Soumis(1991)}]{dumas1991pickup}
Dumas Y, Desrosiers J, Soumis F, 1991 \emph{The pickup and delivery problem
  with time windows}. \emph{European journal of operational research}
  54(1):7--22.

\bibitem[{Feillet et~al.(2004)Feillet, Dejax, Gendreau, \protect\BIBand{}
  Gueguen}]{feillet2004exact}
Feillet D, Dejax P, Gendreau M, Gueguen C, 2004 \emph{An exact algorithm for
  the elementary shortest path problem with resource constraints: Application
  to some vehicle routing problems}. \emph{Networks: An International Journal}
  44(3):216--229.

\bibitem[{Fu(1999)}]{fu1999improving}
Fu L, 1999 \emph{Improving paratransit scheduling by accounting for dynamic and
  stochastic variations in travel time}. \emph{Transportation Research Record}
  1666(1):74--81.

\bibitem[{Fukasawa, He, \protect\BIBand{} Song(2016)}]{fukasawa2016branch}
Fukasawa R, He Q, Song Y, 2016 \emph{A branch-cut-and-price algorithm for the
  energy minimization vehicle routing problem}. \emph{Transportation Science}
  50(1):23--34.

\bibitem[{Gschwind \protect\BIBand{} Irnich(2015)}]{gschwind2015effective}
Gschwind T, Irnich S, 2015 \emph{Effective handling of dynamic time windows and
  its application to solving the dial-a-ride problem}. \emph{Transportation
  Science} 49(2):335--354.

\bibitem[{Guerra et~al.(2017)Guerra, Bolotin, Lim, Heffernan, Deeks, Li,
  \protect\BIBand{} Crowcroft}]{guerra2017basic}
Guerra FM, Bolotin S, Lim G, Heffernan J, Deeks SL, Li Y, Crowcroft NS, 2017
  \emph{The basic reproduction number (r0) of measles: a systematic review}.
  \emph{The Lancet Infectious Diseases} 17(12):e420--e428.

\bibitem[{Gupta et~al.(2010)Gupta, Chen, Miller, \protect\BIBand{}
  Surya}]{gupta2010improving}
Gupta D, Chen HW, Miller LA, Surya F, 2010 \emph{Improving the efficiency of
  demand-responsive paratransit services}. \emph{Transportation research part
  A: policy and practice} 44(4):201--217.

\bibitem[{Ho et~al.(2018)Ho, Szeto, Kuo, Leung, Petering, \protect\BIBand{}
  Tou}]{ho2018survey}
Ho SC, Szeto WY, Kuo YH, Leung JM, Petering M, Tou TW, 2018 \emph{A survey of
  dial-a-ride problems: Literature review and recent developments}.
  \emph{Transportation Research Part B: Methodological} 111:395--421.

\bibitem[{Hosni, Naoum-Sawaya, \protect\BIBand{}
  Artail(2014)}]{hosni2014shared}
Hosni H, Naoum-Sawaya J, Artail H, 2014 \emph{The shared-taxi problem:
  Formulation and solution methods}. \emph{Transportation Research Part B:
  Methodological} 70:303--318.

\bibitem[{Kohl et~al.(1999)Kohl, Desrosiers, Madsen, Solomon, \protect\BIBand{}
  Soumis}]{kohl19992}
Kohl N, Desrosiers J, Madsen OB, Solomon MM, Soumis F, 1999 \emph{2-path cuts
  for the vehicle routing problem with time windows}. \emph{Transportation
  Science} 33(1):101--116.

\bibitem[{Lehu{\'e}d{\'e} et~al.(2014)Lehu{\'e}d{\'e}, Masson, Parragh,
  P{\'e}ton, \protect\BIBand{} Tricoire}]{lehuede2014multi}
Lehu{\'e}d{\'e} F, Masson R, Parragh SN, P{\'e}ton O, Tricoire F, 2014 \emph{A
  multi-criteria large neighbourhood search for the transportation of disabled
  people}. \emph{Journal of the Operational Research Society} 65:983--1000.

\bibitem[{Liu, Luo, \protect\BIBand{} Lim(2015)}]{liu2015branch}
Liu M, Luo Z, Lim A, 2015 \emph{A branch-and-cut algorithm for a realistic
  dial-a-ride problem}. \emph{Transportation Research Part B: Methodological}
  81:267--288.

\bibitem[{Luo, Liu, \protect\BIBand{} Lim(2019)}]{luo2019two}
Luo Z, Liu M, Lim A, 2019 \emph{A two-phase branch-and-price-and-cut for a
  dial-a-ride problem in patient transportation}. \emph{Transportation Science}
  53(1):113--130.

\bibitem[{Luo et~al.(2017)Luo, Qin, Zhu, \protect\BIBand{} Lim}]{luo2017branch}
Luo Z, Qin H, Zhu W, Lim A, 2017 \emph{Branch and price and cut for the
  split-delivery vehicle routing problem with time windows and linear
  weight-related cost}. \emph{Transportation Science} 51(2):668--687.

\bibitem[{Matl, Hartl, \protect\BIBand{} Vidal(2018)}]{matl2018workload}
Matl P, Hartl RF, Vidal T, 2018 \emph{Workload equity in vehicle routing
  problems: A survey and analysis}. \emph{Transportation Science}
  52(2):239--260.

\bibitem[{Meyers(2007)}]{meyers2007contact}
Meyers L, 2007 \emph{Contact network epidemiology: Bond percolation applied to
  infectious disease prediction and control}. \emph{Bulletin of the American
  Mathematical Society} 44(1):63--86.

\bibitem[{Molenbruch, Braekers, \protect\BIBand{}
  Caris(2017)}]{molenbruch2017typology}
Molenbruch Y, Braekers K, Caris A, 2017 \emph{Typology and literature review
  for dial-a-ride problems}. \emph{Annals of Operations Research} 259:295--325.

\bibitem[{{NADTC}(2021)}]{nadtc2021}
{NADTC}, 2021 \emph{Ada \& paratransit}.
  \url{https://www.nadtc.org/about/transportation-aging-disability/ada-and-paratransit/},
  accessed 29 July 2021.

\bibitem[{Nie et~al.(2022)Nie, Qian, Guo, Jones, Doustmohammadi,
  \protect\BIBand{} Anderson}]{nie2022impact}
Nie Q, Qian X, Guo S, Jones S, Doustmohammadi M, Anderson MD, 2022 \emph{Impact
  of covid-19 on paratransit operators and riders: A case study of central
  alabama}. \emph{Transportation research part A: policy and practice}
  161:48--67.

\bibitem[{Overton(1988)}]{overton1988minimizing}
Overton ML, 1988 \emph{On minimizing the maximum eigenvalue of a symmetric
  matrix}. \emph{SIAM Journal on Matrix Analysis and Applications}
  9(2):256--268.

\bibitem[{Paquette, Cordeau, \protect\BIBand{}
  Laporte(2009)}]{paquette2009quality}
Paquette J, Cordeau JF, Laporte G, 2009 \emph{Quality of service in dial-a-ride
  operations}. \emph{Computers \& Industrial Engineering} 56(4):1721--1734.

\bibitem[{Paquette et~al.(2013)Paquette, Cordeau, Laporte, \protect\BIBand{}
  Pascoal}]{paquette2013combining}
Paquette J, Cordeau JF, Laporte G, Pascoal MM, 2013 \emph{Combining
  multicriteria analysis and tabu search for dial-a-ride problems}.
  \emph{Transportation Research Part B: Methodological} 52:1--16.

\bibitem[{Parragh(2011)}]{parragh2011introducing}
Parragh SN, 2011 \emph{Introducing heterogeneous users and vehicles into models
  and algorithms for the dial-a-ride problem}. \emph{Transportation Research
  Part C: Emerging Technologies} 19(5):912--930.

\bibitem[{Perrier, Langevin, \protect\BIBand{}
  Amaya(2008)}]{perrier2008vehicle}
Perrier N, Langevin A, Amaya CA, 2008 \emph{Vehicle routing for urban snow
  plowing operations}. \emph{Transportation Science} 42(1):44--56.

\bibitem[{Qian, Sun, \protect\BIBand{} Ukkusuri(2021)}]{qian2021scaling}
Qian X, Sun L, Ukkusuri SV, 2021 \emph{Scaling of contact networks for epidemic
  spreading in urban transit systems}. \emph{Scientific reports} 11(1):1--12.

\bibitem[{Qian \protect\BIBand{} Ukkusuri(2021)}]{qian2021connecting}
Qian X, Ukkusuri SV, 2021 \emph{Connecting urban transportation systems with
  the spread of infectious diseases: A trans-seir modeling approach}.
  \emph{Transportation Research Part B: Methodological} 145:185--211.

\bibitem[{Qu \protect\BIBand{} Bard(2015)}]{qu2015branch}
Qu Y, Bard JF, 2015 \emph{A branch-and-price-and-cut algorithm for
  heterogeneous pickup and delivery problems with configurable vehicle
  capacity}. \emph{Transportation Science} 49(2):254--270.

\bibitem[{Ropke \protect\BIBand{} Cordeau(2009)}]{ropke2009branch}
Ropke S, Cordeau JF, 2009 \emph{Branch and cut and price for the pickup and
  delivery problem with time windows}. \emph{Transportation Science}
  43(3):267--286.

\bibitem[{Ropke, Cordeau, \protect\BIBand{} Laporte(2007)}]{ropke2007models}
Ropke S, Cordeau JF, Laporte G, 2007 \emph{Models and branch-and-cut algorithms
  for pickup and delivery problems with time windows}. \emph{Networks: An
  International Journal} 49(4):258--272.

\bibitem[{Savelsbergh \protect\BIBand{} Sol(1998)}]{savelsbergh1998drive}
Savelsbergh M, Sol M, 1998 \emph{Drive: Dynamic routing of independent
  vehicles}. \emph{Operations Research} 46(4):474--490.

\bibitem[{Sayarshad \protect\BIBand{} Gao(2018)}]{sayarshad2018scalable}
Sayarshad HR, Gao HO, 2018 \emph{A scalable non-myopic dynamic dial-a-ride and
  pricing problem for competitive on-demand mobility systems}.
  \emph{Transportation Research Part C: Emerging Technologies} 91:192--208.

\end{thebibliography}
\newpage
\begin{APPENDICES}
\setcounter{table}{0}
\renewcommand{\thetable}{A\arabic{table}}
\section{Risk Measures Calculations}\label{sec:risk_measure_calculation}
We show in the following example the measure of the individual exposed risk of requests $1$ and $2$, denoted by $H_1$ or $H_2$, as illustrated in Example~\ref{example:risk_calculation}.

\begin{example}\label{example:risk_calculation} 
The example presents a route covering two passengers originating from nodes 1 and 2 and destined for nodes 3 and 4, respectively. Since we consider the same route, we will omit the $k$ dimension, and the variables $R_{i}$, $Q_{i}$, and $H_{1}$ and $H_{2}$ are associated with Constraints \eqref{eq:on_vehicle_risk_consistency}, \eqref{eq:risk_consistency}, and \eqref{eq:direct_infection}.

For convenience, we adopt the following assumptions: all travel times are fixed at 5 units ($t=5$), service times are assumed to be 0 units, and risk scores are uniformly set to 1 ($r_1=r_2=1$). All other constraints (e.g., time window, capacity, maximum ride time) are not binding. Specifically, at node 1, the total risk score is updated as $R_{1}\leftarrow r_{1}=1$. At node 2, $R_{2}\leftarrow R_{1}+r_{2}=1+1=2$, which is then updated to be $R_{3}\leftarrow R_{2}+r_{3}=2-1=1$ after visiting node 3. For the cumulative risk $Q_{i}$ along the route, $Q_{1}\leftarrow R_{1}\cdot t=5$, $Q_{2}\leftarrow Q_{1}+R_{2}\cdot t=15$, and $Q_{3}\leftarrow Q_{2}+R_{3}\cdot t=20$. The individual exposed risk for passenger 1 is obtained at node 3, which takes the form $H_{1}\leftarrow \left(Q_{3}-Q_{1}\right)-(A_{3}-A_{1})r_{1}=\left(Q_{3}-Q_{1}\right)-2t\cdot r_{1}=(20-5)-(2\times5)\cdot 1 = 15-10=5$. Similarly, the exposed risk for passenger 2 is expressed as $H_{2}\leftarrow \left(Q_{4}-Q_{2}\right)-(A_{4}-A_{2})r_{2} = (15-0)-10\cdot 1=5$.

\begin{figure}[H]
    \centering
    \caption{Example of risk measures calculation}
    \label{fig:risk_calculation}
    \includegraphics[width=\textwidth]{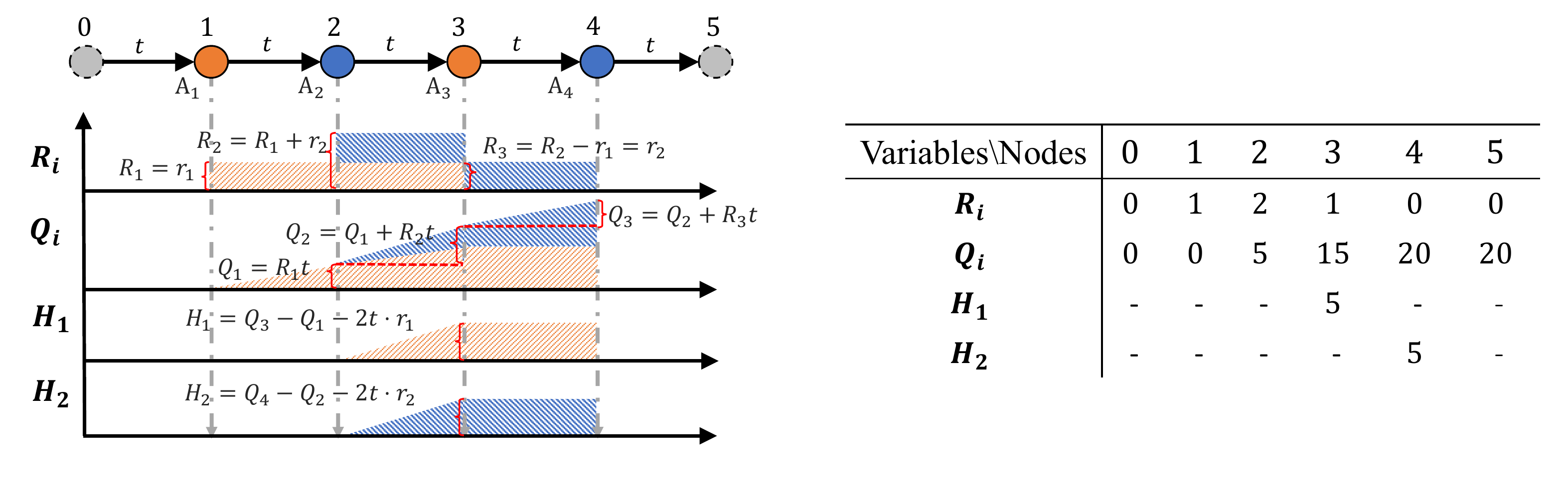}
\end{figure}
\end{example}

\section{Proof of Propositions}
\subsection{Proof of Proposition~\ref{prop:equivalence}}~\label{proof:equivalence}
\textbf{Proposition~\ref{prop:equivalence}} 
\textit{The arc-based and trip-based formulations are equivalent if every route $r\in\Omega$ satisfies the MMR condition.} 
\proof{Proof of Proposition~\ref{prop:equivalence}}
Let $\overline{H}^{\text{arc}}$ and $\overline{H}^{\text{trip}}$ denote the individual exposed risk that aligns with the MMR condition in the objective functions~\eqref{eq:bi_objective_function} and~\eqref{obj:original_master_problem} of the arc-based and trip-based formulations, respectively. In the arc-based formulation, we have the following observation according to Eq.~\eqref{eq:arc_min_max_risk}. 
\begin{subequations}
\begin{equation}
    \overline{H}^{\text{arc}}=\overline{H}=\max\limits_{i\in\mathcal{P}}\sum\limits_{k\in\mathcal{K}}H_{i}^{k}
\end{equation}

Suppose to the contrary that the arc-based and trip-based formulations are equivalent, i.e., $\overline{H}^{\text{arc}}=\overline{H}^{\text{trip}}=\overline{H}$. However, there exists a route $r\in\Omega$ that does not satisfy the MMR condition, where the maximum exposed risk on route $r$ is not minimized.

Without loss of generality, let us assume that there is at least one request $i$ on route $r$ associated with the maximum individual exposed risk among all routes, denoted by $\overline{H}_{ir}$. In the arc-based formulation, this maximum exposed risk is ensured to be minimized, i.e., $\overline{H}_{ir}=\overline{H}^{\text{arc}}$. However, in the trip-based formulation, while generating feasible columns, the exposed risk for request $i$ on route $r$, denoted by $H_{ir}$, is overestimated, leading to an even higher maximum exposed risk among all routes, i.e., $H_{ir}=\max\limits_{j\in\mathcal{P},r\in\Omega} H_{jr} = \overline{H}^{\text{trip}}$.

Therefore, we observe the following:
\begin{equation}
    \overline{H}^{\text{trip}}=H_{ir}>\overline{H}_{ir}=\overline{H}^{\text{arc}}
\end{equation}

This finding contradicts the assumption that $\overline{H}^{\text{arc}}=\overline{H}^{\text{trip}}$. Hence, we can conclude that the maximum exposed risk $\overline{H}$ must satisfy the MMR condition.
\qed
\end{subequations}

\subsection{Proof of Proposition~\ref{prop:R_H_non_increasing}}~\label{proof:R_H_non_increasing}
\textbf{Proposition~\ref{prop:R_H_non_increasing}} 
\textit{Both $\Delta \BFH_{\eta_{l},\eta_{l'}}$ and $\BFH_{\eta_{l'},\star}$ are  non-increasing function with respect to the increase in $\delta_{\eta_{l},\eta_{l'}}$. }
\proof{Proof of Proposition~\ref{prop:R_H_non_increasing}}
Suppose to the contrary that $\Delta \BFH_{\eta_{l},\eta_{l'}}$ increases with a higher value of delay $\delta_{\eta_{l},\eta_{l'}}$ at node $\eta_{l}$. This implies that when the delay at node $\eta_{l}$ increases, the co-ride time $\Delta A_{\eta_{l},\eta_{l'}}^{i}$ of one request $i\in\mathcal{O}_{l}$ also increases, resulting in a higher realized risk increment towards the maximum risk. The co-ride time $\Delta A_{\eta_{l},\eta_{l'}}$ can be written by:
\begin{equation}
    \Delta A_{\eta_{l},\eta_{l'}}^{i} = A_{l'}-(A_{l}+\delta_{\eta_{l},\eta_{l'}}) = t_{\eta_{l},\eta_{l'}}+s_{\eta_{l}} + \omega_{\eta_{l},\eta_{l'}}
\end{equation}
where the actual waiting time $\omega_{\eta_{l},\eta_{l'}}$ is the only adjustable variable in this context. Hence, the increased co-ride time is therefore caused by an increase in the waiting time, and equivalently, a lower level of delay. However, it contradicts the assumption that the increased exposed risk is due to a higher level of delay.

Therefore, we can conclude that the function $\Delta \BFH_{\eta_{l},\eta_{l'}}$ is non-increasing.

Next, we prove the non-increasing property of the function $\BFH_{\eta_{l'},\star}$. We can rewrite $\BFH_{\eta_{l'},\star}$ as below:
\begin{equation}
    \BFH_{\eta_{l'},\star} = \Delta \BFH_{\eta_{l'},\eta_{l''}}+\Delta \BFH_{\eta_{l''},\eta_{l'''}}+\ldots
\end{equation}
where $\eta_{l''}$ and $\eta_{l'''}$ represent the subsequent nodes to extend in sequence. Since we have already proven that the function $\Delta \BFH_{\eta_{l},\eta_{l'}}$ is non-increasing, we can conclude that the sum of non-increasing functions still preserves this property.
\qed

\subsection{Proof of Proposition~\ref{prop:optimality}}~\label{proof:optimality}
\textbf{Proposition~\ref{prop:optimality}} 
\textit{Let $\bar{\delta}_{\eta_{l},\eta_{l'}}$ be the minimum possible delay that minimizes $\Delta \BFH_{\eta_{l},\eta_{l'}}$, and $\delta_{\eta_{l},\eta_{l'}}^{\star}$ be one feasible delay that minimizes Eq.~\eqref{eq:bellman_mmr}. Then it holds true that $\delta_{\eta_{l},\eta_{l'}}^{\star}\geq \bar{\delta}_{\eta_{l},\eta_{l'}}$. Moreover, delaying only $\bar{\delta}_{\eta_{l},\eta_{l'}}$ at node $\eta_{l'}$ and delaying  $\delta_{\eta_{l},\eta_{l'}}^{\star}-\bar{\delta}_{\eta_{l},\eta_{l'}}$ at the subsequent node extended from $\eta_{l'}$ also minimize Eq.~\eqref{eq:bellman_mmr}.}
\proof{Proof of Proposition~\ref{prop:optimality}}

On arc $(\eta_{l},\eta_{l'})$, the potential waiting time at node $\eta_{l'}$, denoted by $\Delta \omega_{\eta_{l},\eta_{l'}}\geq 0$, is the only adjustable variable that can be offset to minimize the exposed risk increment $\Delta \BFH_{\eta_{l},\eta_{l'}}$. Requests tend to utilize the delay resources to minimize waiting time, approaching zero when sufficient delay resources are available, such that the arrival time at node $\eta_{l'}$ is exactly the earliest start-of-service time, denoted by $A_{l}+t_{\eta_l,\eta_{l'}}+s_{\eta_{l}}+\Bar{\delta}_{\eta_{l},\eta_{l'}}=a_{\eta_{l'}}$.

Another possibility to minimize Eq.~\eqref{eq:bellman_mmr} is by considering a further delay in the start-of-service at node $\eta_{l}$, resulting in an arrival time at node $\eta_{l'}$ given by $A_{l}+t_{\eta_l,\eta_{l'}}+s_{\eta_{l}}+\delta_{\eta_{l},\eta_{l'}}^{\star}>a_{\eta_{l'}}$. This delay can also be effective if subsequent nodes have waiting times that can be fully offset by the additional delay of $A_{\eta_{l'}}-a_{\eta_{l'}}$. 
In such cases, we observe $\delta_{\eta_{l},\eta_{l'}}^{\star}\geq \bar{\delta}_{\eta_{l},\eta_{l'}}=A{\eta_{l'}}-a_{\eta_{l'}}\geq 0$, indicating $\delta_{\eta_{l},\eta_{l'}}^{\star}\geq \bar{\delta}_{\eta_{l},\eta_{l'}}$.

However, excessive delay resources ($A_{l'}-a_{\eta_{l'}}>0$) eliminate possible extensions compared to $A_{l'}=a_{\eta_{l'}}$. Thus, waiting time should be minimized such that the arrival time at the extension node is \textit{exactly} the earliest time window ($A_{l'}=a_{\eta_{l'}}$).

We next prove that, delaying only $\bar{\delta}_{\eta_{l},\eta_{l'}}$ at node $\eta_{l'}$ and delaying  $\delta_{\eta_{l},\eta_{l'}}^{\star}-\bar{\delta}_{\eta_{l},\eta_{l'}}$ at the subsequent node extended from $\eta_{l'}$ also minimize Eq.~\eqref{eq:bellman_mmr}.

We consider the determinant of the amount of delay as a sequential decision process. Initially, we utilize the delay resource of $\Bar{\delta}_{\eta_{l},\eta_{l'}}$ to ensure $A_{l'}=a_{\eta_{l'}}$. Subsequently, if there exists waiting time at subsequent nodes, we utilize the remaining delay resource of $\delta_{\eta_{l},\eta_{l'}}^{\star}-\Bar{\delta}_{\eta_{l},\eta_{l'}}$ to maximize the offset of this waiting time. This approach, compared to directly delaying $\delta_{\eta_{l},\eta_{l'}}^{\star}$ at the current node $\eta_l$, maximizes the offset of all waiting times. Consequently, the potential exposed risk remains equal along the node sequence.
\qed

\subsection{Proof of Proposition~\ref{prop:dominance_rule_darp}}\label{proof:dominance_rule_darp}
\textbf{Proposition~\ref{prop:dominance_rule_darp}}
\textit{Given two RDARP-feasible labels $l_{1}$ and $l_{2}$ residing at the same node, i.e., $\eta_{1}=\eta_{2}\in\mathcal{N}$, $l_{1}$ is said to dominate $l_{2}$ if the following conditions are satisfied:}
\begin{align}
&A_{1}\leq A_{2}, W_{1}\leq W_{2},{Q_{1}\leq Q_{2}},\mathcal{V}_{1} \subseteq \mathcal{V}_{2},\mathcal{O}_{1}\subseteq \mathcal{O}_{2}, \mathcal{O}_{1}^{a} \subseteq \mathcal{O}_{2}^{a}~\label{appdx_prop:dominance_rule_darp0}\\
&\Tilde{c}_{1} \leq \Tilde{c}_{2}~\label{appdx_prop:dominance_rule_darp1}\\
&D^{o}_{1}(A_{1}) - A_{1}\geq D^{o}_{2}(A_{2}) - A_{2}, &&\forall o \in \mathcal{O}_{1}~\label{appdx_prop:dominance_rule_darp2}\\ 
&D^{o}_{1}(B^{o}_{1}) \geq D^{o}_{2}(B^{o}_{2}), &&\forall o \in \mathcal{O}_{1}~\label{appdx_prop:dominance_rule_darp3}\\
& d_{1}^{o} \geq d_{2}^{o}, &&\forall o \in \mathcal{O}_{1}^{a}~\label{appdx_prop:dominance_rule_darp4}
\end{align}

\proof{Proof of Proposition~\ref{prop:dominance_rule_darp}}
For simplicity, we differentiate the associated variables of $l_{1}$ and $l_{2}$ by using the subscripts $1$ and $2$.
For any feasible extension of node sequence $\mathcal{R}_{2}$ for $l_{2}$, we construct another extension $\mathcal{R}_{1}$ for $l_{1}$ by removing the drop-off nodes such that $\mathcal{R}_{1}:=\mathcal{R}_{2}\setminus\mathcal{D}$ while keeping the start-of-service time for node $i\in\mathcal{R}_{1}$. 
Given that $l_{2}$ is feasible, Eq.~\eqref{appdx_prop:dominance_rule_darp0} ensures the feasibility of $l_{1}$ with respect to the time window constraints, capacity constraints, pairing constraints, and precedence constraints. 
% By the assumption that $l_{1}$ and $l_{2}$ have the same feasible extension, 
{In this case, the set of possible extensions of $l_{2}$ is a subset of that of $l_{1}$, where $l_{1}$ indicates earlier start-of-service time, fewer load, less exposed risk, fewer visited passengers (and locations), suggesting more flexibility for the potential extensions.}
Eq.~\eqref{appdx_prop:dominance_rule_darp1} indicates that the reduced cost of $l_{1}$ cannot be greater than that of $l_{2}$. Note that the terms $\Tilde{c}_{l_{1}}$ or $\Tilde{c}_{l_{2}}$ includes the minimized incremental risk aggregated by arcs, factored by $h_{l_{1'}}^{o}-h_{l_{1}}^{o}$ or $h_{l_{2'}}^{o}-h_{l_{2}}^{o}$ (see Eq.~\eqref{REF:reduced_cost}), which is proved to have no conflict with the future extensions. 
Eqs.~\eqref{appdx_prop:dominance_rule_darp2} and~\eqref{appdx_prop:dominance_rule_darp3} guarantee that $D^{o}_{1}(\theta)\geq D^{o}_{2}(\theta), \forall o\in\mathcal{O}_{1}, \theta\in[A_{\eta_{l_{2}}}, b_{\eta_{l_{2}}}]$, which implies that it is feasible to deliver passenger $o$ of $l_{1}$ whenever it is feasible for $l_{2}$ (see Proposition 3, 4, and 5 in~\cite{gschwind2015effective}). Finally, we introduce a risk-related dominance rule, denoted as Eq.~\eqref{appdx_prop:dominance_rule_darp4}, which extends the existing DARP dominance conditions. In this rule, we require that $\mathcal{O}_1^a \subseteq \mathcal{O}_2^a$. This rule signifies that label $l^1$ provides more available delay buffer $d_{1}^{o}\geq d_{2}^{o}$ for requests $o \in \mathcal{O}_1 \cup \mathcal{O}_1^a$ to delay their service, allowing for a larger amount of time that can be delayed for each open and associated request $o \in \mathcal{O}_1 \cup \mathcal{O}_1^a$ in future extensions. Consequently, label $l^1$ possesses more flexibility to delay its service and a higher potential to achieve a lower MMR value during label extension compared to label $l^2$.
\qed
\endproof

\subsection{Proof of Proposition~\ref{prop:rdarp_edarp}}\label{proof:rdarp_edarp}
\textbf{Proposition~\ref{prop:rdarp_edarp}}: \textit{Under the EDARP setting, the measure of individual exposed risk $H_{i}$ for a real passenger $i\in\mathcal{P}$ is equivalent to the onboard duration $A_{n+i}-A_{i}$.}
\proof{Proof of Proposition~\ref{prop:rdarp_edarp}} 
We prove by showing the equivalency between exposed risk $H_i$ and contact duration with the dummy passenger $p_0$. Specifically, the calculation of $H_i$ can be expanded based on Constraints~\eqref{eq:direct_infection} as follows (we tentatively omit the dimension $k$).
\begin{align*}
    H_{i} &= Q_{n+i} - Q_{i} - \left(A_{n+i}-A_{i}\right)r_{i}\\
    & = \sum_{(u,v)\in r=\{i,\ldots,n+i\}} R_{u} \left(A_{u}-A_{v}\right) - \left(A_{n+i}-A_{i}\right)\cdot 0\\
    & = \sum_{(u,v)\in r=\{i,\ldots,n+i\}} r_{0} \left(A_{u}-A_{v}\right)\\
    & = A_{n+i} - A_{i}
\end{align*}
This completes the proof.\qed

\section{Exact $\epsilon$-constraint method}\label{sec:Pareto_front_approximation}
We adopt the $\epsilon$-constraint method (ECM) to capture the Pareto front of our RDARP. 
To ensure the exact Pareto optimal solutions~\citep{berube2009exact}, we need to iteratively solve the RMPs $\BFP^{\text{cost}}(\epsilon^{\text{risk}})$ and $\BFP^{\text{risk}}(\epsilon^{\text{cost}})$. As noted by \citet{berube2009exact}, there may exist multiple solutions by solving one single-objective problem (e.g., $\BFP^{\text{cost}}(\epsilon^{\text{risk}})$), where some dominated solutions may be generated. However, by additionally solving the other single-objective problem (e.g., $\BFP^{\text{risk}}(\epsilon^{\text{cost}})$) with a fixed and restrictive cap $\epsilon^{\text{cost}}$), it will lead to the unique Pareto-optimal solution, thus ensuring the exactness of our iterative approach.

The iterative algorithm proceeds as follows. We first express the bi-objective functions as $f^{\text{cost}}$ (total travel cost) and $f^{\text{risk}}$ (maximum risk). Let the Ideal points be $(f_{\text{cost}}^{I},f_{\text{risk}}^{I})$ and Nadir points be $(f_{\text{cost}}^{N},f_{\text{risk}}^{N})$. Starting with the DARP solution $(f_{\text{cost}}^{I}, f_{\text{risk}}^{N})$, we iteratively solve $\BFP^{\text{cost}}(\epsilon^{\text{risk}})$ with the updated $\epsilon^{\text{risk}}$. And the exact solution is further validated by solving $\BFP^{\text{risk}}(\epsilon^{\text{cost}})$. The set of Pareto optimal solutions is considered as the Pareto front, denoted by $\mathcal{F}:=\{(f_{\text{cost}}^{i},f_{\text{risk}}^{i}): i \in 1,2,3,\ldots\}$. A detailed iterative approach to approximate the Pareto front is shown in Algorithm~\ref{alg:Pareto_front_approximation}.

In Algorithm~\ref{alg:Pareto_front_approximation}, we start by solving a DARP case, and the by-product of maximum risk $\overline{H}$ is checked by solving $\BFP^{\text{risk}}(\epsilon^{\text{cost}})$. For each iteration, we update the upper bound of maximum risk with a step size of $\Delta \epsilon^{\text{risk}}$. The iterative procedure terminates if the upper bound is updated to be a negative number or no feasible solution can be obtained under a sufficiently restrictive upper bound $\epsilon^{\text{risk}}$. Finally, we report the set of Pareto-optimal solutions as a Pareto frontier $\mathcal{F}$.

\begin{algorithm}[H]
\footnotesize%\small, \scriptsize, or \tiny
\caption{Pareto front approximation}
\label{alg:Pareto_front_approximation}
\textbf{Input:} Step size for max-risk term $\Delta \epsilon^{\text{risk}}$, an empty set of Pareto front $\mathcal{F}=\{\}$\\
\textbf{Output:} $\mathcal{F}:=\{(f_{\text{cost}}^{i},f_{\text{risk}}^{i}): i = \mathbb{N}\}$.
\begin{algorithmic}[1] %[1] enables line numbers
    \State $i=1$, $\epsilon^{\text{risk}}=\infty$
    \While{$\epsilon^{\text{risk}}-\Delta \epsilon^{\text{risk}}\ge 0 \& \BFP^{\text{cost}}(\epsilon^{\text{risk}})$ is feasible}
        \State \textbf{Solve} the $\BFP^{\text{cost}}(\epsilon^{\text{risk}})$ with $\epsilon^{\text{risk}}$ and \textbf{obtain} $f_{\text{cost}}^{i+1}$ \Comment{Solve the RDARP under an adjusted $\epsilon^{\text{risk}}$}
        \State \textbf{Solve} the $\BFP^{\text{risk}}(\epsilon^{\text{cost}})$ with $\epsilon_{t} = f_{\text{cost}}^{i+1}$ and \textbf{obtain} $f_{\text{risk}}^{i+1}$ \Comment{Guarantee the min-max exposed risk under the optimal travel cost $f_{\text{cost}}^{i+1}$.}
        \State $\epsilon^{\text{risk}} \leftarrow f_{\text{risk}}^{i+1}-\Delta \epsilon^{\text{risk}}$
        \State $\mathcal{F}\leftarrow \mathcal{F}\bigcup\{(f_{\text{cost}}^{i+1},f_{\text{risk}}^{i+1})\}$
    \EndWhile
\end{algorithmic}
\end{algorithm}
\begin{figure}[H]
    \centering
    \caption{An illustration of Pareto front}\label{fig:illustration_pareto}
    \includegraphics[width=0.5\textwidth]{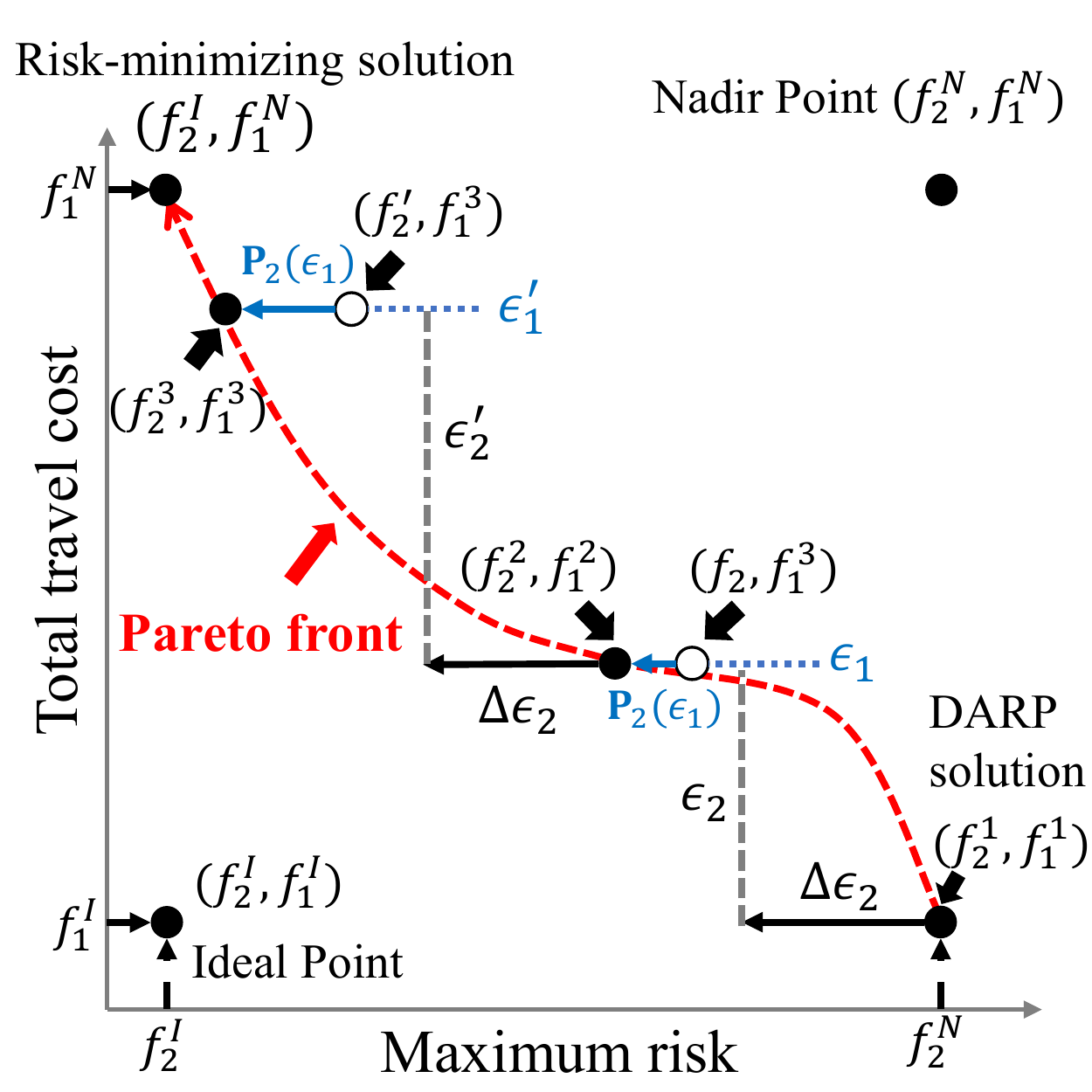}
\end{figure}
An illustration of Algorithm~\ref{alg:Pareto_front_approximation} is shown in Figure~\ref{fig:illustration_pareto}, where we use the subscripts $1$ and $2$ in replace of $cost$ and $risk$ for convenience of presentation. In particular, the red dashed curve shows the Pareto front, the blue arrows indicate that better Pareto-optimal solutions are founded by solving $\BFP^{\text{risk}}(\epsilon^{\text{cost}})$, and the black arrows represent the step size of $\epsilon^{\text{risk}}$.

\section{Description of ride time related REF}\label{appendix:description_ride_time_resources}
\textbf{Initiating:} Without loss of generality, let $o_1$ be a new request at the node $\eta_{l'}$. 
For $o_1$, the latest possible drop-off time with respect to the earliest start-of-service time at node $\eta_{l'}$, $D_{l'}^{o_1}(A_{l'})$, is initiated, which is bounded by either $L_{\max}^{o_{1}}$ or $b_{\eta_{l'}+n}$ (see Eq.~\eqref{REF:ldT}). Similarly, the latest possible drop-off time, {given the breaking point $B_{l'}^{o_1}$ at $\eta_{l'}$}, can be expressed as $D_{l'}^{o_1}(B_{l'}^{o_1})$, which is initiated by the latest drop-off time (see Eqs.~\eqref{REF:B} and~\eqref{REF:ldB}). And it defines the \textit{absolute} latest possible drop-off time. Therefore, $\left[D_{l'}^{o_1}(A_{l'}),D_{l'}^{o_1}(B_{l'}^{o_1})\right]$ defines the DTW for $o_1$'s drop-off node. Note that the latest possible drop-off time $D_{l'}^{o_1}(t_{o_1})$ can be understood as a linear piece-wised function with respect to the start-of-service time $t_{o_1}$ at $o_1$'s pick-up node, which is bounded by the absolute latest drop-off time (i.e., latest drop-off time{$b_{o_1+n}$}). Specifically, given an arbitrary start-of-service time $t_{o_1}$, $D_{l'}^{o_1}(t)$ increases linearly as $t_{o_1}$ delays if $D_{l'}^{o_1}(t)\in [D_{l'}^{o_1}(A_{l'}),D_{l'}^{o_1}(B_{l'}^{o_1})]$. And we particularly denote the breaking point $B_{l'}^{o_1}$ by the earliest pick-up time $t_{o_1}$ such that for any $t_{o_1}\geq B_{l'}^{o_1}$, we have $D_{l'}^{o_1}(t)=D_{l'}^{o_1}(B_{l'}^{o_1})$. It describes the case that any delay at $o_{1}$ will not further extend the latest possible drop-off time due to the latest time window constraints ($b_{o_{1}}$ and $b_{o_{1}+n}$).

\textbf{Updating:} The REF for updating the existing open requests follows a similar procedure. Let $o_2\in\mathcal{O}_{l}$ be an existing open request. The DTW can be further bounded by the time windows implied by the nodes on the transitive closure of the partial route $\{o_{2},\ldots,o_2+n\}$. Consider a partial route $\{o_{2},\ldots,\eta_{l},\eta_{l'}\}$ with $o_{2}\in\mathcal{O}_{l'}$. To properly update the range of DTW for $o_2$, we first update the breaking point $B_{l'}^{o_2}$, representing the earliest start-of-service time at $\eta_{l'}$ such that the $o_2$'s latest possible drop-off time is exactly $D_{l'}^{o_2}(B_{l'}^{o_2})$. Then we consider three cases regarding the relationship among $B_{l}^{o_2}+s_{\eta_{l}}+t_{\eta_{l},\eta_{l'}}, A_{l'}$, and $B_{l'}$ (see the REF in Eq.~\eqref{REF:B}). 
\begin{enumerate}
    \item If $B_{l}^{o_2}+s_{\eta_{l}}+t_{\eta_{l},\eta_{l'}}<A_{l'}$, then $B_{l'}^{o_2}=A_{l'}$. It implies that any start-of-service time at $\eta_{l'}$ is beyond $B_{l'}^{o_2}$, resulting in  $D_{l'}^{o_2}(A_{l'})=D_{l'}^{o_2}(B_{l'}^{o_2})$.
    \item If $B_{l}^{o_2}+s_{\eta_{l}}+t_{\eta_{l},\eta_{l'}}\in [A_{l'},B_{l'}]$, then $D_{l'}^{o_2}(A_{l'})$ increases linearly to the amount of potential waiting time on the arc $(\eta_{l},\eta_{l'})$, denoted by $\max\{A_{l'}-s_{\eta_{l}}-t_{\eta_{l},\eta_{l'}},0\}$ (see Eq.\eqref{REF:ldT}).
    \item If $B_{l}^{o_2}+s_{\eta_{l}}+t_{\eta_{l},\eta_{l'}} > B_{l'}$, $B_{l'}^{o_2}$ is adjusted to be $B_{l'}$ to ensure the time window feasibility. Accordingly, the absolute latest possible drop-off time is given by: 
    \begin{equation*}
        D_{l'}^{o_2}(B_{l'}^{o_2})=D_{l'}^{o_2}\left(B_{l}^{o_2}+s_{\eta_{l}}+t_{\eta_{l},\eta_{l'}}-\left(B_{l}^{o_2}+s_{\eta_{l}}+t_{\eta_{l},\eta_{l'}}-B_{l'}\right)\right)=D_{l'}^{o_2}(B_{l'})
    \end{equation*} 
    Therefore, the updated $D_{l'}^{o_2}\left(B_{l}^{o_2}\right)$ can be understood as the latest possible drop-off time if one starts the service at $B_{l'}$, see Eq.\eqref{REF:ldB}.
\end{enumerate}

\section{Description of Algorithm~\ref{alg:contact_duration_calibriation} on an interlaced route}\label{appendix:risk_REF_update}

Figure~\ref{fig:example_REF} shows an example of the REF updates on an interlaced partial route $(0,i,j,i+n,k,j+n,\ldots)$. For convenience, we set all travel time as 10, service time as 0, and maximum ride time as $L_{\max}^{i}=20,L_{\max}^{j}=L_{\max}^{k}=40$, and all other constraints (e.g., time, capacity, and cumulative risk constraints) are not binding. 

The resources are updated as follows. Set $\mathcal{O}_{l}$ records the new requests $i,j,k$ and removes the element $o\in\mathcal{O}_{l}$ after visiting $o+n$. As all passengers are interlaced, no passengers will be removed from the set of associated requests $\mathcal{O}_{l}^{a}$. 
The earliest possible start-of-service time is updated based on Eq.~\eqref{REF:early_time}. All the latest possible delivery time $B_{l}$ are bounded by $b_{\eta_{l}}$ except the node $j+n$. At node $\eta_{l}=j+n$, $B_{l}$ is bounded by request $i$'s maximum ride time such that the latest start-of-service time at node $j$ is $B_{l}^{j}=30$, leading to the latest possible delivery time $D_{l}^{j}(B_{l}^{j})=B_{l}^{j}+L_{\max}^{j}=30+40=70$. The maximum amount of delay is measured by the time gap between the earliest and latest possible start-of-service time $B_{l}-A_{l}$, updated for each label extension at node $\eta_{l'}=i\in\mathcal{N}$.

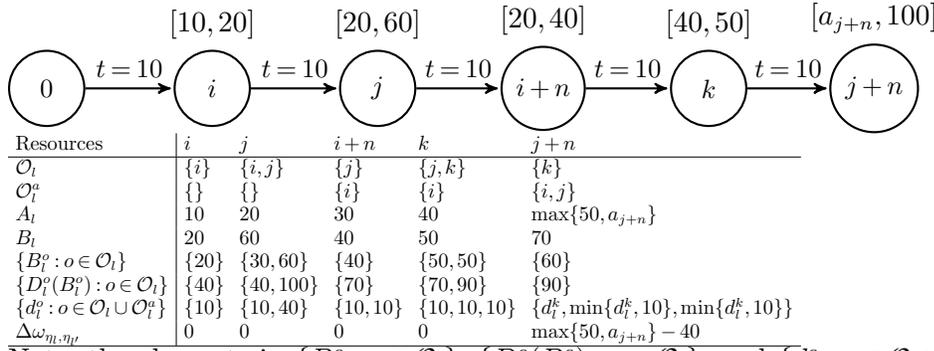
\begin{figure}[H]
    \caption{An illustration of the path $(0,i,j,i+n,k,j+n,\ldots)$}
    \label{fig:example_REF}
    \begin{tikzpicture}[->,>=stealth',shorten >=1pt,auto,node distance=2.2cm,  thick,main node/.style={circle,draw,minimum size=1cm,font=\sffamily\small\bfseries}]
      \node[main node] (1) {$0$};
      \node[main node, label=above:{$[10,20]$}] (2) [right of=1] {$i$};
      \node[main node, label=above:{$[20,60]$}] (3) [right of=2] {$j$};
      \node[main node, label=above:{$[20,40]$}] (4) [right of=3] {$i+n$};
      \node[main node, label=above:{$[40,50]$}] (5) [right of=4] {$k$};
      \node[main node, label=above:{$[a_{j+n},100]$}] (6) [right of=5] {$j+n$};
      \path[every node/.style={font=\sffamily\small}]
        (1) edge node [above]  {$t=10$} (2)
        (2) edge node [above]  {$t=10$} (3)
        (3) edge node [above]  {$t=10$} (4)
        (4) edge node [above]  {$t=10$} (5)
        (5) edge node [above]  {$t=10$} (6);
    \end{tikzpicture}
    \scalebox{0.7}{\begin{tabular}{l|lllllll}
        Resources & $i$ & $j$ & $i+n$ & $k$ & $j+n$ \\
        \hline
        $\mathcal{O}_{l}$ & $\{i\}$ & $\{i,j\}$& $\{j\}$& $\{j,k\}$& $\{k\}$\\
        $\mathcal{O}_{l}^{a}$ & $\{\}$ & $\{\}$& $\{i\}$& $\{i\}$& $\{i,j\}$\\
        $A_{l}$ & 10 & 20 & 30 & 40 & $\max\{50,a_{j+n}\}$ \\
        $B_{l}$ & 20 & 60 & 40 & 50 & 70 \\
        $\{B_{l}^{o}:o\in\mathcal{O}_{l}\}$ & $\{20\}$ & $\{30,60\}$ & $\{40\}$ & $\{50,50\}$ & $\{60\}$ \\
        $\{D_{l}^{o}(B_{l}^{o}):o\in\mathcal{O}_{l}\}$& $\{40\}$ & $\{40,100\}$ & $\{70\}$ & $\{70,90\}$ & $\{90\}$ \\
        $\{d_{l}^{o}:o\in\mathcal{O}_{l}\cup\mathcal{O}_{l}^{a}\}$ & $\{10\}$ & $\{10,40\}$ & $\{10,10\}$ & $\{10,10,10\}$ & $\{d_{l}^{k},\min\{d_{l}^{k},10\},\min\{d_{l}^{k},10\}\}$ \\
        $\Delta \omega_{\eta_l,\eta_{l'}}$ & 0& 0& 0& 0& $\max\{50,a_{j+n}\}-40$\\
        \hline
    \end{tabular}
    }
    \\
    Note: the elements in $\{B_{l}^{o}:o\in\mathcal{O}_{l}\}$, $\{D_{l}^{o}(B_{l}^{o}):o\in\mathcal{O}_{l}\}$, and $\{d_{l}^{o}:o\in\mathcal{O}_{l}\cup\mathcal{O}_{l}^{a}\}$ maintain the same order as the elements in the set $O_l$ or $O_l^a$.
     
\end{figure}

In Figure~\ref{fig:algo2_3conditions}, we aim to calibrate the onboard time $\Delta A_{k,j+n}$ on arc $(k,j+n)=(\eta_{l},\eta_{l'})$ under two scenarios, where the earliest start-of-service time $a_{j+n}$ is set to be $60$ ($\Delta d_{k,j+n}^{o}\ge\Delta\omega_{k,j+n}$ and the waiting time can be completely offset), and $65$ ($\Delta d_{k,j+n}^{o}<\Delta\omega_{k,j+n}$ and the waiting time can be partially offset without violating other constraints), respectively. The maximum allowable delay for nodes $i$ and $k$ in all cases is $10$. In the case of arcs $(i,j)$, the maximum amount of delay for node $j$ is $40$, which is further constrained by the drop-off node $i+n$ with $B_{l'}^{j}-A_{l'}=40-30=10$. The calculations of the minimum onboard duration for the three cases are presented below.

\begin{itemize}
    \item[\textbf{Case (1):}] if $a_{j+n}=60$, $A_{l'} = \max\{a_{k}+t_{k,j+n}+s_{k},a_{j+n}\}=60$. The waiting time is $\Delta \omega_{\eta_{l},\eta_{l'}}=\min\{A_{l'}-s_{\eta_{l}}-t_{\eta_{l},\eta_{l'}},b_{\eta_{l}}\}-A_{l}=\min\{60-10-0,50\}-40 = 50-40=10$, which can be completely offset by a delay of $\Delta d_{k,j+n}^{o}=10$ for all $o\in\mathcal{O}_{l'}^{a}=\{i,j,k\}$. Hence, $\Delta A_{\eta_{l},\eta_{l'}}^{o} = t_{\eta_{l},\eta_{l'}}+s_{\eta_{l}}=10$.
    \item[\textbf{Case (2):}] if $a_{j+n}=65$, $A_{l'} = \max\{a_{k}+t_{k,j+n}+s_{k},a_{j+n}\}=65$, and $\Delta \omega_{\eta_{l},\eta_{l'}}=\min\{A_{l'}-s_{\eta_{l}}-t_{\eta_{l},\eta_{l'}},b_{\eta_{l}}\}-A_{l}=\min\{65-10-0,50\}-40 = 55-40=15$. In this case, for every associated request $o\in\mathcal{O}_{l'}^{a}$, the waiting time can be reduced to $\Delta \omega_{\eta_{l},\eta_{l'}}-\Delta d_{k,j+n}^{o}=15-10=5$ by delaying $\Delta d_{k,j+n}^{o}=10$. Hence, the minimal possible onboard duration is $\Delta A_{\eta_{l},\eta_{l'}}^{o}=t_{ij}+s_i+\Delta \omega_{\eta_{l},\eta_{l'}}-\Delta d_{k,j+n}^{o}=15$.
\end{itemize}      

\begin{figure}[H]
    \centering
    \caption{Comparison between delay buffer $\Delta d_{\eta_{l},\eta_{l'}}$ and waiting time $\Delta \omega_{\eta_{l},\eta_{l'}}$}
    \label{fig:algo2_3conditions}
    \begin{minipage}[c]{0.45\textwidth}
        \includegraphics[width=\textwidth]{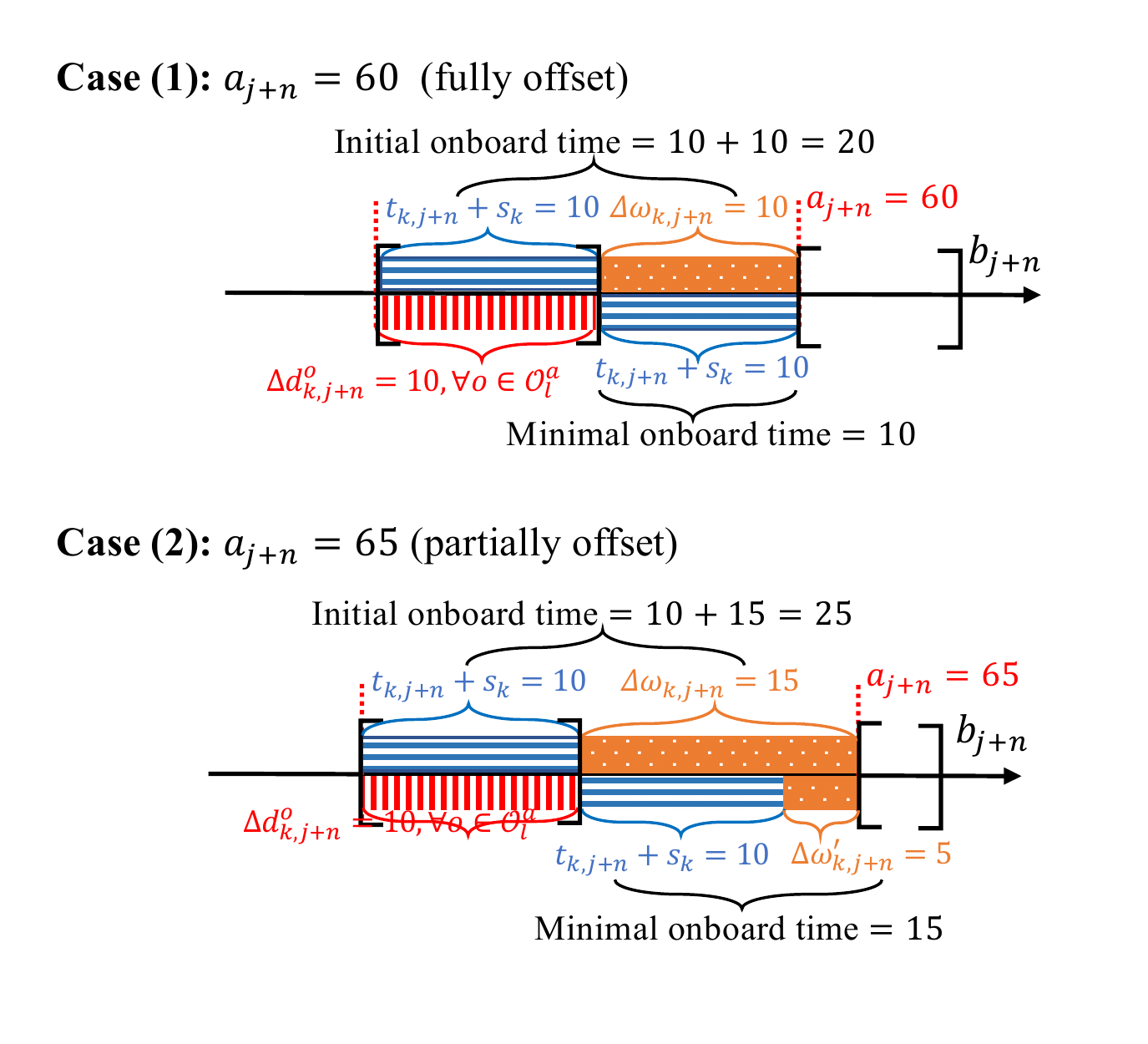} % 3conditions
    \end{minipage}
    \begin{minipage}[c]{0.48\textwidth}
    \centering
        \begin{tabular}{l|ccc} \hline
        Variables $\backslash$ Cases & \textbf{(1)} & \textbf{(2)} \\\hline
        $t_{ij}+s_{i}$ & \multicolumn{2}{c}{$10+0=10$}\\\hline
        $\mathcal{O}_{l}$ & \multicolumn{2}{c}{$\{j,k\}$}\\\hline
        $\mathcal{O}_{l}^{a}$ & \multicolumn{2}{c}{$\{i\}$}\\\hline
        $\{d_{l}^{i},d_{l}^{j},d_{l}^{k}\}$ & \multicolumn{2}{c}{$\{10,10,10\}$} \\\hline
        $\Delta \omega_{k,j+n}$ &  $10$ & $15$\\\hline
        $\max\{\Delta \omega_{k,j+n}-\Delta d_{k,j+n}^{o} ,0\}$& 0 & 5\\\hline
        $\Delta A_{k,j+n}$ & $10$ & $15$\\\hline
        \end{tabular}
    \end{minipage}
\end{figure}

\section{Characteristic function for associated requests}\label{sec:char_risk_associated}

We introduce the function $\BFG_{l}^{i}(\delta)$ to model the piece-wise increase of exposed risk based on the service order and the risk scores and the maximum amount of delay of co-riders. Consider an associated request $i\in\mathcal{O}_{l}^{a}$ with $n-1$ co-riders denoted as $j_1, j_2, \dots, j_{n-1}$ who were exposed to the associated request $i$ and are still onboard, where ${j_1, j_2, \dots, j_{n-1}}:=\{j\in\mathcal{O}_{l}:\Gamma_{l}^{i}<\Gamma_{l}^{j}<\Gamma_{l}^{i+n}\}$. In this scenario, the function $\BFG_{l}^{o}(\delta)$ is defined as a piecewise linear function with $n-1$ linear segments. Each segment represents the maximum amount of delay for a specific co-rider between the pick-up node $i$ and drop-off node $i+n$, and the slope of the segment corresponds to the total risk score of all associated and open requests that precede that co-rider.

\begin{equation}
    \frac{\partial \BFG_{l}^{i}(\delta)}{\partial\delta}:=\begin{cases}
        0 & \text{if~} \delta \leq d_{l}^{i}\\
        \left(\sum\limits_{o\in\mathcal{O}_{l}\cup\mathcal{O}_{l}^{a},\Gamma_{l}^{o}<\Gamma_{l}^{j_{1}}<\Gamma_{l}^{o+n}} r_{o}\right)-r_{i} & \text{if~}d_{l}^{i}<\delta \leq d_{l}^{j_1}\\
        \left(\sum\limits_{o\in\mathcal{O}_{l}\cup\mathcal{O}_{l}^{a},\Gamma_{l}^{o}<\Gamma_{l}^{j_{2}}<\Gamma_{l}^{o+n}} r_{o}\right)-r_{i} & \text{if~}d_{l}^{j_1}<\delta\leq d_{l}^{j_2}\\
        \ldots& \\
        \left(\sum\limits_{o\in\mathcal{O}_{l}\cup\mathcal{O}_{l}^{a},\Gamma_{l}^{o}<\Gamma_{l}^{j_{N}}<\Gamma_{l}^{o+n}} r_{o}\right)-r_{i} & \text{if~}d_{l}^{j_{n-2}}< \delta \leq d_{l}^{j_{n-1}}\\
    \end{cases}
\end{equation}

\section{Detour Rate Calculations}\label{appendix:detour_rate_calculation}

\begin{example}\label{example:edarp_calculation}
To better understand the equivalency, we illustrate an EDARP-version example of risk measure in Figure~\ref{fig:edarp_calculation}, where $t=5,s_i=0,r_0=1,r_1=r_2=-r_3=-r_4=0$. Following similar logic as in Example~\ref{example:risk_calculation}, we report that $H_{1}=A_{3}-A_{1}=10$ and $H_{2}=A_{4}-A_{2}=10$.
\end{example}
\begin{figure}[H]
    \centering
    \caption{Example of risk measure in EDARP}
    \label{fig:edarp_calculation}
    \includegraphics[width=0.9\textwidth]{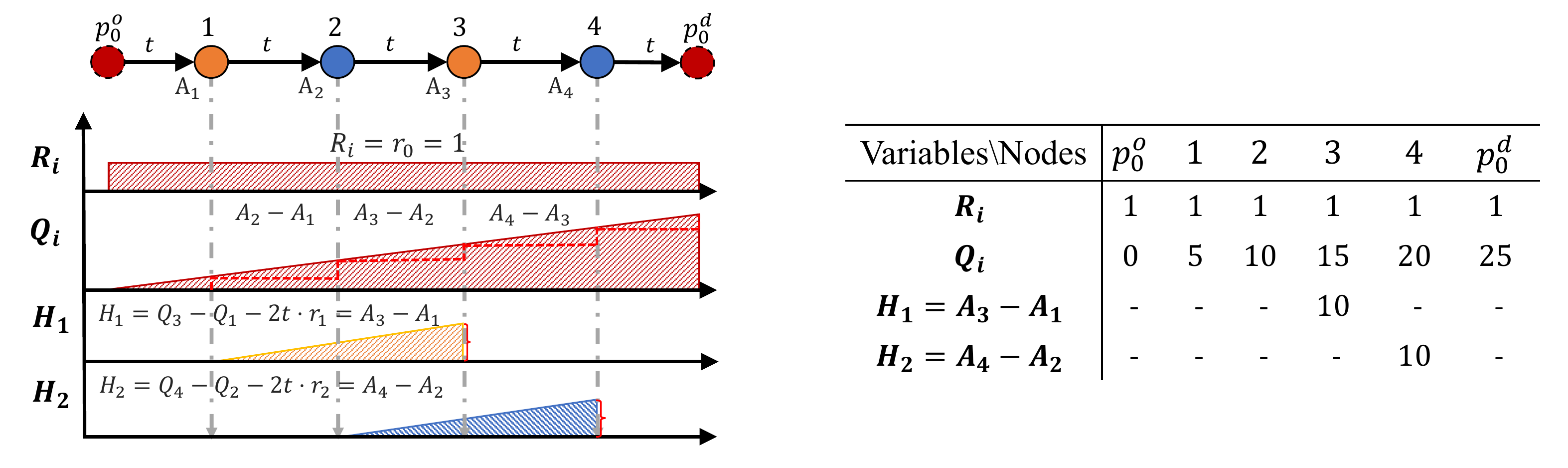}
\end{figure}

\section{Valid inequalities}\label{appendix:valid_inequality}

\subsection{Infeasible path elimination constraints}
We first introduce the IPEC, which is also known as the tournament constraints. For each RDARP-infeasible (partial) trip $I\in\mathcal{I}$, the IPEC are enforced to RMP, which takes the form:

\begin{equation}
    \sum_{r\in\Omega'}\sum_{(i,j)\in r} \beta_{ij,I}{\lambda_{r}}\le |I| - 1
\end{equation}
where parameter $\beta_{ij,I}$ denotes the number of times that arc $(i,j)$ is traversed by $I$, and $|I|$ is the number of arcs in the (partial) trip $I$. Moreover, \cite{cordeau2006branch} suggested a strengthened version of IPEC, which requires that the partial path $r$ starts with a pick-up node $i$ and ends up with $n+i$. It implies that for the particular passenger $i$, both $i$ and $n+i$ cannot be connected with the transitive closure between $i$ and $n+i$. In this case, the strengthened IPEC can be expressed as:
\begin{equation}
    \sum_{r\in\Omega'}\sum_{(i,j)\in r} \beta_{ij,I_s}{\lambda_{r}}\le |I_s| - 2
\end{equation}
where $I_s\in\mathcal{I}$ denotes the infeasible paths starting with $i\in\mathcal{P}$ and ending with $n+i$.

{The separation heuristics of infeasible path elimination constraints (IPEC) are conducted under an enumeration procedure (see Appendix C in \citet{gschwind2015effective}).} 

\subsection{Two-path inequality}
Two-path inequality was first introduced by \cite{kohl19992} to solve the VRPTW, and was later adopted to solve the PDPTW~\citep{ropke2009branch}. The procedure is detailed as follows. Consider a subset of nodes $\mathcal{S}\subset \mathcal{P}\bigcup\mathcal{D}$, and the nodes in $\mathcal{S}$ cannot be served by one vehicle on the same trip. The two-path valid inequality can be expressed as:

\begin{equation}
    \sum_{r\in\Omega'}\sum_{i\in \mathcal{S}} \sum_{j\in \mathcal{N} \setminus \mathcal{S}} \beta_{ij,r}\lambda_{r} \geq 2\label{eq:2path_inequality}
\end{equation}
where parameter $\beta_{ij,r}$ is the number of times that route $r$ traverses arc $(i,j)\in\mathcal{A}$. 

{The separation heuristic regarding the 2-path cuts are introduced in Section 6.2.2 in \citet{ropke2009branch}. 
For our implementation, we truncate the size of candidate sets $\mathcal{S}$ to inspect, where $|\mathcal{S}|\le 4$. In this case, some long sequences of nodes may be omitted. However, we consider that most infeasible long sequences are consisting of infeasible short segments, which can be well covered by a subset of $|\mathcal{S}|\le 4$.}

\subsection{Rounded capacity inequality}
Finally, we introduce the rounded capacity inequalities. The rounded capacity constraints propose a lower bound of the number of times that one vehicle must enter or leave one set $\mathcal{S}$. Let the predecessors and successors of $\mathcal{S}$ denote by $\sigma^{-}(\mathcal{S})=\{i\in\mathcal{P}\setminus \mathcal{S}:n+i\in \mathcal{S}\}$ and $\sigma^{+}(\mathcal{S})=\{n+i\in\mathcal{D}\setminus \mathcal{S}:i\in \mathcal{S}\}$, respectively. The rounded capacity is given by:

\begin{equation}
    \sum_{r\in\Omega'}\sum_{i\in \mathcal{S}} \sum_{j\in \mathcal{N} \setminus \mathcal{S}} \beta_{ij,r}\lambda_{r} \geq \max\{1,\lceil \frac{W\left(\sigma^{-}(\mathcal{S})\right)}{W_{\max}} \rceil,\lceil - \frac{W\left(\sigma^{+}(\mathcal{S})\right)}{W_{\max}} \rceil \}\label{eq:RC_inequality}
\end{equation}
where $W\left(\sigma^{-}(\mathcal{S})\right)$ is the total load of the set of nodes in $\sigma^{-}(\mathcal{S})$, denoted by $\sum\limits_{i\in \sigma^{-}(\mathcal{S})} w_{i}$, and similarly, $W\left(\sigma^{+}(\mathcal{S})\right)=\sum\limits_{i\in \sigma^{+}(\mathcal{S})} w_{i}$.

{For the rounded capacity constraint, we adopt the constructive heuristic as described in Section 6.2.1 in \citet{ropke2009branch}.}

\subsection{Handling Dual Variables in Branching Strategies}\label{appendix:branching_cuts}

{We adopt two types of branching strategies: (1) branching on number of vehicles and (2) branching on the outflow of node pair.

Specifically, we first branch on the number of vehicles, calculated as $\sum\limits_{r\in\Omega'}\Bar{\lambda_{r}}$, where $\Bar{\lambda_{r}}$ denote the LP solution. If the number of vehicles is fractional, the two branches can be expressed as: 

\begin{equation}\label{eq:branching_num_veh}
    \begin{cases}
    \sum\limits_{r\in\Omega'}{\lambda_{r}} \le \lfloor \sum\limits_{r\in\Omega'}\Bar{\lambda_{r}} \rfloor\\
    \sum\limits_{r\in\Omega'}{\lambda_{r}} \ge \lceil \sum\limits_{r\in\Omega'}\Bar{\lambda_{r}} \rceil
    \end{cases}
\end{equation}
which follows similar forms as Constraints~\eqref{cons:fleet_size_constraint}. Let $\mu^{-}\le0$ and $\mu^{+}\ge0$ be the associated dual variables in Eq.~\eqref{eq:branching_num_veh}. The reduced cost function is then updated as follows:
\begin{equation}
    \begin{cases}
        \tilde{c}_{r}'\leftarrow\tilde{c}_{r}-\mu^{+}\\
        \tilde{c}_{r}'\leftarrow\tilde{c}_{r}-\mu^{-}
    \end{cases}
\end{equation}
where $\tilde{c}_{r}'$ is the updated reduced cost function and $\tilde{c}_{r}$ can be $\Tilde{c}_{r}^{\text{cost}}$ or $\Tilde{c}_{r}^{\text{risk}}$ depending on the specific problem.

The second branching strategy is to branch on the node pair $\Lambda$, denoted by $
\delta^{+}\left(\Lambda\right)=\{(i,j)\in\mathcal{A}:i\in \Lambda, j\in \mathcal{N}\setminus\Lambda\}$. Let $\Bar{x}_{ij}$ be an LP solution. The set $\Lambda$ is of priority to be selected if the outflow $\sum\limits_{(i,j)\in \delta^{+}(\Lambda)}\Bar{x}_{ij}$ is closest to 1.5, and two branches are presented as follows:
\begin{equation}
    \begin{cases}
        \sum\limits_{(i,j)\in \delta^{+}(\Lambda)}{x}_{ij}\le 1\\
        \sum\limits_{(i,j)\in \delta^{+}(\Lambda)}{x}_{ij}\ge 2
    \end{cases}
\end{equation}
To align with the structure in the trip-based formulation, we can translate the equation above as follows:
\begin{equation}
    \begin{cases}
        \sum\limits_{(i,j)\in \delta^{+}(\Lambda)}\beta_{ij,r}\lambda_{r}\le 1, \forall \Lambda \in \Gamma\\
        \sum\limits_{(i,j)\in \delta^{+}(\Lambda)}\beta_{ij,r}\lambda_{r}\ge 2, \forall \Lambda \in \Gamma
    \end{cases}
\end{equation}
where $\Gamma$ denotes the set of selected node pairs. Note that the inequalities above can be directly added to the RMP without any structural changes. Without loss of generality, let $\gamma_{\Lambda}$ be the dual variables associated with the $\le 1$ branch. The updated reduced cost function $\Tilde{c}_{r}'$ takes the form:
\begin{equation}
    \Tilde{c}_{r}'\leftarrow \Tilde{c}_{r} - \sum_{\Lambda\in\Gamma} \gamma_{\Lambda} \sum_{(i,j)\in\delta^{+}(\Lambda)}\beta_{ij,r}\label{eq:branching_reduced_cost}
\end{equation}}
\subsection{Handling Dual Variables in Valid Inequalities}\label{appendix:dual_valid_inequalities}

\textbf{IPEC:} We first introduce the way to handle the dual variables in the IPEC. For the inequality $\sum_{r\in\Omega'}\sum_{(i,j)\in r} \beta_{ij,I} \lambda_{r}\le |I| - 1$, we denote by $\{\phi_I\le 0:I\in\mathcal{I}\}$ the associated dual variables. And the reduced cost function is updated as follows.
\begin{equation}
    \Tilde{c}_{r}'\leftarrow \Tilde{c}_{r} - \sum_{I\in\mathcal{I}} \phi_{I} \sum_{(i,j)\in I}\beta_{ij,r}
\end{equation}

\textbf{Two-path inequality and rounded capacity:} we note that the reduced cost regarding the two-path and rounded capacity inequalities takes similar forms in the Eq.~\eqref{eq:branching_reduced_cost} with respect to the set $\mathcal{S}$. Specifically, let $\{\gamma_{\mathcal{S}}:\mathcal{S}\in \BFS\}$ be the dual variables in Eqs.~\eqref{eq:2path_inequality} or \eqref{eq:RC_inequality}. The updated reduced cost function is expressed as follows.
\begin{equation}
    \Tilde{c}_{r}'\leftarrow \Tilde{c}_{r} -  \sum_{\mathcal{S}\in \BFS} \gamma_{\mathcal{S}} \sum_{i\in \mathcal{S}} \sum_{j\in \mathcal{N} \setminus \mathcal{S}} \beta_{ij,r}
\end{equation}

\section{Impact of Valid Inequalities}\label{appendix:impact_of_valid_inequality}

\begin{table}[H]
\TABLE
{Impacts of Cuts on the Extended DARP Benchmark Instances 
\label{table:impact_of_valid_ineuqalities}}
{\footnotesize
\begin{tabular}{ccccccccccccccccc}\toprule
    \multirow{2}{*}{Instance} & \multirow{2}{*}{UB} & \multicolumn{5}{c}{LB} & \multicolumn{3}{c}{\# of Cuts}& \multicolumn{5}{c}{CPU Time (s)} \\\cmidrule(lr){3-7} \cmidrule(lr){8-10} \cmidrule(lr){11-15} 
    && No Cuts & IPEC & 2PC & RC & All& IPEC & 2PC & RC & No Cuts & IPEC & 2PC & RC & All  \\ \hline
a2-16 &  294.25 &  100.00 &    * &      * &    * &      * &       198 &    250 &                55 &       0.2 &     0.2 &   0.2 &   0.2 &    0.3 \\
a2-20 &  344.83 &  100.00 &    * &      * &    * &      * &       295 &    230 &                43 &       1.0 &     1.0 &   1.0 &   1.3 &    1.1 \\
a2-24 &  441.06 &   99.78 &  0.0 &    0.0 &  0.0 &    0.0 &       443 &      0 &                55 &       1.9 &     2.0 &   1.9 &   2.4 &    2.4 \\
a3-24 &  344.83 &  100.00 &    * &      * &    * &      * &       498 &    122 &                36 &       1.3 &     1.4 &   1.3 &   1.3 &    1.5 \\
a3-30 &  494.84 &  100.00 &    * &      * &    * &      * &       643 &    348 &                97 &       3.2 &     3.7 &   4.0 &   2.8 &    2.9 \\
a3-36 &  595.86 &   98.79 &  0.0 &  31.81 &  0.0 &  31.81 &       999 &    118 &                45 &       6.2 &     9.4 &  10.5 &   7.0 &   10.8 \\
a4-32 &  487.43 &  100.00 &    * &      * &    * &      * &       763 &    223 &                54 &       1.7 &     1.8 &   1.9 &   1.7 &    1.7 \\
a4-40 &  557.68 &  100.00 &    * &      * &    * &      * &      1325 &    330 &                69 &       8.0 &     8.4 &  10.5 &   8.8 &    9.0 \\
a4-48 &  678.58 &   99.22 &  0.0 &    0.0 &  0.0 &    0.0 &      1953 &    216 &                34 &      26.8 &    25.0 &  26.6 &  23.3 &   28.6 \\
a5-40 &  499.16 &  100.00 &    * &      * &    * &      * &      1337 &    206 &                48 &       5.4 &     6.9 &   6.8 &   7.2 &    8.4 \\
a5-50 &  686.62 &   99.62 &  0.0 &    0.0 &  0.0 &    0.0 &      2142 &    272 &                31 &      17.2 &    18.6 &  15.5 &  18.3 &   16.0 \\
a5-60 &  817.44 &   99.73 &  0.0 &    0.0 &  0.0 &    0.0 &      2963 &    295 &                43 &      58.8 &    66.2 &  62.0 &  61.1 &   81.1 \\
a6-48 &  607.94 &  100.00 &    * &      * &    * &      * &      2007 &    228 &                35 &      13.5 &    15.4 &  15.2 &  12.1 &   16.7 \\
a6-60 &  833.90 &   99.99 &  0.0 &    0.0 &  0.0 &    0.0 &      3007 &    404 &                33 &      44.9 &    40.1 &  42.1 &  41.8 &   41.9 \\
a6-72 &  922.97 &   99.89 &  0.0 &    0.0 &  0.0 &    0.0 &      4401 &    287 &                41 &     212.1 &   201.7 & 196.3 & 201.1 &  207.6 \\
a7-56 &  725.23 &  100.00 &    * &      * &    * &      * &      2576 &    419 &               132 &      20.9 &    21.9 &  24.0 &  20.5 &   23.7 \\
a7-70 &      NaN &     NaN &    - &      - &    - &      - &      4046 &    359 &                65 &     143.6 &   139.2 & 147.9 & 139.6 &  154.8 \\
a7-84 & 1042.00 &   99.35 &  0.0 &   4.01 &  0.0 &   4.01 &      6069 &    404 &                44 &     248.5 &   254.0 & 270.9 & 285.0 &  299.2 \\
a8-64 &  756.35 &   99.93 &  0.0 &    0.0 &  0.0 &    0.0 &      3581 &    379 &                51 &      43.3 &    48.7 &  45.5 &  47.1 &   50.3 \\
a8-80 &     NaN &     NaN &    - &      - &    - &      - &      5578 &    179 &                18 &     190.8 &   179.0 & 199.1 & 186.8 &  197.1 \\
a8-96 &     NaN &     NaN &    - &      - &    - &      - &      - &    - &                - &     - &   - & - & - &  - \\\bottomrule
\end{tabular}}
{-: CPU time limit;
*: Optimality at the root node.}
\end{table}

\section{Assessment of risk score}\label{appendix:risk_assessment}

Specifically, the daily-updated county-level COVID risk level can be accessed from COVID Act Now (\url{https://covidactnow.org}), an online risk/vaccine level tracker, where levels 2 and 3 are converted to the risk scores of 0.2 and 0.3. We further manually quantify the risk levels of the passengers based on the travel purposes and age groups. Based on these information, the calculation of the individual risk level for passenger $i$ can be expressed as below:

\begin{equation}
    r_{i} = r_{i}^{1} + r_{i}^{2} + r_{i}^{3}
\end{equation}
where $r_{i}^{1}$, $r_{i}^{2}$, and $r_{i}^{3}$ represents the risk level of passenger $i$ regarding the factors of travel purpose, age group, and origin county, respectively. 
The detailed specification of individual risk level is summarized as below.

\begin{enumerate}
    \item $r^1_i$, Travel purpose:
    \begin{enumerate}
        \item 0.1: `Personal'
        \item 0.2: `Education', `Employment', `Workshop', `Recreation', `Shopping'
        \item 0.3: `Dialysis', `Medical', `Nutrition' (visiting senior center)
    \end{enumerate}
    \item $r^2_i$, Age group:
    \begin{enumerate}
        \item 0.1: 18 - 44
        \item 0.2: 44 - 65
        \item 0.3: 65 and over.
    \end{enumerate}
    \item $r^3_i$, County-level risk:
    \begin{enumerate}
        \item 0.2: Pick-up locations at Walker County
        \item 0.3: Pick-up locations at Jefferson and Shelby Counties.
    \end{enumerate}
\end{enumerate}

We recognize that our risk score estimation has limitations in representing the most accurate information for each passenger. However, it is important to note that passengers can voluntarily provide such information to the DRT agency. Our risk estimation is used to differentiate groups of passengers, with higher risk passengers being separated from others. For example, two seniors traveling for dialysis may have higher risk scores than three teenagers traveling to an activity center, and our goal is to split those seniors into different routes while considering the admissible exposed risk among the three teenagers.

Even in the absence of precise risk scores, our model can still support passengers with the same level of risk. Our RDARP aims to minimize the contact duration between homogeneous passengers, which provides meaningful risk-aware operation strategies. Moreover, this study focuses on a general model for a risk-aware operation that can be easily adopted by different demand-responsive transit operators with basic reservation information. Therefore, the accuracy of risk assessment or the legal applicability is outside the scope of our research.

\section{Computational Results}
\subsection{Sensitivity analyses on time window}

\begin{table}[H]
\TABLE
{Computational results of Pareto fronts for the real-world instances 
\label{table:pareto_front_real-world_TW_V1}}
{\footnotesize
\begin{tabular}{rrr|rr|rr|rr|rr|rr}\toprule
    \multirow{2}{*}{Instance} & \multirow{2}{*}{$b_i-a_i$} &  \multirow{2}{*}{\# Sol} & \multicolumn{2}{c}{Historical} &\multicolumn{2}{c}{DARP sol} & \multicolumn{2}{c}{Min-Risk RDARP sol} & \multicolumn{2}{c}{DARP vs. RDARP} & \multirow{2}{*}{$T_{MP}$} & \multirow{2}{*}{$T_{SP}$} \\\cmidrule(lr){4-5}\cmidrule(lr){6-7}\cmidrule(lr){8-9}\cmidrule(lr){10-11}
    &&& Cost & $\overline{H}$ & Cost & $\overline{H}$ & Cost & $\overline{H}$ & Cost (\%) & $\overline{H}$(\%) & & \\ \hline
  AM-25-5 &  30 &     9 &            1093.53 &        138.10 & 747.60 & 42.50 &   1065.78 &      0.00 &     42.56 &    100.00 & 1256.9 &  1930.1 \\
  AM-25-5 &  45 &     9 &            1093.53 &        138.10 & 739.84 & 42.50 &   1043.94 &      0.00 &     41.10 &    100.00 &  341.2 &  1361.9 \\
  AM-25-5 &  60 &    11 &            1093.53 &        138.10 & 713.05 & 42.50 &    951.75 &      0.00 &     33.48 &    100.00 &   41.2 &   922.5 \\
  AM-32-7 &  30 &    13 &            1314.34 &        138.10 & 885.23 & 60.46 &   1185.55 &      0.00 &     33.93 &    100.00 & 7415.8 & 11348.3 \\
  AM-32-7 &  45 &    11 &            1314.34 &        138.10 & 863.72 & 51.46 &   1158.64 &      0.00 &     34.15 &    100.00 & 7275.0 &  9749.7 \\
  AM-32-7 &  60 &    13 &            1314.34 &        138.10 & 850.62 & 48.87 &   1151.56 &      0.00 &     35.38 &    100.00 &   25.6 &  7682.4 \\
  PM-25-5 &  30 &     8 &            1151.16 &         26.18 & 819.88 & 39.75 &    953.15 &      0.00 &     16.25 &    100.00 &   43.1 &   350.1 \\
  PM-25-5 &  45 &    10 &            1151.16 &         26.18 & 803.84 & 50.77 &    946.12 &      0.00 &     17.70 &    100.00 &   26.0 &   551.5 \\
  PM-25-5 &  60 &     8 &            1151.16 &         26.18 & 770.34 & 45.49 &    912.02 &      2.18 &     18.39 &     95.21 &    1.7 &   401.1 \\
  PM-32-7 &  30 &     9 &            1471.93 &         49.23 & 985.84 & 38.18 &   1179.01 &      0.00 &     19.59 &    100.00 & 7320.0 &  8877.7 \\
  PM-32-7 &  45 &     8 &            1471.93 &         49.23 & 978.45 & 38.18 &   1152.62 &      0.00 &     17.80 &    100.00 & 7363.2 & 11579.5 \\
  PM-32-7 &  60 &    10 &            1471.93 &         49.23 & 966.92 & 29.75 &   1123.47 &      0.00 &     16.19 &    100.00 &    5.0 &  6054.9 \\
  \bottomrule
\end{tabular}}{}
\end{table}

\end{APPENDICES}

\end{document}